\documentclass[12pt]{article}
\usepackage{stmaryrd,amssymb,amsmath}

\begin{document}

\newtheorem{theorem}{Theorem}[section]
\newtheorem{lemma}[theorem]{Lemma}
\newtheorem{definition}[theorem]{Definition}
\newtheorem{example}[theorem]{Example}
\newtheorem{corollary}[theorem]{Corollary}

\title{
Coarse-scale representations and smoothed Wigner transforms
}

\author{
Agissilaos G. ATHANASSOULIS%
    \footnote{DMA, ENS, 45 rue d'Ulm,
F-75005 Paris and INRIA Paris-Rocquencourt, EPI Bang, France, agis.athanassoulis@gmail.com},
Norbert J.\ MAUSER%
\footnote{
Wolfgang Pauli Inst. c/o  Fak.\ f.\ Mat., Univ. Wien, Nordbergst.\ 15,
    A-1090 Wien, mauser@courant.nyu.edu},
Thierry PAUL%
    \footnote{CNRS, DMA, ENS, 45 rue d'Ulm,
F-75005 Paris, paul@dma.ens.fr}
}
\date{}
\maketitle
\begin{abstract}
Smoothed Wigner transforms   have been used in signal processing, as a regularized version of the Wigner transform, and have been proposed as an alternative to it in the homogenization and / or semiclassical limits of wave equations.

We derive explicit, closed formulations for the coarse-scale representation of the action of pseudodifferential operators.
The resulting ``smoothed operators'' are in general of infinite order.
 The formulation of an appropriate framework, resembling the Gelfand-Shilov spaces, is necessary.

Similarly we treat the ``smoothed Wigner calculus''. In particular this allows us to reformulate any linear equation, as well 
as certain nonlinear ones (e.g. Hartree and cubic non-linear Schr\"odinger), as coarse-scale phase-space equations (e.g. smoothed Vlasov), 
with spatial and spectral resolutions controlled by two free parameters. 
Finally, it is seen that the smoothed Wigner calculus can be approximated,  uniformly on phase-space, by differential operators in the semiclassical regime. This improves the respective weak-topology approximation result for the Wigner calculus.

\end{abstract}
\vskip 1.7cm

\tableofcontents


\section{Introduction and statement of the main results}

Homogenization is an increasingly important and diverse paradigm of applied mathematics.  It is fair to say that it often consists in the reduction of a ``complicated'' problem, typically involving multiple scales, to an effective problem which describes correctly certain coarse-scale features, while appropriately averaging ``less important'' ones, without having to keep track of them explicitly. The simpler, effective problem is then more amenable to analytical and / or numerical treatment.

The Wigner transform (WT) has been used extensively in the homogenization of wave problems, and notably in semiclassical  asymptotics for the Schr\"odinger equation. As we briefly mention in the abstract, this is where the motivation for the smoothed Wigner transform (SWT) and for the study of ``smoothed operators'' comes from. However, the development of the ``smoothed calculus'' is sophisticated enough on its own, and completely independent from the specifics of the WT or the SWT on the technical level. 

Because of that, and in order to make the presentation friendlier to the reader, the paper and its introduction have two parts. First we discuss the smoothed calculus: in Section \ref{Intro1} we motivate the result, and outline its simplest formulation. 
In Section \ref{sec1_2} we briefly outline the application of the smoothed calculus to a SWT-based scheme for the reformulation of general classes of wave equations to coarse-scale kinetic problems in phase-space, i.e. in a space of position $x$ and momentum / wavenumber $k$. 


\subsection{Deriving smoothed equations}
\label{Intro1}

On the technical level, the derivation of equations for the SWT consists in ``smoothing'' the well known Wigner equations. That is, starting from a well-defined system of 
equations\footnote{written here in a symbolic form, assuming that the operator $L$ contains all initial/boundary conditions etc.}
\begin{equation}
\label{eqA}
Lw=0 
\end{equation}
for an appropriate function $w$, we want to derive the equations governing a {\em smoothed version of $w$}, symbolically $\Phi w$. In other words -- and more generally -- we want to smooth equation (\ref{eqA}) and commute correctly the smoothing with the operator,
\begin{equation}
\label{eqB}
Lw=0  \,\, \Leftrightarrow \,\, \Phi L w=0 \,\, \Leftrightarrow \,\, \tilde{L} (\Phi w)=0,
\end{equation}
so as to get a closed problem for the smoothed function $\tilde{w}=\Phi w$. As we will see in the sequel, building a useful, practical theory using the ``smoothed" problem
\begin{equation}
\label{eqCc}
\tilde{L} \tilde{w}=0,
\end{equation}
involves  more work than just deriving it. 

In the case of smoothed Wigner transforms, which is our concrete motivation, it must be emphasized that there are sound mathematical and physical reasons to believe that the smoothed equations are useful (at least for certain problems). In addition, computational aspects (in particular the treatment of concrete problems, with comparisons to exact and / or independent full numerical solutions) of smoothed Wigner equations have already been examined in \cite{Ath,Ath0}, with very encouraging first results. 

The smoothed dynamics can be expressed as ``convolution-deconvolution sandwiches'' $\tilde{L}=\Phi L \Phi^{-1}$, 
\begin{equation}
\label{eq0}
\Phi L \Phi^{-1} : \Phi w \mapsto \Phi L w.
\end{equation}
So the first question we treat in this paper is the explicit computation / representation of these ``convolution-deconvolution sandwiches'', in ways amenable to analysis and computation. This comes together with the need for a basic framework, since the new operators are of infinite order in general.

The core result can be outlined as follows (see Section \ref{AppA} for the notations and
conventions we use for the Fourier transform and the Weyl pseudodifferential
calculus): 
\begin{theorem}[Smoothed calculus]
\label{thrm1}
Let $f(x) \in \mathcal{S}(\mathbb{R}^n)$, $L(x,k) \in \mathcal{S}'(\mathbb{R}^{2n})$ and $L$ be the operator with total Weyl symbol $L(x,k)$, i.e. $Lf(x)=\int_{y,k\in\mathbb{R}^n} e^{2\pi ik(x-y)} L(\frac{x+y}2,k) f(y)dkdy$ . Denote moreover by $\Phi$ the smoothing operator
\begin{equation}
\label{Phi}
\Phi: f(x)\mapsto \mathcal{F}^{-1}_{X \shortrightarrow x} \left[{e^{-\frac{\pi}{2}\sigma^2 X^2 }  \mathcal{F}_{a \shortrightarrow X} \left[{ f(a) }\right] }\right],
\end{equation}
(where $\mathcal{F}$  is the Fourier transform) and $w=\Phi f$. Then the operator $\tilde{L}=\Phi L \Phi^{-1}$ can be expressed as
\begin{equation}
\label{eqC}
\tilde{L}w(x)=2^n \int\limits_{k,u \in \mathbb{R}^n} {
e^{ 2\pi  (ix-\sigma^2 k) u+ 2 \pi i k x } \hat{w}(k-u) \hat{L}_1 (2u,k)dudk},
\end{equation}
where
\begin{equation}
\label{eqCa}
\hat{L}_1 (u,k)=\mathcal{F}_{x \shortrightarrow u} \left[{ L(x,k) }\right].
\end{equation}
Another, equivalent formulation is
\begin{eqnarray}
\label{eqCb}
\tilde{L}  w (x) =  \int\limits_{y,k \in \mathbb{R}^n} { e^{2\pi i (x-y)k -2\pi \sigma^2 k^2} L \left( \frac{x+y}2,k \right) w(y-i\sigma^2 k) dydk}.
\end{eqnarray}
\end{theorem}

Naturally, before we can prove Theorem \ref{thrm1} we have to show that equations (\ref{eqC}) and (\ref{eqCb}) make sense. This is achieved by showing that smoothed functions $w(x)$ have properties very closely resembling those of Gelfand-Shilov functions of type $\mathcal{S}^{\frac{1}{2},B}$ \cite{Gel}, as we see in more detail in the body of the paper (section \ref{sec2_3}). \\ \vspace{0.5mm}

A natural question to ask is ``what is the Weyl symbol of a convolution-deconvolution sandwich $\Phi L \Phi^{-1}$?'' To motivate the answer, let us first consider the free Schr\"odinger evolution
 \begin{equation}\label{sch}
 \psi_t(x):=(e^{\frac{i t}{2}\Delta}\psi_0)(x)=\frac1 {\sqrt{2\pi i t}}\int e^{\frac{-(x-y)^2}{2it}}\psi_0(y)dy.
 \end{equation}
It is well known that the Wigner transform of a wavefunction $\psi$ is the Weyl symbol of its orthogonal projector,
\[
W[\psi](x,k)= \sigma_{Weyl}\left( |\psi\rangle \langle \psi | \right).
\]
It is also well known that the free-space Schr\"odinger evolution (\ref{sch}) is pushed on the Wigner function level by 
\begin{equation}
\label{frew}
W[\psi_t](x,k)=W[\psi_0](x+2\pi tk,k),
\end{equation}
which leads to
\begin{equation}
\label{frew2}
W[\psi_0](x+2\pi tk,k)=W[\psi_t](x,k)=e^{\frac{i t}{2}\Delta} \, W[\psi_0](x,k) \, e^{-\frac{i t}{2}\Delta}.
\end{equation}
In fact this extends immediately to any Weyl operator composed by free Schr\"odinger evolution, i.e. in general
\begin{equation}
\sigma_{Weyl}\left(e^{\frac{i t}{2}\Delta} \, L \, e^{-\frac{i t}{2}\Delta} \right)(x,k)=\sigma_{Weyl}(L)(x-2\pi tk,k).
\end{equation}

By noting that $\Phi=e^{\frac{\sigma^2}{4\pi}\Delta}$ which is {\it formally}\rm\  equal to $e^{i\frac{\sigma^2}{4\pi i} \Delta}$ we can expect that, in the case of analytic symbol, we should have
\begin{equation}
\label{eqD}
\sigma_{Weyl}\left( \Phi L \Phi^{-1} \right)(x,k) =\sigma_{Weyl}\left( L \right) (x+ \frac{i \sigma^2 \, k}{2},k).
\end{equation}

This is in fact the case, as we prove in Section \ref{sec2_1}:

\begin{theorem}[Smoothed calculus: Case of analyticity in $x$, Weyl symbol]
\label{thrmwscdd}
Let  $L(x,k) \in \mathcal{S}'(\mathbb{R}^{2n})$.
Moreover denote $\hat{L}_1 (u,k)=\mathcal{F}_{x \shortrightarrow u} \left[{ L(x,k) }\right]$ and assume that $ \exists M>0$ such that $supp \, \hat{L}_1 (u,k) \subseteq [-M,M]^n \times \mathbb{R}^n$, i.e. $\hat{L}_1 (u,k)$ has compact support in $u$. (In particular it follows then that for each $k \in \mathbb{R}^n$, $L(x,k)$ is an entire analytic function of $x$).
Then the Weyl symbol of $\Phi \, L(x,\partial_x) \, \Phi^{-1}$ is
\begin{equation}
\label{eq120a}
\tilde{L} (x,k)=L(x+ \frac{ i \sigma^2 \, k}{2},k).
\end{equation}
\end{theorem}

In the case that the symbol is analytic in $k$ instead of $x$, we can also have a simplified version of Theorem \ref{thrm1}, namely
\begin{theorem}[Smoothed calculus: Case of analyticity in $k$]
\label{thrmwscdif}
Let  $L(x,k) \in \mathcal{S}'(\mathbb{R}^{2n})$. In addition assume that $L(x,k)$ is a continuous function of $(x,k)$, and $\forall{x \in \mathbb{R}^n}$ $L_x(k)=L(x,k)$ is (the restriction to the real numbers of) an entire-analytic function, and moreover $\forall x\in \mathbb{R}^n \,\, G(k,y)=L \left( \frac{x+y}2,k+\frac{i(x-y)}{\sigma^2} \right) \in \mathcal{S}'(\mathbb{R}^n)$. For example differential operators, $L(x,k)=\sum\limits_{m=0}^N A_m(x) k^m$, fall in this category.

Then
\begin{eqnarray}
\label{eq66bINTRO}
\Phi L \Phi^{-1} w (x) =  
\int\limits_{k \in \mathbb{R}^n} { F(x,k) e^{ - 2\pi \sigma^2 k^2} w(x-i\sigma^2 k) dk},
\end{eqnarray}
where
\begin{equation}
\label{eq66b4}
F(x,k)=\int\limits_{y \in \mathbb{R}^n} {
e^{-2\pi i (x-y)k } L(\frac{x+y}2,k+\frac{i(x-y)}{\sigma^2}) dy }.
\end{equation}

\end{theorem}


Finally it is important to emphasize that the explicit computation of the ``sandwich'' does provide us with important partial cancellations: applying $\Phi L \Phi^{-1}$ as three distinct operators passes from a deconvolution, i.e. a Fourier multiplier of growth $e^{C k^2}$. Applying the result of Theorem \ref{thrm1} passes from imaginary translations, i.e. from multipliers of growth $e^{C k}$, but no deconvolutions.

\subsection{Homogenization in terms of the smoothed Wigner transform}
\label{sec1_2}

In this Subsection we briefly outline the idea of SWT-based homogenization and its motivation, and state the main results in that direction.\\ \vspace{0.5mm}

Consider a problem of the form
\begin{equation}
\label{eq1}
\begin{array}{c}
u_t+L(x,\partial_x)u=0, \\
u(x,0)=u_0(x),
\end{array}
\end{equation}
for a wavefunction $u(x,t): \mathbb{R}^{n+1} \shortrightarrow \mathbb{C}^d$, where $L(x,\partial_x)$ is a $d \times d$ matrix of pseudodifferential operators with matrix-valued Weyl symbol $L(x,k)$. We will assume that equation (\ref{eq1}) describes the propagation of waves (e.g. Schr\"{o}dinger, acoustics, Maxwell's equations etc).

The physical observables of the wavefunction are the scalar functions of time
\begin{equation}
\label{eq2}
\mathbb{M}(t)=\int\limits_{y \in \mathbb{R}^n} {
\bar{u}^T(y,t) \, M u(y,t)dy
},
\end{equation}
corresponding to operators $M=M(x,\partial_x)$ from an appropriate class (e.g. with polynomials or Schwartz test functions as their Weyl symbols). Physically, bilinear observables describe e.g. energy and energy flux in many cases, including the examples mentioned earlier\footnote{In most linear problems the natural energy functional is quadratic in the wavefunction; this can be understood better e.g. in the context of variational formulations.}. Auxiliary quantities of interest in this context are the physical- and Fourier-space densities for the observables,
\begin{equation}
\label{eq3}
\begin{array}{c}
\mathbb{M}(x,t)=\bar{u}(x,t)^T Mu(x,t) , \\
\mathbb{M}_\mathcal{F}(k,t)=\hat{\bar{u}}(k,t)^T \widehat{Mu}(k,t) .
\end{array}
\end{equation}
They are called densities because
\begin{equation}
\label{eq4}
\int\limits_{x \in \mathbb{R}^n} {\mathbb{M}(x,t)dx}=\int\limits_{k \in \mathbb{R}^n} {\mathbb{M}_\mathcal{F}(k,t)dk}=\mathbb{M}(t)
\end{equation}

That is, equation (\ref{eq1}) is seen here more as a book-keeping mechanism; the object of interest is not the point values of the wavefunction $u(x,t)$, but a collection of observables and, to some extent, their densities\footnote{Another way to understand this point is that the problem (\ref{eq1}) describes the microscopic dynamics, but we may only be interested in a macroscopic view of the problem.}. This is often a satisfactory framework, most notably in quantum mechanics.

The SWT-based homogenization approach consists in simplifying the book-keeping problem (\ref{eq1}) while keeping track exactly of the observables.
The name hints towards the fact that the SWT results from the well-known Wigner transform, after the latter is convolved with an appropriate kernel. This is a very natural idea, and several variants have appeared in many contexts, most notably time-frequency analysis \cite{Fl1,Jan} and semiclassical limits (see the more detailed discussion after Theorem \ref{thrmwswt}). The SWT is defined as the sesquilinear transform
\begin{equation}
\label{eq5}
\begin{array}{c}
\tilde{W}: f,g \mapsto
\tilde{W}_{i,j}[f,g](x,k)=\\ 
=\left({ \frac{\sqrt{2}}{\sigma_x} }\right)^n
\int\limits_{u,y \in \mathbb{R}^n} {
e^{ -2 \pi i k y - \frac{\pi \sigma_k^2 y^2}{2} -\frac{2 \pi (u-x)^2}{\sigma_x^2} } f_i(u+\frac{y}{2}) \bar{g}_j(u-\frac{y}{2}) du dy }.
\end{array}
\end{equation}
We will work mostly with its quadratic (and time-dependent) version
\[
\tilde{W}[u](x,k,t)=\tilde{W}[u(\cdot,t),u(\cdot,t)](x,k): \mathbb{R}^{2n+1} \shortrightarrow \mathbb{C}^{d \times d}.
\]
Despite the obvious increase in dimensionality introduced by the SWT, it can be used for compression, because it doesn't exhibit oscillations, in contrast to the wavefunction $u$. This is achieved by an appropriate smoothing in phase-space, controlled by the parameters $\sigma_x,\sigma_k$\footnote{Calibrating the smoothing is pretty well understood, but it isn't central here. A rule of thumb is that $\sigma^2_x$ must be comparable to the wavelengths of $u(x,t)$, and $\sigma^2_k$ to the wavelengths of $\hat{u}(k,t)$; see also Section \ref{sec3_1}, and Section \ref{secSemi} for problems in the semiclassical scaling.}. So the basic idea is switching an oscillatory wavefunction for a smooth phase-space density which lives on a twice-dimensional space.

The simplest paradigm for this homogenization approach consists in two steps:
\begin{itemize}
\item
the derivation of exact equations for the evolution in time of $\tilde{W}[u](x,k,t)$, and
\item
the derivation of a ``smoothed trace formula'', expressing the bilinear observables directly in terms of $\tilde{W}[u](x,k,t)$.
\end{itemize}
(One more step is necessary for the treatment of systems, namely decomposing $\tilde{W}[u](x,k)$ on an appropriate matrix basis).

The derivation  of both the equations for the evolution of the SWT and  the smoothed trace formula follows along the same lines as Theorem \ref{thrm1}, and rests on the following core computation:
\begin{theorem}[Smoothed Wigner calculus]
\label{thrm2}
Let $f(x),g(x) \in \mathcal{S}(\mathbb{R}^n)$, $L(x,k) \in \mathcal{S}'(\mathbb{R}^{2n})$ and $L$ be the operator with $L(x,k)$ as its Weyl symbol. Moreover, denote $\, w(x,k)=W[f,g](x,k)$.

Then
\begin{equation}
\label{eq6}
\tilde{W}[Lf,g](x,k)=\tilde{\mathcal{L}} \tilde{W}[f,g](x,k)
\end{equation}
where
\begin{equation}
\label{eq7}
\begin{array}{l}
\tilde{\mathcal{L}} w(x,k)=\\
=2^{2n}
\int\limits_{X,K,S,T} {
e^{2 \pi i \left[{ S(x-K+i \sigma_x^2X) + T(k+X+i\sigma_k^2 K) +xX+kK }\right] }
\hat{L}(2S,2T) \hat{w}(X-S,K-T) dSdTdXdK
}.
\end{array}
\end{equation}
or, equivalently\begin{equation}
\begin{array}{l}
\tilde{\mathcal{L}} w(x,k)= \\
=\int\limits_{S,T \in \mathbb{R}^n} {
\hat{L}(S,T) 
e^{2\pi i (Sx+Tk)-\frac{\pi}{2}\left({ \sigma_x^2 S^2 +\sigma^2_k T^2 }\right)} w(x+\frac{T+i\sigma_x^2 S}{2},k-\frac{S-i\sigma_k^2 T}{2})dSdT}.
\end{array}
\end{equation}
\end{theorem}

Like before, we can compute the Weyl symbol under appropriate assumptions: 
\begin{theorem}[Weyl symbols for the smoothed Wigner calculus]
\label{thrmwswt}
Consider $f(x),g(x) \in \mathcal{S}(\mathbb{R}^n)$, and $L(x,k)$ to be the Fourier transform of a compactly supported tempered distribution. (In particular it follows that it is  the restriction to the real numbers of an entire analytic function). Then the operator $\tilde{\mathcal{L}}$, defined in equation (\ref{eq6}), has Weyl symbol
\begin{equation}
\label{eq8}
\begin{array}{c}
\tilde{\mathcal{L}}(x,k,X,K)
=L\left({
x-\frac{ K-i\sigma_x^2X}{2} , k+\frac{X+i\sigma_k^2K}{2}  }\right)
\end{array}
\end{equation}
in the sense that
\begin{equation}
\label{eq8a}
\begin{array}{c}
\tilde{\mathcal{L}}w(x,k)=
\int\limits_{a,b,X,K \in \mathbb{R}^n} { e^{2\pi i \left[ X(x-a)+K(k-b) \right]} \tilde{\mathcal{L}} \left( \frac{x+a}{2},\frac{k+b}{2},X,K \right) w(a,b)dadbdXdK }.
\end{array}
\end{equation}
\end{theorem}

Observe how simple and intuitive is the passage to phase-space in terms of the Weyl symbols: Theorem \ref{thrmwswt} can be automatically guessed (and proved, if its assumptions hold) from Theorem \ref{thrmwscdd}. \\ \vspace{1mm}

The concept that certain bilinear functionals, and not the point values of the wavefunction, carry the ``important information'' (the ``physical observables''),  originates in quantum mechanics (and has found applications in other contexts as well). Introduced in 1932 \cite{Wig}, the Wigner transform (WT) appeared in the 90's  as an important tool for homogenization of wave propagation. The concept of semiclassical measures was extensively studied; see for example  \cite{LionsPaul,GL93,Ge96,BurqBourbaki,ZZM}; adaptations to the case of Schr\"odinger operators  with periodic coefficients   and applications of the method to vector problems were carried out \cite{Ge91,MMP1,GMMP}; and applications to stochastic problem were also studied, see e.g. \cite{Ryz}.

In many of the works mentioned above the idea of smoothed Wigner transforms (often under the name ``Husimi functions") appears as a technical device, e.g. for proving the positivity of the Wigner measure or for using interpolation estimates from classical kinetic theory, e.g. \cite{MM1}. The link to coherent states (i.e. abstract wavelet transforms) has also been pointed out, e.g. in \cite{LionsPaul}; however ``a theory of smoothed Wigner transforms'' has not been tackled, essentially because of the problem of dealing with ``convolution-deconvolution sandwiches" and formulating explicit smoothed Wigner equations. 

This problem is solved here, and the formulation of exact smoothed Wigner equations for a broad class of problems, together with a smoothed trace formula for the recovery of the observables, is carried out  in Section \ref{secCSPS}. As concrete examples, we work out the smoothed Wigner equations for the linear Schr\"odinger equation, the cubic non-linear Sch\"odinger equation, and the Hartree equation.  \\ \vspace{1mm}

Virtually all the existing work with WTs is in the semiclassical regime, therefore it is appropriate that we look at a semiclassical application. We do so in Section \ref{secSemi}; more specifically, we formulate and prove the following 
\begin{theorem}[Semiclassical finite-order approximations to the smoothed Wigner calculus]
\label{thrm15intro}
Let $N \in \mathbb{N}$, $V(x):\mathbb{R}^n \shortrightarrow \mathbb{R}$.  
Consider a ``semiclassical family of wavefunctions'' $\lbrace f^\varepsilon \rbrace \subset \mathcal{S}(\mathbb{R}^n)$  for which $\exists M_0>0, \varepsilon_0 \in (0,1)$ such that
\begin{equation}
\label{eq149intro}
||f^\varepsilon ||_{L^2 (\mathbb{R}^n) } \leqslant M_0 \,\,\,\, \forall \, \varepsilon \in(0, \varepsilon_0).
\end{equation}

According to Theorem \ref{thrm2},
\begin{equation}
\label{eq141intro}
\tilde{W}^\varepsilon [V f^\varepsilon, f^\varepsilon](x,k)=
\int\limits_{S \in \mathbb{R}^n} 
{ e^{ 2\pi i S x - \frac{\varepsilon \pi}{2} \sigma_x^2 S^2 } \hat{V}(S) \tilde{w}^\varepsilon(x+\frac{i\varepsilon \sigma_x^2}{2}S,k-\frac{\varepsilon}{2}S)dS } .
\end{equation}
Under appropriate assumptions for the potential $V(x)$
\footnote{Of the form 
\[
\exists C_1,M_1>0 \,:\,\, |\hat{V}(k)| \leqslant C_1 |k|^{-M_1} \, \,\, \forall \, 0<|k|\leqslant 1,
\]
and, there is an appropriate $M_2<-3$ (satisfying additional constraints) such that
\[
\exists C_2>0 \,:\,\, |\hat{V}(k)| \leqslant  C_2 |k|^{M_2} \,\,\,\, \forall \, |k|>1.
\]
}, and for $\sigma_x^2 \leqslant 2$, this expression can be approximated by differential operators
\begin{equation}
\label{eq1510intro}
\begin{array}{l}
\tilde{W}^\varepsilon [V f^\varepsilon, f^\varepsilon](x,k)= 
\sum\limits_{m=0}^N \varepsilon^m  P_m \,\, \tilde{W}^\varepsilon[f^\varepsilon](x,k) 
+r_\varepsilon(x,k),
\end{array}
\end{equation}
where $P_m$ are homogeneous differential operators of order $m$, with coefficients depending on the $m$-order derivatives of the smoothed potential $\tilde{V}(x)=
\left( \frac{\sqrt{2}}{\sqrt{\varepsilon} \sigma_x} \right)^n \int\limits_{x' \in \mathbb{R}^n} {e^{-\frac{2\pi |x-x'|^2}{\varepsilon \sigma_x^2}} V(x')dx'}$, and
\begin{equation}
\label{errintro}
|| r_\varepsilon ||_{L^\infty (\mathbb{R}^{2n})} \, = \, O\left( \varepsilon^{\frac{N+1}{2}-n} \right).
\end{equation}
\end{theorem}

It should be pointed out that the -- somehow surprising -- possibly divergent bound of equation (\ref{errintro}) has to be compared with the sharp estimate $|| \tilde{W}^\varepsilon [f^\varepsilon]||_{L^\infty (\mathbb{R}^{2n})}=O\left( \frac{1}{\varepsilon^n} \right)$, as we will see in the case of WKB ansatz; see the first of the remarks after Theorem \ref{thrm15} in  Section \ref{secSemi}. \\ \vspace{0.5mm}

This result can be used to construct PDE approximations to the smoothed Wigner equations corresponding to semiclassical problems, in analogy to the construction of asymptotic equations for the Wigner measure in the works mentioned earlier (e.g. \cite{LionsPaul,GMMP}). Such approximations were first proposed in \cite{Ath,Ath0}, where it was seen that they provide an efficient computational method for semiclassical problems, recovering $\varepsilon$-dependent information that the Wigner measure cannot keep track of. 

%
%
%

\subsection{Proofs}

Theorem \ref{thrm1} is, verbatim, a concatenation of Theorem \ref{thrm3}, equation (\ref{eq13}) in Sections \ref{sec2_1}, and Theorem \ref{thrm11}, equation (\ref{eq66}) in Section \ref{sec2_4}. 

Theorems \ref{thrmwscdd} and \ref{thrmwscdif} are Theorems \ref{thrm3b} and \ref{thrm11diff} in Sections \ref{sec2_1} and \ref{sec2_4}, respectively.

Theorems \ref{thrm2} and \ref{thrmwswt} are Theorems \ref{thrm4} and \ref{thrm4b} respectively, both found in Section \ref{sec3_1}.

Theorem \ref{thrm15intro} is Theorem \ref{thrm15} of Section \ref{secSemi}.

\subsection{Organization of the paper}

In Section \ref{sec2} we develop the necessary prerequisites, prove the core result and look into several equivalent formulations. We look into an elementary formulation, with minimal prerequisites, in Section \ref{sec2_1}. As we mentioned earlier, establishing certain properties of smoothed functions is a necessary step; this is done in Section \ref{sec2_3}. A very helpful tool (although not an exact match) comes from Gelfand-Shilov spaces, briefly reviewed in Section \ref{sec2_2}. Other formulations of the smoothed calculus are formulated and proved in Section \ref{sec2_4}. Indeed, these more sophisticated formulations will prove particularly useful in the sequel. The references to the specifics of Wigner transforms are kept to a minimum throughout these Sections. 

In Section \ref{sec3_1} the fundamental calculus for the SWT is developed, as an application of Section \ref{sec2}, and a general-purpose phase-space reformulation scheme is outlined in Section \ref{secCSPS}. An implementation of the phase-space reformulation of wave problems, as outlined earlier, is presented in Section \ref{secCSPS}. The smoothed Wigner equations for virtually any linear system, and the smoothed trace formula are formulated making use of the smoothed Wigner calculus. An interesting point is that closed phase-space equations for the cubic non-linear Schr\"odinger and Hartree equations can be obtained, essentially with no additional work. In Section \ref{secSemi} we examine the semiclassical asymptotics for the smoothed Wigner calculus. This allows for a quantitative comparison to the respective WT-based results. \\ \vspace{1mm}

\section{Definitions and notations}
\label{AppA}

The Fourier transform is defined as
\begin{equation}
\hat{f}(k)=\mathcal{F}_{x \shortrightarrow k}\left[f(x)\right]
=\int\limits_{x\in \mathbb{R}^n} {e^{-2\pi ikx}f(x)dx}.
\end{equation}
Inversion is given by
\begin{eqnarray}
\check{f}(k)=\mathcal{F}^{-1}_{x \shortrightarrow k}\left[f(x)\right]
=\int\limits_{x\in \mathbb{R}^n} {e^{2\pi ikx}f(x)dx}, \\
\mathcal{F}_{b \shortrightarrow x}\left[ \mathcal{F}^{-1}_{a \shortrightarrow b}\left[ f(a) \right] \right]=
\mathcal{F}^{-1}_{b \shortrightarrow x}\left[ \mathcal{F}_{a \shortrightarrow b}\left[ f(a) \right] \right]=
f(x).
\end{eqnarray}

An operator $L$ is denoted $L(x,\partial_x)$ and said to have Weyl symbol $L(x,k)$ when
\begin{equation}
Lf(x)=\int\limits_{y,k \in \mathbb{R}^n} {
e^{2\pi i k(x-y)} L^(\frac{x+y}{2}, k) f(y)dydk
}.
\end{equation}
An equivalent expression which we also use is
\begin{equation}
\label{eqA03}
Lf(x)= \mathcal{F}^{-1}_{k \shortrightarrow x} \left[{ \,
\int\limits_{u \in \mathbb{R}^n} {
\hat{f}\left({ k-\frac{u}{2} }\right) e^{\pi i x u} \hat{L}_1(u,k) du
} \, }\right],
\end{equation}
where $\hat{L}_1 (u,k)=\mathcal{F}_{x \shortrightarrow u} \left[{ L(x,k) }\right]$.

The trace formula,
\begin{equation}
\int\limits_{x,k \in \mathbb{R}^n} {
L(x,k) W[f,g](x,k)dxdk}=
\int\limits_{y \in \mathbb{R}^n} { Lf(y) \, \bar{g}(y)dy },
\end{equation}
can be seen as an equivalent definition of the Wigner transform (WT),
\begin{equation}
W[f,g](x,k)=\int\limits_{y \mathbb{R}^n} {
e^{-2\pi i k y} f\left({ x+\frac{y}{2}}\right) \bar{g}\left({ x-\frac{y}{2}}\right) dy
},
\end{equation}
 associating it with the Weyl calculus.

Sometimes a scaled version of the Weyl calculus is used. This can be motivated e.g. from WKB functions, i.e. functions of the form $f(x)=A(x) e^{\frac{2\pi i}{\varepsilon} S(x)}$. The scaled Weyl calculus is defined consistently with the scaled WT through the trace formula, i.e. $L(x,\varepsilon \partial_x)$ is defined through
\begin{equation}
\int\limits_{y \in \mathbb{R}^n} { L(x,\varepsilon \partial_x)f(y) \, \bar{g}(y)dy }=
\int\limits_{x,k \in \mathbb{R}^n} {
L(x,k) W^\varepsilon [f,g](x,k)dxdk},
\end{equation}
where
\begin{equation}
W^\varepsilon[f,g](x,k)=\int\limits_{y \mathbb{R}^n} {
e^{-2\pi i k y} f\left({ x+\frac{ \varepsilon y}{2}}\right) \bar{g}\left({ x-\frac{\varepsilon y}{2}}\right) dy
}.
\end{equation}

It is often the case that the Weyl symbol itself depends on $\varepsilon$, $L^\varepsilon(x,k)$ \footnote{e.g. $L^\varepsilon(x,k)=P(x,k)+\varepsilon Q(x,k)$.}. In that case $L^\varepsilon (x,\varepsilon \partial_x)$ is defined as
\begin{equation}
\int\limits_{x,k \in \mathbb{R}^n} {
L^\varepsilon (x,k) W^\varepsilon [f,g](x,k)dxdk}=
\int\limits_{y \in \mathbb{R}^n} { L^\varepsilon (x,\varepsilon \partial_x)f(y) \, \bar{g}(y)dy }.
\end{equation}
This can easily be seen to be equivalent to
\begin{equation}
L^\varepsilon (x,\varepsilon \partial_x) f(a)=\int\limits_{y,k \in \mathbb{R}^n} {
e^{2\pi i k(a-y)} L^\varepsilon (\frac{x+y}{2}, \varepsilon k) f(y)dydk
}.
\end{equation}

For a function $F(x,k)$ satisfying appropriate conditions 
(and for any tempered distribution), 
the operator $F(x,\frac{ \partial_x}{2\pi i})$ is defined, in terms of the Weyl calculus, 
as the operator with Weyl symbol $F(x,k)$. It must be noted that usually Weyl symbol classes are taken to be more restricted than $\mathcal{S}'$. Imposing some more assumptions will probably be necessary in certain contexts; however for our purposes the more general choice works well. 

The order $\delta$ of an operator $L(x,\partial_x)$ is defined as the smallest 
$\delta>0$ s.t. $\forall \, \alpha, \beta \in \mathbb{N} \cup \lbrace 0 \rbrace \, 
\exists C_{\alpha,\beta}>0 \, : \, \forall \, l,m=1,...,n$
\begin{equation}
| \frac{\partial^{\alpha+\beta}}{\partial x_l^\alpha \partial k_m^\beta} L(x,k) | < 
C_{\alpha,\beta} (1+|k|)^{\delta-\beta}.
\end{equation}

Finite order PDOs are well defined on Sobolev spaces of the same order \cite{Hor}. 
Typical examples of infinite order operators include deconvolutions and 
imaginary translations.\\ \vspace{0.5mm}

\noindent  {\bf Remark on notation:} Please note that our conventions and notations for the Fourier transform and the Weyl calculus, clearly stated here, are used throughout the text without additional explanation. \\

\section{Smoothed calculus}
\label{sec2}

\subsection{Explicit formulation of convolution-deconvolution sandwiches}
\label{sec2_1}

In this Subsection we will present the derivation of the elementary formulations 
of Theorem \ref{thrm1} (i.e. equations (\ref{eqC}) ), focusing on the mechanics of the derivation, 
as well as basic interpretation and application issues. 
The pseudodifferential and other operator-theoretic aspects are kept to a minimum here.

\begin{definition}[The smoothing operator]
\label{DefSmoo1}
The operator $\Phi$ is defined as in equation (\ref{Phi}), i.e.
\begin{equation}
\label{PHI}
\begin{array}{c}
\Phi: f(x)\mapsto \mathcal{F}^{-1}_{k \shortrightarrow x} \left[{e^{-\frac{\pi}{2}\sigma^2 k^2 }  \mathcal{F}_{x' \shortrightarrow k} \left[{ f(x') }\right] }\right]
=\frac{2^\frac{n}{2}}{\sigma^n} \int\limits_{x \in \mathbb{R}^n} { e^{-2\pi \frac{(x-x')^2}{\sigma^2}} f(x')dx'}.
\end{array}
\end{equation}
\end{definition}

\noindent {\bf Remark:} {\em
 The notation $e^{-\frac{\pi}{2}\sigma^2 k^2 }$ is used interchangeably with $e^{-\frac{\pi}{2}\sigma^2 |k|^2 }$, i.e. $k^2=k \cdot k$. A natural generalization of Definition \ref{DefSmoo1} would be
\begin{equation}
\label{PHIa}
\begin{array}{c}
\Phi_{\sigma_1,...,\sigma_n}: f(x)\mapsto \mathcal{F}^{-1}_{k \shortrightarrow x} \left[{e^{-\frac{\pi}{2}
\sum\limits_{l=1}^{n} {\sigma_l^2 k_l^2} }  \mathcal{F}_{x' \shortrightarrow k} \left[{ f(x') }\right] }\right]= \\ { } \\
=\frac{2^\frac{n}{2}}{ \mathop{\Pi}\limits_{l=1}^{n} \sigma_l} \int\limits_{x \in \mathbb{R}^n} { e^{-2\pi \sum\limits_{l=1}^{n} {\frac{(x_l-x'_l)^2}{\sigma_l^2}}} f(x')dx'}.
\end{array}
\end{equation}
In fact we will use a a smoothing like that later, but most of the time it isn't worth the notational inconvenience --  and all our results are generalized to the anisotropic case in a straightforward manner.}

Operators like $\Phi$ are very common, and are often called mollifiers. Observe that $\Phi$ is translation invariant, i.e. a Fourier multiplier. Indeed, its action is very intuitively seen in the Fourier domain (it damps ``the high wavenumbers'' with a Gaussian weight). It is also straightforward to observe that it is one-to-one, with inverse
\begin{equation}
\label{PHI1}
\begin{array}{c}
\Phi^{-1}: f(x)\mapsto \mathcal{F}^{-1}_{X \shortrightarrow x} \left[{e^{\frac{\pi}{2}\sigma^2 X^2 }  \mathcal{F}_{x' \shortrightarrow X} \left[{ f(x') }\right] }\right].
\end{array}
\end{equation}
Understanding the image by $\Phi$ of the Schwartz test-functions $\mathcal{S}$, as well as other basic spaces will also be important. In particular, it is fair to say that smoothed functions are restrictions to a real space of entire-analytic functions. 

Observe moreover that  $\Phi^{-1}$ is very hard (often impossible) to implement in practice, e.g. numerically. This is the basic reason why we want to compute explicitly the ``sandwich'' $\Phi L \Phi^{-1}$, looking for some sort of mutual (partial) cancellation of $\Phi$ and $\Phi^{-1}$.

Let us start with a very simple observation:

\begin{lemma}[Smoothed polynomial calculus] $\forall m \in \mathbb{N}, \, i \in \lbrace 1,...,n \rbrace$
\label{lmm1}
\begin{eqnarray}
\label{eq9}
\Phi x_i^m \Phi^{-1} & = & \left({x_i + \frac{\sigma^2 \partial_{x_i}}{4 \pi} }\right)^m, \\
\label{eq10}
\Phi \partial_{x_i}^m \Phi^{-1} & = & \partial_{x_i}^m.
\end{eqnarray}
\end{lemma}

\noindent {\bf Proof:} The way to interpret and prove any expression of the form $\Phi L \Phi^{-1} =\tilde{L}$ is by checking that $\forall f(x) \in \mathcal{S}(\mathbb{R}^n)$
\begin{equation}
\label{eq10a}
\tilde{L}\Phi f= \Phi L f.
\end{equation}
Indeed, to prove equation (\ref{eq9}) for $m=1$ it suffices to check that
\begin{equation}
\begin{array}{c}
\label{eq11}
\left({x + \frac{\sigma^2 \partial_x}{4 \pi} }\right) \Phi f=
x \Phi f + \frac{2^\frac{n}{2}}{\sigma^n} \frac{\sigma^2 \partial_x}{4 \pi} \int\limits_{x \in \mathbb{R}^n} { e^{-2\pi \frac{(x-x')^2}{\sigma^2}} f(x')dx'}=\\
=x \Phi f + \frac{2^\frac{n}{2}}{\sigma^n} \frac{\sigma^2 }{4 \pi} \int\limits_{x \in \mathbb{R}^n} { \partial_x e^{-2\pi \frac{(x-x')^2}{\sigma^2}} f(x')dx'} =\\
=x \Phi f + \frac{2^\frac{n}{2}}{\sigma^n} \frac{\sigma^2 }{4 \pi} \int\limits_{x \in \mathbb{R}^n} {  \frac{-4\pi (x-x')}{\sigma^2}   e^{-2\pi \frac{(x-x')^2}{\sigma^2}} f(x')dx'} =\\
=\frac{2^\frac{n}{2}}{\sigma^n} \int\limits_{x \in \mathbb{R}^n} {  e^{-2\pi \frac{(x-x')^2}{\sigma^2}} x' f(x')dx'} = \Phi (xf).
\end{array}
\end{equation}

In order to prove equation (\ref{eq10}) it is easier to work in the Fourier domain:
\begin{equation}
\begin{array}{c}
\label{eq12}
\mathcal{F}_{x \shortrightarrow X} \left[ \partial_{x_i} \Phi f \right]=
e^{-\frac{\pi}{2} \sigma^2 X^2} 2 \pi X_i \hat{f}(X)=
\mathcal{F}_{x \shortrightarrow X} \left[\Phi  \partial_{x_i} f \right].
\end{array}
\end{equation}

The generalization for $m>1$ for either case is obvious. The proof is complete.

\begin{theorem}[Elementary formulation of the smoothed calculus]
\label{thrm3}
Let $f(x) \in \mathcal{S}(\mathbb{R}^n)$, $L(x,k) \in \mathcal{S}'(\mathbb{R}^{2n})$, $L$ be the operator with Weyl symbol $L(x,k)$ and 
\[
w(x)=\Phi f(x).
\]
Then 
\begin{equation}
\label{eq13}
\Phi L \Phi^{-1} w (x)=2^n \int\limits_{k,u \in \mathbb{R}^n} {
e^{ 2\pi  (ix-\sigma^2 k) u+ 2 \pi i k x } \hat{w}(k-u) \hat{L}_1 (2u,k)dudk
},
\end{equation}
where $\hat{L}_1 (u,k)=\mathcal{F}_{x \shortrightarrow u} \left[{ L(x,k) }\right]$.
\end{theorem}

\noindent {\bf Proof: } Let us check first of all that the integral of equation (\ref{eq13}) indeed converges. To that end observe that the {\em rhs} of equation (\ref{eq13}) is equal to
\begin{equation}
\label{eq14}
\begin{array}{c}
2^n \int\limits_{k,u \in \mathbb{R}^n} {
e^{ 2\pi  (ix-\sigma^2 k) u+ 2 \pi i k x } \hat{w}(k-u) \hat{L}_1 (2u,k)dudk
}=\\
=2^n \int\limits_{k,u \in \mathbb{R}^n} {
e^{ 2\pi  (ix-\sigma^2 k) u+ 2 \pi i k x } e^{-\frac{\pi}{2}\sigma^2(k-u)^2} \hat{f}(k-u) \hat{L}_1 (2u,k)dudk
}=\\
=2^n \int\limits_{k,u \in \mathbb{R}^n} {
e^{ 2\pi  ix(k+u)- \frac{\pi}{2}\sigma^2(k+u)^2 } \hat{f}(k-u) \hat{L}_1 (2u,k)dudk
}=\\
=\int\limits_{k,u \in \mathbb{R}^n} {F_x(k,u)\hat{L}_1 (2u,k)dudk}.
\end{array}
\end{equation}
It suffices to show that
\[
F_x(k,u)=2^n e^{ 2\pi  ix(k+u)- \frac{\pi}{2}\sigma^2(k+u)^2 } \hat{f}(k-u) \, \in \, \mathcal{S}(\mathbb{R}^{2n}),
\]
because if that is true, then the integral exists (for each $x \in \mathbb{R}^n$) as a duality pairing between $F_x(k,u) \, \in \, \mathcal{S}(\mathbb{R}^{2n})$ and $\hat{L}_1 (2u,k) \, \in \, \mathcal{S}'(\mathbb{R}^{2n})$. 

Denote
\begin{equation}
\label{eq15}
\mathcal{T}:F(k,u) \mapsto F(k-u,k+u).
\end{equation}
$\mathcal{T}$ is essentially a rotation, and it is clear that $\mathcal{T}\left({ \mathcal{S}(\mathbb{R}^{2n}) }\right) \subseteq \mathcal{S}(\mathbb{R}^{2n})$. Now observe that
\begin{equation}
\label{eq16}
F_x(k,u)=2^n \mathcal{T}\left({ f(k) e^{-\frac{\pi}{2}\sigma^2 u^2 + 2\pi i xu} }\right) \, \in \, \mathcal{S}(\mathbb{R}^{2n}),
\end{equation}
since $\forall x \in \mathbb{R}^n$
\begin{equation}
\label{eq16aa}
f(k) e^{-\frac{\pi}{2}\sigma^2 u^2 + 2\pi i xu}  \, \in \, \mathcal{S}(\mathbb{R}^{2n})
\end{equation}
and therefore equation (\ref{eq14}) can be cast as a Schwartzian duality pairing (the independent variables in equation (\ref{eq16aa}) are $k,u$; $x$ plays the role of a parameter).

Morally, the point here is that despite the real exponential term in equation (\ref{eq13}), the integral exists because we act on smoothed functions (which have Gaussian decay in the Fourier domain). \\

Let us now prove equation (\ref{eq13}). Starting from the {\em lhs} one observes that (we use equation (\ref{eqA03}) from Section \ref{AppA} for the implementation of the Weyl calculus) 
\begin{equation}
\begin{array}{c}
\label{eq14a}
\Phi L \Phi^{-1} w (x)=\Phi L f(x)=\\
=2^n \Phi \left[{ \mathcal{F}^{-1}_{k \shortrightarrow x} \left[{ \,
\int\limits_{u \in \mathbb{R}^n} {
\hat{f}\left({ k-u }\right) e^{2\pi i x u} \hat{L}_1(2u,k) du
} \, }\right] }\right]=\\
=2^n \mathcal{F}_{l \shortrightarrow x}^{-1} \left[{
e^{-\frac{\pi}{2}\sigma^2 l^2 }  \mathcal{F}_{x' \shortrightarrow l} \left[{
\mathcal{F}^{-1}_{k \shortrightarrow x'} \left[{ \,
\int\limits_{u \in \mathbb{R}^n} {
\hat{f}\left({ k-u }\right) e^{2\pi i x' u} \hat{L}_1(2u,k) du
} \, }\right]
}\right]
}\right] =\\
=2^n \int\limits_{u,k,x',l \in \mathbb{R}^n}{
e^{2\pi i \left[{xl-x'l+kx'+x'u}\right] -\frac{\pi}{2}\sigma^2 l^2 } \hat{f}(k-u) \hat{L}_1(2u,k)dudkdx'dl }=\\
=2^n \int\limits_{u,k,x',l \in \mathbb{R}^n}{
e^{2\pi i x' \left({-l+k+u}\right)}dx' \,\, e^{2\pi i xl -\frac{\pi}{2}\sigma^2 l^2 } \hat{f}(k-u) \hat{L}_1(2u,k)dudkdl }=\\
=2^n \int\limits_{u,k,l \in \mathbb{R}^n}{
\delta(k+u-l) e^{2\pi i xl -\frac{\pi}{2}\sigma^2 l^2 }dl \,\, \hat{f}(k-u) \hat{L}_1(2u,k)dudk }=\\
=2^n \int\limits_{u,k, \in \mathbb{R}^n}{
e^{2\pi i x(k+u) -\frac{\pi}{2}\sigma^2 (k+u)^2 } \hat{f}(k-u) \hat{L}_1(2u,k)dudk }.
\end{array}
\end{equation}
But we have already seen the last expression above to be equal to the {\em rhs} of equation (\ref{eq13}), in equation (\ref{eq14}). The proof is complete. \\

It is very easy to check the consistency of Theorem \ref{thrm3} with the Weyl calculus for $\sigma=0$, and with Lemma \ref{lmm1}. Indeed the comparison with Lemma \ref{lmm1} reveals nicely the ``miraculous cancellation'' that takes place for polynomials: if $L(x,k)$ is a polynomial, then $\mathop{supp}_{u} \hat{L}_1(u,k) = \lbrace 0 \rbrace $, and it doesn't allow the real exponential term to give rise to an infinite order operator. For example, if $L(x,k)=(2\pi i x)^m$, equation (\ref{eq13}) becomes
\begin{equation}
\label{eq17}
\Phi L \Phi^{-1} w (x)= \int\limits_{k,u \in \mathbb{R}^n} {
e^{ \pi  (ix-\sigma^2 k) u} \hat{w}\left({k-\frac{u}{2}}\right) \delta^{(m)}(u-0) du \,\, e^{2 \pi i k x} dk
},
\end{equation}
generating of course the same end result as Lemma \ref{lmm1}. In general however contributions from $u$ away from zero will give rise to operators not of finite order (e.g. for $L(x,k)$ a Gaussian). 

Observe finally that our result is, under appropriate conditions,  equivalent to
\begin{equation}
\label{eq18}
\Phi L \Phi^{-1}=L(x+ \frac{ \sigma^2 \, \partial_x}{4 \pi},\frac{\partial_x}{2\pi i}).
\end{equation}

\begin{theorem}[Weyl symbols for the smoothed calculus]
\label{thrm3b}
Let  $L(x,k) \in \mathcal{S}'(\mathbb{R}^{2n})$.
Moreover denote $\hat{L}_1 (u,k)=\mathcal{F}_{x \shortrightarrow u} \left[{ L(x,k) }\right]$, and assume that $ \exists M>0$ such that $supp \, \hat{L}_1 (u,k) \subseteq [-M,M]^n \times \mathbb{R}^n$, i.e. $\hat{L}_1 (u,k)$ has compact support in $u$. (In particular it follows then that for each $k \in \mathbb{R}^n$, $L(x,k)$ is an entire analytic function of $x$).
Then the Weyl symbol of $\Phi \, L(x,\partial_x) \, \Phi^{-1}$ is
\begin{equation}
\label{eq120}
\tilde{L} (x,k)=L(x+ \frac{ i \sigma^2 \, k}{2},k).
\end{equation}
\end{theorem}

\noindent {\bf Remark:} {\em Observe that the Weyl symbols produced this way need not be tempered distributions (take e.g. $L(x,k)=e^{ix}$). This is related to the fact that the smoothed operators are of infinite order in general. Understanding better what this means, both in terms of analysis and applications, is one of the objectives of this paper. In any case it hints towards the use of more general test-function / distribution theories as a natural direction. We believe that in that context the assumptions on $L(x,k)$ for this Theorem can be probably relaxed. \\}

\noindent {\bf Proof: }  We will directly verify that the operator with the Weyl symbol of equation (\ref{eq120}) (equivalently the operator of equation (\ref{eq18})) is the same as that of equation (\ref{eq13}):
\begin{equation}
\label{eq19}
\begin{array}{c}
\tilde{L} w(x)=
\int\limits_{y,k \in \mathbb{R}^n} {
e^{2\pi i k (x-y)}L\left({ \frac{x+y}{2} + \frac{i \sigma^2 k}{2},k }\right) w(y)dydk}=\\
=2^n \int\limits_{y,k \in \mathbb{R}^n} {
e^{2\pi i k (x-y)} \mathcal{F}^{-1}_{u \shortrightarrow x} \left[{ e^{2\pi i u (x+y+i \sigma^2 k )} \hat{L}_1(2u,k) }\right] w(y)dydk}=\\
=2^n \int\limits_{y,k \in \mathbb{R}^n} {
e^{2\pi i \left[{ k (x-y) +u (x+y+i \sigma^2 k ) }\right] } w(y) dy \,\,   \hat{L}_1(2u,k)  dudk}=\\
=2^n \int\limits_{y,k \in \mathbb{R}^n} {
e^{-2\pi i y(k-u)  } w(y) dy \,\,e^{2\pi i \left[{ kx +u (x+i \sigma^2 k ) }\right] }   \hat{L}_1(2u,k)  dudk}=\\
=2^n \int\limits_{y,k \in \mathbb{R}^n} {
\hat{w}(k-u)  \,\,e^{2\pi i \left[{ kx +u (x+i \sigma^2 k ) }\right] }   \hat{L}_1(2u,k)  dudk}.
\end{array}
\end{equation}
The first step which needs more explanation is the equality
\begin{equation}
\label{eq121}
\begin{array}{c}
2^{-n} L\left({ \frac{x+y}{2} + \frac{i \sigma^2 k}{2},k }\right) =
\mathcal{F}^{-1}_{u \shortrightarrow x} \left[{ e^{2\pi i u (x+y+i \sigma^2 k )} \hat{L}_1(2u,k) }\right].
\end{array}
\end{equation}
To see that, first of all recall that $L(x,k) \in \mathcal{S}'(\mathbb{R}^{2n})$, and therefore $\hat{L}_1(2u,k)=\mathcal{F}_{x \shortrightarrow U} \left[{L(x,k)}\right] |_{U=2u}$ is well defined. Now
\begin{equation}
\label{eq122}
\begin{array}{c}
\int\limits_{k \in \mathbb{R}^n} { e^{2\pi i u \, \cdot \, (x+y+i \sigma^2 k)} \hat{L}_1(2u,k) du}=
\int\limits_{k \in \mathbb{R}^n} { \sum\limits_{l=0}^{\infty} {\frac{ \left[{2 \pi i u \cdot  (x+y+i \sigma^2 k)}\right]^l }{l!}} \hat{L}_1(2u,k)du}=\\ { } \\
=2^{-n} 
\sum\limits_{l=0}^{\infty} {\frac{ \left[ (x+y+i \sigma^2 k)  \, \cdot \,  \partial_x \right] ^l }{ l!} L(0,k) }=
2^{-n} L\left({ \frac{x+y}{2} + \frac{i \sigma^2 k}{2},k }\right).
\end{array}
\end{equation}
The series is a Taylor expansion of $L(x,k)$ in the $x$ variable for each $k$; it converges absolutely following our assumption. This also justifies the interchange of the order of summation and integration in equation (\ref{eq122}) through dominated convergence. (See also lemma \ref{lmm4}).

Another step in equation (\ref{eq19}) that needs justification is the interchange of the $du$ and $dy$ integrations, passing from the second to the third line. There it suffices to observe that $e^{-2 \pi \sigma^2 u k}$ can be replaced by $e^{-2 \pi \sigma^2 u k} \, \chi_{[-2M,2M]^n}(u)$, in which case the result follows by the standard tempered distribution calculus.

The proof is complete.  \\ \vspace{0.5mm}

As we commented briefly earlier, it seems reasonable that $L(x+ \frac{ i \sigma^2 \, k}{2},k)$ can be defined precisely, in an appropriate (weak) sense for any $L(x,k) \in \mathcal{S}'(\mathbb{R}^{2n})$. Gelfand-Shilov spaces of test functions and their duals (ultra-distributions) seem like an appropriate framework in that direction. 

Before we get to that question, however, there are some more basic results about smoothed functions that we need; these too can be easily formulated as modifications of standard results for Gelfand-Shilov test-functions. In the following Section we go briefly over the basic facts concerning Gelfand-Shilov test-functions, to set the stage for our results concerning additional formulations of the smoothed calculus.

\subsection{A natural framework for the smoothed calculus: Gelfand-Shilov spaces}
\label{sec2_2}

Gelfand-Shilov spaces (and ultra-distributions) is a well-defined topic which has been attracting increasing attention recently, with several dedicated monographs and papers. It is clear that a full presentation of them is completely outside the scope of this work. What is however necessary, is to outline some basic facts and also some motivation as to why these -- somewhat unusual -- spaces are well suited for the study of our smoothed calculus. Here we will focus on presenting the necessary background; in the next Subsection we will focus on its use and application. \\ \vspace{0.5mm}

Usual test-functions, i.e. Schwartz test functions, are defined by
\begin{equation}
\label{eq29}
f(x) \in \mathcal{S}(\mathbb{R}^n) \,\,\Leftrightarrow \,\, \mathop{sup}\limits_{x \in \mathbb{R}^n} | x^p \partial_x^q f(x) |< \infty \,\,\, \forall p,q \in \left( \mathbb{N} \cup \lbrace 0 \rbrace \right)^n,
\end{equation}
i.e., morally, by decay of any order of derivatives faster than any power of $\frac{1}{|x|+1}$  \cite{Gel}. It is important to note that this condition is symmetric with respect to the Fourier transform, i.e.
\begin{equation}
\label{eq30}
\mathcal{F}\left({ \mathcal{S}(\mathbb{R}^n) }\right)= \mathcal{S}(\mathbb{R}^n) .
\end{equation}
 This allows for a theory of distributions ``of some polynomial order of growth'' in space and Fourier domain.
In addition, the Weyl calculus allows the transfer of function-space and function-theory results to operators, providing very strong tools for many problems -- in particular problems involving differential operators. However, there are two classes of operators, very important in this study, that are not contained in this framework: imaginary translations (i.e. extensions to the complex plane)
\begin{equation}
\label{eq32}
T_{i\sigma}:f(x) \mapsto f(x+i\sigma)=\mathcal{F}^{-1}_{k \shortrightarrow x} \left[{ e^{2\pi \sigma k}\hat{f}(k) }\right],
\end{equation}
and deconvolutions
\begin{equation}
\label{eq33}
\Phi^{-1}:f(x)\mapsto \mathcal{F}^{-1}_{k \shortrightarrow x} \left[{e^{\frac{\pi}{2}\sigma^2 k^2 }  \hat{f}(k) }\right].
\end{equation}

So the basic motivation is very simple: we need to accommodate Fourier multipliers and / or Weyl symbols of faster than polynomial growth. (Of course things are actually more subtle than that in the end). This brings on the question of smaller test-function spaces -- after all deconvolutions and imaginary translations are {\em not} well defined on all of $\mathcal{S}(\mathbb{R}^n)$. 

These spaces of ``very nice functions'' should then be examined with respect to (nontriviality first of all!), closedness under elementary operations, behaviour under the Fourier transform etc. This is contained in the Gelfand-Shilov theory of test-functions and distributions. \\ \vspace{0.5mm}

At this point we will go over the Definition of Gelfand-Shilov test-functions and some basic facts that are relevant. All results are quoted from chapter IV of \cite{Gel}, and we use the same notation. The implications of these core facts for the problem at hand will be seen shortly. Many facts about them (after some adaptation) will be very helpful for working with smoothed functions.

For simplicity we will work in the case $n=1$, i.e. for functions of $\mathbb{R}$. The generalization to higher-dimensional spaces is straightforward.

\begin{definition}[Spaces of type $\mathcal{S}$]  
\label{def4}
Let $\alpha, \beta \geqslant 0$. \\
\begin{enumerate}
\item
$f(x) \in \mathcal{S}_\alpha$ iff \, $\exists A,C_q$ depending on $f$ s.t.
\begin{equation}
\label{eq34}
| x^k \partial_x^q f(x) | \leqslant C_q A^k k^{k \alpha} \,\,\, \forall k,q\in\mathbb{N} \cup \lbrace 0 \rbrace
\end{equation}

\item
$f(x) \in \mathcal{S}^\beta$ iff \, $\exists B,C_k$ depending on $f$ s.t.
\begin{equation}
\label{eq35}
| x^k \partial_x^q f(x) | \leqslant C_k B^q q^{q \beta} \,\,\, \forall k,q\in\mathbb{N} \cup \lbrace 0 \rbrace
\end{equation}

\item
$f(x) \in \mathcal{S}_\alpha^\beta$ iff \, $\exists A,B,C$ depending on $f$ s.t.
\begin{equation}
\label{eq36}
| x^k \partial_x^q f(x) | \leqslant C A^k B^q k^{k \alpha} q^{q \beta} \,\,\, \forall k,q\in\mathbb{N} \cup \lbrace 0 \rbrace
\end{equation}

\end{enumerate}
\end{definition}

\begin{theorem}[Equivalent Definitions of spaces of type $\mathcal{S}$]
\label{thrm5}
\begin{enumerate}
\item
$f(x) \in \mathcal{S}_\alpha$ iff $\, \forall q\in\mathbb{N} \cup \lbrace 0 \rbrace \,$ $\exists a,C_q$ depending on $f$ s.t.
\begin{equation}
\label{eq37}
| \partial_x^q f(x) | \leqslant C_q e^{ -a|x|^{\frac{1}{\alpha}} } \,.
\end{equation}
The parameter $a$ is related to $A$ of the respective Definition \ref{def4}.1 explicitly, and more specifically $a=a(A)=\frac{\alpha}{e A^{\frac{1}{\alpha}}}$.

\item
$f(x) \in \mathcal{S}^\beta$ iff \, $\forall k\in\mathbb{N} \cup \lbrace 0 \rbrace$ \, $\exists b,C'_k$ depending on $f$ s.t. $f(x)$ can be (uniquely) extended to an entire function $f(x+iy)$ satisfying the estimates
\begin{equation}
\label{eq38}
| x^k  f(x+iy) | \leqslant C'_k  e^{ b |y|^\frac{1}{1-\beta} } \,.
\end{equation}
The parameter $b$ is related to $B$ of the respective Definition \ref{def4}.2 explicitly, and more specifically $b=b(B)=\frac{1-\beta}{e}(Be)^{\frac{1}{1-\beta}}$.

\item
$f(x) \in \mathcal{S}_\alpha^\beta$ iff \, $\exists a,b,C$ depending on $f$ s.t. $f(x)$ can be (uniquely) extended to an entire function $f(x+iy)$ satisfying the estimates
\begin{equation}
\label{eq39}
| f(x+iy) | \leqslant C  e^{ -a|x|^{\frac{1}{\alpha}}+b|y|^{\frac{1}{1-\beta}} }.
\end{equation}
The parameters $a,b$ are explicitly related to the parameters $A,B$ of Definition \ref{def4}.3.
\end{enumerate}
\end{theorem}

\begin{theorem}[Fourier transforms]
\label{thrm6}
\begin{enumerate}
\item
Fourier transforms of $\mathcal{S}_\alpha$
\begin{equation}
\label{eq40}
\mathcal{F} \left({ \mathcal{S}_\alpha }\right)=\mathcal{S}^\alpha.
\end{equation}

\item
Fourier transforms of $\mathcal{S}^\beta$
\begin{equation}
\label{eq41}
\mathcal{F} \left({ \mathcal{S}^\beta }\right)=\mathcal{S}_\beta.
\end{equation}

\item
Fourier transforms of $\mathcal{S}_\alpha^\beta$
\begin{equation}
\label{eq42}
\mathcal{F} \left({ \mathcal{S}_\alpha^\beta }\right)=\mathcal{S}^\alpha_\beta.
\end{equation}

\end{enumerate}
\end{theorem}

\begin{theorem}[Nontriviality]
\label{thrm7}
The space $\mathcal{S}_\alpha$ is nontrivial (i.e. contains a nonzero function) for all $\alpha > 0$. The space $\mathcal{S}^\beta$ is nontrivial for all $\beta >0$. The space $\mathcal{S}_\alpha^\beta$ is nontrivial iff
\begin{itemize}
\item[ ]
$\alpha+\beta \geqslant 1, \,\,\, \alpha>0, \, \beta>0$, or 
\item[ ]
$\alpha=0, \, \beta >1$, or
\item[ ]
$\beta=0, \, \alpha >1$.
\end{itemize}
\end{theorem}

\begin{theorem}[Closedness under elementary operations]
\label{thrm8}
Let $\alpha, \beta \geqslant 0$. Each of the spaces $\mathcal{S}_\alpha^\beta, \mathcal{S}_\alpha, \mathcal{S}^\beta$ is closed under translation
\[
f(x) \mapsto f(x+\lambda), \,\,\, \lambda \in \mathbb{R}^n,
\]
modulation
\[
f(x) \mapsto e^{2\pi i \,\, \lambda \cdot x}f(x), \,\,\, \lambda \in \mathbb{R}^n,
\]
dilation
\[
f(x) \mapsto f(\lambda x), \,\,\, \lambda \in \mathbb{R},
\]
differentiation
\[
f(x) \mapsto \partial_x f(x),
\]
and multiplication by $x$
\[
f(x) \mapsto xf(x).
\]
\end{theorem}

Before we go on to formulate a more precise form of the smoothed calculus, we need to introduce a final family of spaces.

\begin{definition}[Countably normed spaces of type $\mathcal{S}$]  
\label{def5}
Let $\alpha, \beta \geqslant 0$, \, $A,B>0$. \\
\begin{enumerate}
\item
$f(x) \in \mathcal{S}_{\alpha,A}$ iff \, $\forall  q \in \mathbb{N} \cup \lbrace 0 \rbrace, \, \bar{A}>A, \,\,\, \exists C_{q,\bar{A}}$ depending on $f$ s.t. $\forall k \in \mathbb{N} \cup \lbrace 0 \rbrace$
\begin{equation}
\label{eq44}
| x^k \partial_x^q f(x) | \leqslant C_{q,\bar{A}} \bar{A}^k k^{k \alpha},
\end{equation}
or equivalently $\forall q \in \mathbb{N}\cup \lbrace 0 \rbrace , \, \delta > 0, \, \,\, \exists C_{q,\delta}$,
\begin{equation}
\label{eq45}
| \partial_x^q f(x) | \leqslant C_{q,\delta} e^{-(a(A)-\delta)|x|^{\frac{1}{\alpha}}},
\end{equation}
where $a(A)=\frac{\alpha}{e A^{\frac{1}{\alpha}}}$.

\item
$f(x) \in \mathcal{S}^{\beta,B}$ iff \, $\forall  k \in \mathbb{N} \cup \lbrace 0 \rbrace, \,
\bar{B}>B \,\,\, \exists C_{k,\bar{B}}$ depending on $f$ s.t. $\forall q \in \mathbb{N} \cup \lbrace 0 \rbrace$
\begin{equation}
\label{eq46}
| x^k \partial_x^q f(x) | \leqslant C_{k,\bar{B}} \bar{B}^q q^{q \beta}.
\end{equation}

\item
$f(x) \in \mathcal{S}_{\alpha,A}^{\beta,B}$ iff \, $\forall \bar{A}>A, \bar{B}>B  \,\,\, \exists C_{\bar{A},\bar{B}} $ depending on $f$ s.t. $\forall q,k \in \mathbb{N} \cup \lbrace 0 \rbrace$
\begin{equation}
\label{eq47}
| x^k \partial_x^q f(x) | \leqslant C_{\bar{A},\bar{B}} \bar{A}^k \bar{B}^q k^{k \alpha} q^{q \beta}
\end{equation}

\end{enumerate}
\end{definition}

The new spaces are related to $\mathcal{S}_\alpha, \mathcal{S}^\beta, \mathcal{S}_\alpha^\beta$ by
\begin{eqnarray}
\mathcal{S}^{\beta} & = & \mathop{\cup}\limits_{B>0} {\mathcal{S}^{\beta,B}},\\
\mathcal{S}_{\alpha} & = & \mathop{\cup}\limits_{A>0} {\mathcal{S}_{\alpha,A}},\\
\mathcal{S}_{\alpha}^{\beta} & = & \mathop{\cup}\limits_{A,B>0} {\mathcal{S}_{\alpha,A}^{\beta,B}}.
\end{eqnarray}
Moreover, they are ordered by
\begin{eqnarray}
A_1 \, < \, A_2 \,\, & \Rightarrow & \,\, {\mathcal{S}_{\alpha,A_1}} \subseteq {\mathcal{S}_{\alpha,A_2}}, \\
\label{eqT1}
B_1 \, < \, B_2 \,\, & \Rightarrow & \,\, {\mathcal{S}^{\beta,B_1}} \subseteq {\mathcal{S}^{\beta,B_2}}, \\
\label{eqT2}
B_1 \, < \, B_2, \, A_1 \, < \, A_2 \,\,\, & \Rightarrow & \,\, {\mathcal{S}_{\alpha,A_1}^{\beta,B_1}} \subseteq {\mathcal{S}_{\alpha,A_2}^{\beta,B_2}}.
\end{eqnarray}

Let us now summarize the generalizations of Theorems \ref{thrm6}, \ref{thrm7}, \ref{thrm8} for the new family of spaces here:

\begin{theorem}[Basic properties of $\mathcal{S}_{\alpha,A}, \mathcal{S}^{\beta,B}, \mathcal{S}_{\alpha,A}^{\beta,B}$ ] 
\label{thrm9} 
\begin{enumerate}
\item
{\bf Fourier transforms:}
\begin{equation}
\label{eq48}
\mathcal{F} \left({ \mathcal{S}_{\alpha,A} }\right)=\mathcal{S}^{\alpha,A}.
\end{equation}

\begin{equation}
\label{eq49}
\mathcal{F} \left({ \mathcal{S}^{\beta,B} }\right)=\mathcal{S}_{\beta,B}.
\end{equation}

\begin{equation}
\label{eq50}
\mathcal{F} \left({ \mathcal{S}_{\alpha,A}^{\beta,B} }\right)=\mathcal{S}^{\alpha,A}_{\beta,B}.
\end{equation}

\item
{\bf Nontriviality:} The spaces $\mathcal{S}_{\alpha,A}, \mathcal{S}^{\beta,B}$ are nontrivial for all $\alpha, \beta, A, B >0$. The space $\mathcal{S}_{\alpha,A}^{\beta,B}$ is nontrivial iff
\begin{itemize}
\item[ ]
$\alpha+\beta > 1, \,\,\, \alpha>0, \, \beta>0$, and $A,B > 0$; or 
\item[ ]
$\alpha=0, \, \beta>1, \,\, A,B>0$; or
\item[ ]
$\alpha>1, \, \beta=0, \,\, A,B>0$; or
\item[ ]
$\alpha+\beta=1, \,\, AB>\gamma(\alpha,\beta)$ for some appropriate $\gamma(\alpha,\beta)>0$.
\end{itemize}

\item
{\bf Closedness under elementary operations:} each of $\mathcal{S}_{\alpha,A}, \mathcal{S}^{\beta,B}, \mathcal{S}_{\alpha,A}^{\beta,B}$ is closed under translation, modulation, differentiation 
and multiplication by $x$ \footnote{The operations are precisely defined in the statement of Theorem \ref{thrm8}}. Dilations
scale obviously; if
\[
D_\lambda : f(x) \mapsto f(\lambda x)
\]
then
\begin{eqnarray}
D_\lambda ( \mathcal{S}_{\alpha,A} ) = \mathcal{S}_{\alpha,\frac{A}{\lambda}}, \\
D_\lambda ( \mathcal{S}^{\beta,B} ) = \mathcal{S}^{\beta,\lambda B}, \\
D_\lambda ( \mathcal{S}^{\beta,B}_{\alpha,A} ) = \mathcal{S}^{\beta,\lambda B}_{\alpha,\frac{A}{\lambda}}.
\end{eqnarray}

\end{enumerate}
\end{theorem}

\subsection{Spaces of smoothed functions}
\label{sec2_3}

\noindent {\bf Remark:} {\em In this Section we will work in the general $x \in \mathbb{R}^n$ setup, in contrast to the previous Section, where we only examined $n=1$. This is necessary, since the Wigner transform (which we want to apply our results to) doubles the number of independent variables, e.g. takes $1$-dimensional problems to $2$-dimensional ones. We go on to the more general $n$-dimensional case since it presents no essential difficulties. (Indeed Gelfand and Shilov also present briefly the $n$-dimensional generalization of their theory in Section 9 of chapter IV of \cite{Gel}, after a more detailed study of the $1$-dimensional case).} \footnote{ For example, Definition \ref{def5}, part 1, is generalized as follows:
$f(x) \in \mathcal{S}_{\alpha,(A_1,...,A_n)}$ iff \, $\forall  k,q \in \left({ \mathbb{N} \cup \lbrace 0 \rbrace }\right)^n, \,\, \bar{A}=(\bar{A}_1,...,\bar{A_n}),$  $ \bar{A_l}>A_l \, \forall l=1,...,n, \,\,\, \exists C_{q,\bar{A}}$ depending on $f$ s.t.
\[
| x^k \partial_x^q f(x) | \leqslant C_{q,\bar{A}} \mathop{\Pi}\limits_{l=1}^{n} \bar{A_l}^{k_l} {k_l}^{k_l \alpha},
\]
where of course $x^k=\left({x_l^{k_l}}\right), \, \partial_x^q=\partial_{x_l}^{q_l}$. (Further generalization with $\alpha=(\alpha_1,...,\alpha_n)$ is also possible, and straightforward, but it will not be of interest in this work).
}

\begin{lemma}[The range of $\Phi$, I]
\label{lmm2}
Consider a Schwartz test-function $f(x) \in \mathcal{S}(\mathbb{R}^n)$. Then its image under the smoothing operator belongs to a space of type $\mathcal{S}^{\frac{1}{2},B}$,
\begin{equation}
\label{eq103a}
f(x) \in \mathcal{S}(\mathbb{R}^n) \,\,\, \Rightarrow \,\,\,
\Phi f(x) \in \mathcal{S}^{\frac{1}{2},\frac{1}{\sigma \sqrt{\pi e}}},
\end{equation}
or, more generally,
\begin{equation}
\label{eq103}
f(x) \in \mathcal{S} \,\,\, \Rightarrow \,\,\,
\Phi_{\sigma_1,...,\sigma_n} f(x) \in \mathcal{S}^{\frac{1}{2},\frac{1}{\sqrt{\pi e}} \cdot \left({ \frac{1}{\sigma_1},...,\frac{1}{\sigma_n} }\right) },
\end{equation}
\end{lemma}

\noindent {\bf Proof:} We begin with the proof of equation (\ref{eq103a}); we will prove that $\widehat{\Phi f}(k) \in \mathcal{S}_{\frac{1}{2},\frac{1}{\sigma \sqrt{\pi e}}}$; then the result follows making use of Theorem \ref{thrm9}. 

Observe that, if we denote $g(k)=\widehat{\Phi f}(k)$,
\begin{equation}
\label{eq43}
g(k)=\widehat{\Phi f}(k)=\hat{f}(k) e^{-\frac{\pi}{2} \sigma^2 k^2} \,\,\, \Rightarrow \,\,\, k^p \partial_k^q g(k) =e^{-\frac{\pi}{2}\sigma^2 k^2} \,\,\, \sum\limits_{l=1}^{q} { P_{p,q;l}(k) \hat{f} ^{(l)}(k) } ,
\end{equation}
for certain appropriate polynomials $P_{p,q;l}(k)$ of degree at most $p+q$. $f(x) \in \mathcal{S}$ implies that
\begin{equation}
\label{eq51}
\forall p,q \in \mathbb{N} \cup \lbrace 0 \rbrace \,\,\,\, \exists C_{p,q} \,\, : \,\,
| \sum\limits_{l=1}^{q} { P_{p,q;l}(k) \hat{f} ^{(l)}(k) }  | \leqslant C_{p,q}.
\end{equation}
Equations (\ref{eq43}), (\ref{eq51}) together imply
\begin{equation}
\label{eq52}
| k^p \partial_k^q g(k) | \leqslant C_{p,q} e^{-\frac{\pi}{2} \sigma^2 |k|^2}.
\end{equation}
Recall that equation (\ref{eq45}) states $g(x) \in \mathcal{S}_{\frac{1}{2},\frac{1}{\sigma \sqrt{e \pi}}}$ iff $\forall q \in \mathbb{N}\cup \lbrace 0 \rbrace , \, \delta > 0, \, \,\, \exists C_{q,\delta}$,
\begin{equation}
\label{eq53}
| \partial_x^q g(x) | \leqslant C_{q,\delta} e^{-(\frac{\pi}{2}\sigma^2-\delta)|x|^2};
\end{equation}
so equation (\ref{eq52}) is actually {\em stronger} than $\widehat{\Phi f}(k) \in \mathcal{S}_{\frac{1}{2},\frac{1}{\sigma \sqrt{\pi e}}}$.

The operator $\Phi_{\sigma_1,...,\sigma_n}$ is defined in equation (\ref{PHIa}). Equation (\ref{eq103}) follows in the same way as equation (\ref{eq103a}).

The proof is complete. \\ \vspace{0.3mm}

Observe also that, using equation (\ref{eqT1}), it follows that $\widehat{\Phi f}(k) \in \mathcal{S}_{\frac{1}{2},\frac{1}{\lambda \sqrt{\pi e}}} \,\, \forall \lambda < \sigma$, and accordingly $\Phi f(x) \in \mathcal{S}^{\frac{1}{2},\frac{1}{\lambda \sqrt{\pi e}}} \,\, \forall \lambda < \sigma$. \\ \vspace{1mm}

At this point it is clear that there is one more family of spaces we will use -- and it is closely related to the Gelfand-Shilov test-functions:

\begin{definition}[$G_{\frac{1}{2},\sigma^2}, \,\, G^{\frac{1}{2},\sigma^2}$]
\label{def6}
We will say that $f(x) \in G_{\frac{1}{2},\sigma^2}$ iff $\forall p,q \in \left({ \mathbb{N} \cup \lbrace 0 \rbrace }\right)^n \,\,\, \exists C_{p,q}>0$ s.t.
\begin{equation}
\label{eq54}
| x^p \partial_x^q f(x) | \leqslant C_{p,q} e^{-\frac{\pi}{2}\sigma^2 |x|^2},
\end{equation}
and $f(x) \in G^{\frac{1}{2},\sigma^2}$ iff 
\begin{equation}
\label{eq55}
\hat{f}(k) \in G_{\frac{1}{2},\sigma^2}.
\end{equation}
\end{definition}

\noindent {\bf Remark:} {\em In the anisotropic case $\sigma=(\sigma_1,...,\sigma_n)$ (i.e. the smoothing operator defined in equation (\ref{PHIa}) ), we will use the same notation, namely:

$f(x) \in G_{\frac{1}{2},\sigma^2}$ iff $\forall p,q \in \left({ \mathbb{N} \cup \lbrace 0 \rbrace }\right)^n \,\,\, \exists C_{p,q}>0$ s.t.
\begin{equation}
\label{eq54a}
| x^p \partial_x^q f(x) | \leqslant C_{p,q} e^{-\frac{\pi}{2} \sum\limits_{l=1}^{n} {\sigma_l^2 x_l^2} },
\end{equation}
and $f(x) \in G^{\frac{1}{2},\sigma^2}$ iff 
\begin{equation}
\label{eq55a}
\hat{f}(k) \in G_{\frac{1}{2},\sigma^2}.
\end{equation}
We will not comment on the anisotropic generalization explicitly from now on, since it is straightforward. It lies basically in the recasting of $\sigma^2 x^2$ as $\sum\limits_{l=1}^{n} {\sigma_l^2 x_l^2}$ in the decay conditions. \\
 }

The letter $G$ is chosen to emphasize that these spaces are intimately tied to the Gaussian smoothing we use.

Now it follows that

\begin{lemma}[The range of $\Phi$, II]
\label{lmm3}
\begin{equation}
\label{eq56}
f(x) \in \mathcal{S} \,\,\, \Leftrightarrow \,\,\,
\Phi f(x) \in G^{\frac{1}{2},\sigma^2} \,\,\, \Leftrightarrow \,\,\,
\widehat{\Phi f}(k) \in G_{\frac{1}{2},\sigma^2}
\end{equation}
and moreover
\begin{equation}
\label{eq57a}
\begin{array}{c}
\forall \lambda_1 > \sigma, \, 0<\lambda_2 < \sigma \,:\,\,\,
\mathcal{S}^{\frac{1}{2},\frac{1}{\lambda_1 \sqrt{\pi e}}}  \subsetneq G^{\frac{1}{2},\sigma^2} \subsetneq \mathcal{S}^{\frac{1}{2},\frac{1}{\lambda_2 \sqrt{\pi e}}}, \\ { } \\
G^{\frac{1}{2},\sigma^2} \subseteq \mathcal{S}^{\frac{1}{2},\frac{1}{\sigma \sqrt{\pi e}}},
\end{array}
\end{equation}
where of course $\frac{1}{\lambda_1}=\left({\frac{1}{\lambda_{1,1}},...,\frac{1}{\lambda_{1,n}}}\right)$ and similarly for the other indices.
\end{lemma}

\noindent {\bf Proof:} Equation (\ref{eq56}) essentially follows from the proof of Lemma \ref{lmm2}. Observe that equation (\ref{eq54}) (which holds for any $f(x) \in G^{\frac{1}{2},\sigma^2}$) implies that if $f(x) \in G_{\frac{1}{2},\sigma^2} \Rightarrow e^{\frac{\pi}{2} \sigma^2}f(x) \in \mathcal{S}$, and therefore, by a Fourier transform, if $f(x) \in G^{\frac{1}{2},\sigma^2}$ then $\Phi^{-1} f(x) \in \mathcal{S}$. This shows $\Phi(\mathcal{S}) \supseteq  G^{\frac{1}{2},\sigma^2}$; $\Phi(\mathcal{S}) \subseteq  G^{\frac{1}{2},\sigma^2}$ is straightforward (in other words equations (\ref{eq52}) and (\ref{eq54}) are the same).

We will prove equation (\ref{eq57a}) in stages; first of all we will show that 
\begin{equation}
\label{eq57b}
\begin{array}{c}
\forall \lambda_1 > \sigma \,:\,\,\,
\mathcal{S}^{\frac{1}{2},\frac{1}{\lambda_1 \sqrt{\pi e}}}  \subsetneq G^{\frac{1}{2},\sigma^2}.
\end{array}
\end{equation}
Observe that (as we saw earlier in Definition \ref{def5}, and equation (\ref{eq53}) ),
\begin{equation}
\label{eq99}
\begin{array}{c}
f(x) \in \mathcal{S}^{\frac{1}{2},\frac{1}{\lambda_1 \sqrt{\pi e}}} \,\, \Leftrightarrow \,\, \\ 
\forall q \in \left({ \mathbb{N}\cup \lbrace 0 \rbrace }\right)^n, \, \delta > 0, \, \,\, \exists C_{q,\delta} : \,\,\,
| \partial_x^q \hat{f}(x) | \leqslant C_{q,\delta} e^{-(\frac{\pi}{2}\lambda_1^2-\delta)x^2} \,\, \Leftrightarrow \,\, \\ 
\forall p,q \in \left({ \mathbb{N}\cup \lbrace 0 \rbrace }\right)^n , \, \delta > 0, \, \,\, \exists C_{q,\delta} : \,\,\,
| x^p \partial_x^q \hat{f}(x) | \leqslant C_{p,q,\delta} e^{-(\frac{\pi}{2}\lambda_1^2-\delta)x^2},
\end{array}
\end{equation}
where $C_{p,q,\delta}$ can be chosen not larger than $C_{p,q,\delta} \leqslant C_{q,\frac{\delta}{2}} \cdot \mathop{sup}\limits_{x \in \mathbb{R}^n} \left[{ |x|^p \, e^{-\frac{\delta}{2}x^2 } }\right]$. Therefore
\begin{equation}
\label{eq100}
\begin{array}{c}
f(x) \in \mathcal{S}^{\frac{1}{2},\frac{1}{\lambda_1 \sqrt{\pi e}}}, \, \lambda_1>\sigma \,\, \Rightarrow \,\, \\
\forall p,q  \in \left({ \mathbb{N}\cup \lbrace 0 \rbrace }\right)^n \exists C'_{p,q} \, : \,\,
|x^p \partial_x^q \hat{f}(x)| \leqslant C'_{p,q} e^{-\frac{\pi}{2}\sigma^2 x^2}
\end{array}
\end{equation}
for $C'_{p,q} \leqslant C_{p,q,\pi \frac{\lambda_1^2-\sigma^2}{2}}$ \footnote{that is $C'_{p,q} \leqslant C_{p,q,\delta}$ where $C_{p,q,\delta}$ are the constants of the same notation in equation (\ref{eq99}) for $\delta = \pi \frac{\lambda_1^2-\sigma^2}{2}$},
and therefore, for $\lambda_1>\sigma$, $\mathcal{S}^{\frac{1}{2},\frac{1}{\lambda_1 \sqrt{\pi e}}}  \subseteq G^{\frac{1}{2},\sigma^2}$. On the other hand $g_1(x)=\mathcal{F}^{-1}_{k \shortrightarrow x} \left[{ e^{-\frac{\pi}{2}\sigma^2 k^2} }\right] \in G^{\frac{1}{2},\sigma^2} \setminus \mathcal{S}^{\frac{1}{2},\frac{1}{\lambda_1 \sqrt{\pi e}}} $, and now equation (\ref{eq57b}) follows.

That
\begin{equation}
\label{eq101}
\begin{array}{c}
\forall \lambda_2 \leqslant \sigma \,:\,\,\,
 G^{\frac{1}{2},\sigma^2} \subseteq \mathcal{S}^{\frac{1}{2},\frac{1}{\lambda_2 \sqrt{\pi e}}}  
\end{array}
\end{equation}
was shown in Lemma \ref{lmm2} (it was shown there for $\lambda_2=\sigma$; for $\lambda_2 < \sigma$ it follows making use of equation (\ref{eqT1})). Also,
$\forall \lambda_2 < \sigma \,:\,\,\,
 G^{\frac{1}{2},\sigma^2} \subsetneq \mathcal{S}^{\frac{1}{2},\frac{1}{\lambda_2 \sqrt{\pi e}}}$ is obvious. 

The proof is complete \footnote{It is not clear at this point whether $\mathcal{S}^{\frac{1}{2},\frac{1}{\sigma \sqrt{\pi e}}}$ is strictly larger than $ G^{\frac{1}{2},\sigma^2}$, or if the two spaces are equal.}. \\ \vspace{1mm}

Lemma \ref{lmm3} highlights 
the close relation of smoothed functions with Gelfand-Shilov spaces. This relation allows us to use (after some adaptation) a lot of existing theory -- most notably extensions in the complex plane. A very useful result has to do with the rate of growth of $G^{\frac{1}{2},\sigma^2}$ functions on the complex plane: since $\Phi(\mathcal{S}) \subseteq  \mathcal{S}^{\frac{1}{2},B} \subseteq \mathcal{S}^\beta$, estimates like the one in equation (\ref{eq38}) hold automatically. However, as we saw smoothed functions satisfy somewhat stronger conditions, and accordingly we can construct somewhat stronger estimates (which turn out to be necessary in the sequel). 

First of all however, we need a technical Lemma:

\begin{lemma}[Fourier-domain representation of entire functions]
\label{lmm4}
Let $f(x) \in \mathcal{S}^{\frac{1}{2}}(\mathbb{R}^n)$. 
We know that then $f(x)$ can be extended to an entire function on the complex domain, $f(x+iy)$. In addition, the following representation is valid
\begin{equation}
\label{eq75}
f(x+iy)=\int\limits_{k \in \mathbb{R}^n} { e^{2\pi i k \, \cdot \, (x+iy)} \hat{f}(k)dk}.
\end{equation}
\end{lemma}

\noindent {\bf Remark: } {\em As we just saw (lemma \ref{lmm2}) smoothed test-functions belong in $\mathcal{S}^{\frac{1}{2}}$; therefore Lemma \ref{lmm4} applies to them.}

\noindent {\bf Proof:} Denote $z=x+iy$; the proof lies with the computation
\begin{equation}
\label{eq76}
\begin{array}{c}
\int\limits_{k \in \mathbb{R}^n} { e^{2\pi i k \, \cdot \, (x+iy)} \hat{f}(k) dk}=
\int\limits_{k \in \mathbb{R}^n} { \sum\limits_{l=0}^{\infty} {\frac{ \left[{2 \pi i k \, \cdot \, (x+iy)}\right]^l }{l!}} \hat{f}(k)dk}=\\ { } \\
=\int\limits_{k \in \mathbb{R}^n} { \sum\limits_{l=0}^{\infty} \left({2 \pi i}\right)^l {\frac{
\left[{ \sum\limits_{s=1}^{n} { k_s (x_s+iy_s)}  }\right]^l }{l!}} \hat{f}(k)dk}=\\ { } \\
=\int\limits_{k \in \mathbb{R}^n} {\sum\limits_{l=0}^{\infty} 
{\left({2 \pi i}\right)^l  \frac{\sum\limits_{|r|=l} 
{ \binom{l}{r_1,...,r_n} \, k^{r_1}_1 \left({ x_1 + i y_1 }\right)^{r_1}
\,...\, k^{r_n}_n \left({ x_n + i y_n }\right)^{r_n} }}{l!} }
\hat{f}(k)dk}=\\ { } \\
=\sum\limits_{l=0}^{\infty}
{
\sum\limits_{|r|=l}
{ \frac{ \left({2 \pi i}\right)^l}{r!} \, z^r \,
\int\limits_{k \in \mathbb{R}^n}
{ 
k^r
\hat{f}(k)dk}}}=
\sum\limits_{l=0}^{\infty}
{
\sum\limits_{|r|=l}
{ \frac{ z^r}{r!}
\partial_x^r f(0)
}}=f(z).
\end{array}
\end{equation}
We have used the standard multi-index notation in the computations \footnote{ $r \in \left({ \mathbb{N} \cup \lbrace 0 \rbrace }\right)^n$ is a multi-index; $|r|=r_1+...+r_n$, $r!=r_1! \cdot ... \cdot r_n!$, $k^r=k_1^{r_1} \cdot  ... \cdot k_n^{r_n}$, $z^r=\left({ x_1^{r_1} +i y_1^{r_1} }\right)  \cdot  ... \cdot \left({ x_n^{r_n} +i y_n^{r_n} }\right)$, $\partial_x^r=\partial_{x_1}^{r_1} \cdot  ... \cdot \partial_{x_n}^{r_n}$. } . 

In equation (\ref{eq76}) above it is first of all seen that the bulky expressions coming from the high-dimensional character of the problem can be nicely summarized as
\begin{equation}
\int\limits_{k \in \mathbb{R}^n} { e^{2\pi i k \, \cdot \, (x+iy)} \hat{f}(k) dk}=
\int\limits_{k \in \mathbb{R}^n} { \sum\limits_{l \in \left({ \mathbb{N} \cup \lbrace 0 \rbrace }\right)^n} {\frac{ 2 \pi i k^l (x+iy)^l }{l!}} \hat{f}(k)dk}
\end{equation}
etc.

The condition we have to check to justify the interchange of summation and integration (in the more compact notation) is
\begin{eqnarray}
\label{eq77}
\forall x,y \in \mathbb{R} \, : \,\,\, \sum\limits_{l \in \left({ \mathbb{N} \cup \lbrace 0 \rbrace }\right)^n} {\frac{  |x+iy|^l }{l!} \int\limits_{k \in \mathbb{R}^n}  {  | \left({2\pi i k}\right)^l \hat{f}(k) | dk } } < \infty.
\end{eqnarray}

When we first expand in a series, in the first line of equation (\ref{eq76}), it is the Taylor series of the exponential; when we summed analytically the series in the last line, this was the Taylor expansion of an entire-analytic function. That is, all the series converge, and the commutation of the series and the integral follows from the dominated convergence Theorem.  \\ \vspace{1mm}

\noindent {\bf Remark:} {\em Lemma \ref{lmm4} essentially is already justified for any entire-analytic function. However, since it is not too long, and to provide more insight, we will give here a more detailed proof of the result. In particular this highlights how the estimates on the entire function are inherited in this series -- something we will return to.} \\ \vspace{1mm}

Observe that $f(x) \in \mathcal{S}^{\frac{1}{2}} \,\, \Rightarrow \,\, \hat{f}(k) \in \mathcal{S}_{\frac{1}{2}}$, (according to Theorem \ref{thrm6}, part 2), and therefore (according to the high-dimensional version of Definition \ref{def4}), $\exists C,A=(A_1,...,A_n)>0 \, : \, \forall l \in \left({ \mathbb{N} \cup \lbrace 0 \rbrace }\right)^n$
\begin{eqnarray}
\label{eq78}
|\left({ 2\pi ik }\right)^l \hat{f}(k)| \leqslant C A^{l} l^{\frac{l}{2}}=C \mathop{\Pi}\limits_{s=1}^{n} {A_s^{l_s} l_s^{\frac{l_s}{2}}}
\end{eqnarray}

Now we have
\begin{equation}
\begin{array}{c}
\label{eq79}
\sum\limits_{ l \in \left({ \mathbb{N} \cup \lbrace 0 \rbrace }\right)^n } {\frac{  |x+iy|^l }{l!} \int\limits_{k \in \mathbb{R}^n}  {  | \left({2\pi i k}\right)^l \hat{f}(k) | dk } } = \\
=  \sum\limits_{ l \in \left({ \mathbb{N} \cup \lbrace 0 \rbrace }\right)^n } {
\frac{  |x+iy|^l }{l!} \left[{ \int\limits_{|k| \leqslant \frac{1}{2\pi} }  {  | \left({2\pi i k}\right)^l \hat{f}(k) | dk } + \int\limits_{|k| > \frac{1}{2\pi}}  {  | \left({2\pi i k}\right)^l \hat{f}(k) | dk } }\right] } \leqslant \\ { } \\
\leqslant C_1 e^{ |x+iy| } + \sum\limits_{ l \in \left({ \mathbb{N} \cup \lbrace 0 \rbrace }\right)^n } {
\frac{  |x+iy|^l }{l!} \left[{ \int\limits_{|k| > \frac{1}{2\pi} }  {  | \left({2\pi i k}\right)^{(l+2)} \hat{f}(k) | \frac{1}{|x^2|} dk }}\right] } \leqslant \\ { } \\
\leqslant C_1 e^{ |x+iy| } + C_2 \left[{ \int\limits_{|k| > \frac{1}{2\pi}}  { \frac{1}{|x^2|} dk }}\right] \, \sum\limits_{ l \in \left({ \mathbb{N} \cup \lbrace 0 \rbrace }\right)^n } {
\frac{  |x+iy|^l }{l!} A^{(l+2)} (l+2)^{\frac{l+2}{2}}   },
\end{array}
\end{equation}
where in the last step we made use of equation (\ref{eq78}). So now we only have to check the absolute convergence of the series
\begin{equation}
\label{eq80}
\sum\limits_{ l \in \left({ \mathbb{N} \cup \lbrace 0 \rbrace }\right)^n } {
\frac{  |x+iy|^l }{l!} A^{(l+2)} (l+2)^{\frac{l+2}{2}}  }=
A^2
\sum\limits_{ l \in \left({ \mathbb{N} \cup \lbrace 0 \rbrace }\right)^n } {
\frac{ (l+2)^{\frac{l+2}{2}}  }{l!} |(x+iy)A|^{l}    },
\end{equation}
or, equivalently, find the radius of convergence of the power series
\begin{equation}
\label{eq80a}
\begin{array}{c}
\sum\limits_{ l \in \mathbb{N}^n } {
a_l z^l   }=\sum\limits_{ l \in \left({ \mathbb{N} \cup \lbrace 0 \rbrace }\right)^n } {
\mathop{\Pi}\limits_{s=1}^{n} a_{l_s} z^{l_s}  } \, ; \\{ } \\
a_l = \frac{ (l+2)^{\frac{l+2}{2}}  }{l!}=\mathop{\Pi}\limits_{s=1}^{n}  {\frac{ (l_s+2)^{\frac{l_s+2}{2}}  }{l_s!}},
\end{array}
\end{equation}
to be infinite. \\

\noindent {\bf Case $n=1$:}
Observe that, making use of the Stirling approximation, it follows that
\begin{equation}
\label{eq82a}
(l+2)^{\frac{l+2}{2}} \approx \sqrt{ \frac{(l+2)! \, e^{(l+2)}} {\sqrt{2\pi (l+2)}} }=
\frac{e}{(2\pi)^\frac{1}{4}} \, \left({ \sqrt{e} }\right)^{l} \, (l+2)^{\frac{1}{4}} \, (l+1)^{\frac{1}{2}} \,  (l!) ^{\frac{1}{2}},
\end{equation}
and therefore, using the Stirling approximation once more,
\begin{equation}
\begin{array}{c}
\label{eq83a}
a_l \approx \frac{e}{(2\pi)^\frac{1}{4}} \frac{  \left({ \sqrt{e} }\right)^{l} \, (l+2)^{\frac{1}{4}} \, (l+1)^{\frac{1}{2}}  }{(l!) ^{\frac{1}{2}}} \approx
\frac{e}{(2\pi)^\frac{1}{2}} \frac{e^l \, (l+2)^{\frac{1}{4}} \, (l+1)^{\frac{1}{2}} \, }{ \sqrt{l}^l } =\\ { } \\
=\frac{e}{(2\pi)^\frac{1}{2}} \, (l+2)^{\frac{1}{4}} \, (l+1)^{\frac{1}{2}} \, \left({ \frac{e}{\sqrt{l}} }\right)^l .
\end{array}
\end{equation}

It is now obvious that $\forall R_0>0 \,\, \exists C=C(R_0)$ such that 
\begin{equation}
\begin{array}{c}
\label{eq83b}
|a_l| \leqslant C \cdot R_0^{-l}.
\end{array}
\end{equation}
But an estimate like this implies that the radius of convergence for the power series with coefficients $a_l$ is at least $R_0$; therefore the radius of convergence is infinite, and the series of equation (\ref{eq80}) always converges. \\

\noindent {\bf General case, $n \in \mathbb{N}$:}
This is a straightforward generalization of the previous computation. First of all, using the Stirling formula like earlier it follows that 
\begin{equation}
\begin{array}{c}
\label{eq83d}
a_l \approx 
\mathop{\Pi}\limits_{s=1}^n {
\frac{e}{(2\pi)^\frac{1}{2}} \, (l_s+2)^{\frac{1}{4}} \, (l_s+1)^{\frac{1}{2}} \, \left({ \frac{e}{\sqrt{l_s}} }\right)^{l_s} },
\end{array}
\end{equation}
and therefore
$\forall R_1,..,R_n > 0 \,\, \exists C=C(R_1,...,R_n)$ such that 
\begin{equation}
\begin{array}{c}
\label{eq83c}
|a_l| \leqslant C \cdot \mathop{\Pi}\limits_{s=1}^n { R_s^{-l_s} }.
\end{array}
\end{equation}
It is a standard (and easy to show) Lemma that equation (\ref{eq83c}) implies that the power series with coefficients $a_l$ (i.e. the series of equation (\ref{eq80a})) converges whenever $|z_s| < R_s$. Since the $R_s$'s can be chosen arbitrarily, it follows that the series of equation (\ref{eq80}) converges always in the multidimensional case as well. 

The proof is complete. \\

Let us also remark that integrals like the one of equation (\ref{eq75}) have also been studied under the name two-sided Laplace transforms \cite{Zem}.

\begin{theorem}[Behaviour of smoothed functions on the complex plane]
\label{cor1}
\label{thrm10}
 Let $f(x) \in G^{\frac{1}{2},\sigma^2}$. As we saw earlier, it can be extended to an entire function on $\mathbb{C}^n$; moreover
\begin{equation}
\label{eq58a}
g(x,y)=e^{- \frac{2 \pi}{\sigma^2} y^2} f(x+iy)   \,\, \in \,\, \mathcal{S}(\mathbb{R}^{2n}).
\end{equation}
\end{theorem}

\noindent {\bf Proof: } First of all observe the following elementary identity
\begin{equation}
\label{eq58a1}
\mathcal{F}_{x \shortrightarrow k} \left[ (2\pi i x)^p \partial_x^q f(x) \right] = \left( -\partial_k \right)^p \left( -2\pi i x \right)^q \hat{f}(k).
\end{equation}

We will use the observation of equation (\ref{eq58a1}) and Lemma \ref{lmm4} to prove equation (\ref{eq58}). We begin from a slightly different point, i.e.
\begin{equation}
\label{eq59n}
\begin{array}{c}
\left( 2\pi i x\right) ^p \left( 2\pi i(x+iy)\right) ^q \partial_x^r \partial_y^s  f(x+iy) = \\ { } \\
=\left( 2\pi i(x+iy)\right)^q \partial_y^s 
\int\limits_{k \in \mathbb{R}^n} { e^{2\pi i k (x+iy)} (-\partial_k)^p (-2\pi i k)^r \hat{f}(k)dk }= \\ { } \\
=\left( 2\pi i(x+iy)\right)^q 
\int\limits_{k \in \mathbb{R}^n} { e^{2\pi i k (x+iy)} (-2\pi k)^s (-\partial_k)^p (-2\pi i k)^r \hat{f}(k)dk }= \\ { } \\
=\int\limits_{k \in \mathbb{R}^n}
{ \left[ \partial_k^q e^{2\pi i k (x+iy)} \right] (-2\pi k)^s (-\partial_k)^p (-2\pi i k)^r \hat{f}(k)dk }= \\ { } \\
=(-1)^q  \int\limits_{k \in \mathbb{R}^n}
{  e^{2\pi i k (x+iy)} \partial_k^q (-2\pi k)^s (-\partial_k)^p (-2\pi i k)^r \hat{f}(k)dk }.
\end{array}
\end{equation}
It follows now, making use of the Definition \ref{def6} of $G^{\frac{1}{2},\sigma^2}$, that there is a constant $C'_{p,q,r,s}$ such that
\begin{equation}
\label{eq59an}
\begin{array}{c}
| \left( 2\pi i x\right) ^p \left( 2\pi i(x+iy)\right) ^q \partial_x^r \partial_y^s  f(x+iy) | \, \leqslant \,
C'_{p,q,r,s} \int\limits_{k \in \mathbb{R}^n} {e^{-2\pi k \cdot y-\frac{\pi}{2}\sigma^2 k \cdot k} dy } = C'_{p,q,r,s} \left(  \frac{\sqrt{2}}{\sigma} \right) e^{\frac{2\pi}{\sigma^2} |y|^2 }.
\end{array}
\end{equation}

The commutation of the $dk$ integral and the $\partial_k$ derivative in equation (\ref{eq59n}) follow from the dominated convergence Theorem and the bounds for $G^{\frac{1}{2},\sigma^2}$ functions.

It is obvious how the following inequality follows from equation (\ref{eq59an}): \\
$\forall \,\, p,q,r,s \in \left( \mathbb{N} \cup \lbrace 0 \rbrace \right)^n \,\,\, \exists C_{p,q,r,s}$ such that
\begin{equation}
\label{eq58}
| \left( 2\pi i x\right) ^p \left( 2\pi iy\right) ^q \partial_x^r \partial_y^s  f(x+iy) | \leqslant C_{p,q,r,s} \,  e^{ \frac{2 \pi}{\sigma^2} y^2}
\end{equation}

In order to complete the proof, observe that for any differential operator with polynomial coefficients $\mathbb{P}$, equation (\ref{eq58}) implies that there is a constant $C(\mathbb{P}) > 0$ such that
\begin{equation}
\label{eq58ab}
| \mathbb{P}  f(x+iy) | \leqslant C(\mathbb{P}) \,  e^{ \frac{2 \pi}{\sigma^2} y^2}.
\end{equation}
Observe in addition, that for each differential operator with polynomial coefficients $\mathbb{P}$ there exists a different differential operator with polynomial coefficients $\tilde{\mathbb{P}}$ such that
\[
\mathbb{P} \left( e^{-\frac{2\pi}{\sigma^2}} g(x,y) \right) = e^{-\frac{2\pi}{\sigma^2}} \tilde{\mathbb{P}} \left( \, g(x,y) \, \right).
\]
Moreover if $\mathbb{P}$ is of order $s$ as a differential operator, and its coefficients are polynomials of degree up to $t$, then $\tilde{\mathbb{P}}$ will be still of order $s$ as a differential operator, and its coefficients will be polynomials of degree up to $s+t$.

Setting $\mathbb{P}=\left( 2\pi i x\right) ^p \left( 2\pi iy\right) ^q \partial_x^r \partial_y^s$ we get
\begin{equation}
\label{eq58aa}
\begin{array}{c}
| \mathbb{P} \left( e^{- \frac{2 \pi}{\sigma^2} y^2} f(x+iy) \right) | =
e^{ -\frac{2 \pi}{\sigma^2}y^2} | \tilde{\mathbb{P}} \left( f(x+iy) \right) | \leqslant 
 C(\tilde{\mathbb{P}}) 
\end{array}
\end{equation}

The proof is complete. \\ \vspace{1mm}

\begin{corollary}
\label{cor2}
 Let $f(x)$ be a smoothed test-function, i.e. $f(x) \in G^{\frac{1}{2},\sigma^2}$. Take a fixed $x \in \mathbb{R}^n$; then
\begin{equation}
\label{eq61a}
g(y)=e^{-2\pi \sigma^2 y^2} f(x+i \sigma^2 y) \, \in  \, \mathcal{S}
\end{equation}
as a function of $y$. Moreover, for any fixed $y$,
\begin{equation}
\label{eq64}
g(x)= f(x+iy) \, \in  \, \mathcal{S}
\end{equation}
as a function of $x$.
\end{corollary}

\subsection{Equivalent formulations of the smoothed calculus} 
\label{sec2_4}

First of all, let us make a remark concerning the Weyl symbols for the smoothed calculus:

\begin{lemma}[Imaginary translations of distributions]
\label{lmm5}
Let $L(x) \in \mathcal{S}'$, $y \in \mathbb{R}^n$. Then $L(x+iy)$ is a well defined functional on $G^{\frac{1}{2},\sigma^2}$.
\end{lemma}

\noindent {\bf Proof:} Take $f(x) \in G^{\frac{1}{2},\sigma^2}$. Now
\begin{equation}
\label{eq85}
\int\limits_{x \in \mathbb{R}^n} {L(x+iy)f(x)dx}=
\int\limits_{x \in \mathbb{R}^n} {L(x)f(x-iy)dx},
\end{equation}
which is well defined since $f(x-iy) \in \mathcal{S}$, as we saw in Corollary \ref{cor2}. The proof is complete.

So, basically, the idea is that since we act on ``very nice functions'' we can have more operations on our distributions, which will be interpreted weakly. Observe that the point of Lemma \ref{lmm5} has {\em absolutely nothing} to do with actually extending $L(x)$ into the complex plane. Giving meaning to the Weyl symbols of convolution-deconvolution sandwiches is somewhat more interesting, since we don't take a fixed imaginary translation, but go over an imaginary axis. 

In any case, it seems tempting to ask whether $L(x+ \frac{ i \sigma^2 \, k}{2},k)$ simply belongs to an ultra-distribution space (i.e. to the dual of some space of the type e.g. $\mathcal{S}^{\beta,B}$). It seems probable that smoothed operators can be cast in a satisfactory framework simply as operators with ultra-distributional Weyl symbols. \\ \vspace{0.5mm}

Now we go on to the equivalent formulations of the smoothed calculus, making use of the properties of smoothed functions that we just proved.

\begin{theorem}[Smoothed calculus]
\label{thrm11}
Let $f(x) \in \mathcal{S}(\mathbb{R})$, $L(x) \in \mathcal{S}'(\mathbb{R}^{2n})$, $L$ be the operator with Weyl symbol $L(x,k)$ and $w(x)=\Phi f(x)$.
Then 
\begin{eqnarray}
\label{eq66}
\Phi L \Phi^{-1} w (x) =  \int\limits_{y,k \in \mathbb{R}^n} { e^{2\pi i (x-y)k -2\pi \sigma^2 k^2} L \left( \frac{x+y}2,k \right) w(y-i\sigma^2 k) dydk}.
\end{eqnarray}
\end{theorem}

\noindent {\bf Remarks:} {\em Before we go on to the proof, some comments should be made:
\begin{itemize}
\item
Equation (\ref{eq66}) is well defined, through Theorem \ref{cor1}.
\item
This form makes clear what we gain by computing explicitly the convolution-deconvolution sandwich, as opposed to applying $\Phi L \Phi^{-1}$ successively as three different operators: in order to implement the sandwich we need to compute / implement imaginary translations $w(x) \mapsto w(x+iy)$, which correspond to a Fourier multiplier $e^{-2\pi ky}$,  but not deconvolutions, which correspond to a much stronger Fourier multiplier, $e^{\frac{\pi}{2}\sigma^2 k^2}$.
\end{itemize} }

\noindent {\bf Proof:} It suffices to show that equation (\ref{eq66}) defines the same operator as equation (\ref{eq13}),
\begin{equation}
\label{eq67}
\Phi L \Phi^{-1} w (x)=2^n \int\limits_{k,u \in \mathbb{R}^n} {
e^{ 2\pi  (ix-\sigma^2 k) u+ 2 \pi i k x } \hat{w}(k-u) \hat{L}_1 (2u,k)dudk
}, 
\end{equation}
where always $\hat{L}_1 (u,k)=\mathcal{F}_{x \shortrightarrow u} \left[{ L(x,k) }\right]$. We start from the expression of equation (\ref{eq66}):
\begin{equation}
\begin{array}{c}
\label{eq67a}
\int\limits_{y,k \in \mathbb{R}^n} { e^{2\pi i (x-y)k -2\pi \sigma^2 k^2} L \left( \frac{x+y}2,k \right) w(y-i\sigma^2 k) dydk}=\\ { } \\

=\int\limits_{y,k \in \mathbb{R}^n} { e^{2\pi i (x-y)k -2\pi \sigma^2 k^2} L \left( \frac{x+y}2,k \right) 
\mathcal{F}_{u \shortrightarrow y} \left[ e^{2\pi \sigma^2ku} \hat{w}(u) \right] dydk}=\\ { } \\

=\int\limits_{u,y,k \in \mathbb{R}^n} { e^{2\pi i (x-y)k -2\pi \sigma^2 k^2 + 2\pi i uy +2\pi\sigma^2 k u} L \left( \frac{x+y}2,k \right)  \hat{w}(u) dudydk}=\\ { } \\

=\int\limits_{u,y,k \in \mathbb{R}^n} { e^{-2\pi i y(k-u)} L \left( \frac{x+y}2,k \right) dy \,
e^{2\pi i xk+2\pi \sigma^2 k(u-k)}   \hat{w}(u) dudydk}=\\ { } \\

=2^n \int\limits_{u,y,k \in \mathbb{R}^n} { e^{2\pi i (k-u)x} \hat{L}_1 \left( 2(k-u),k \right) 
e^{2\pi i xk+2\pi \sigma^2 k(u-k)}   \hat{w}(u) dudydk}=\\ { } \\

=2^n \int\limits_{u,y,k \in \mathbb{R}^n} { e^{2\pi i \left[ ux+xk+i\sigma^2uk \right]} \hat{L}_1 \left( 2u,k \right)  
   \hat{w}(k-u) dudydk}= \\ { } \\
=\Phi L \Phi^{-1} w (x).
\end{array}
\end{equation}

One more equivalent formulation exists when the Weyl symbol $L$ is a differential operator, which is somewhat simpler:
\begin{theorem}[Smoothed calculus, a reformulation for differential operators]
\label{thrm11diff}
Consider the same assumptions for $f$, $L$ as in Theorem \ref{thrm11} above, and in addition let us suppose that $L(x,k)$ is a continuous function of $(x,k)$, and $\forall{x \in \mathbb{R}^n}$ $L_x(k)=L(x,k)$ is (the restriction to the real numbers of) an entire-analytic function, and moreover $\forall x\in \mathbb{R}^n \,\, G(k,y)=L \left( \frac{x+y}2,k+\frac{i(x-y)}{\sigma^2} \right) \in \mathcal{S}'(\mathbb{R}^n)$. For example differential operators, 
\begin{equation}
\label{ddop}
L(x,k)=\sum\limits_{m=0}^N A_m(x) k^m,
\end{equation}
fall in this category. Denote also $w=\Phi f$. Then
\begin{eqnarray}
\label{eq66b}
\Phi L \Phi^{-1} w (x) =  
\int\limits_{k \in \mathbb{R}^n} { F(x,k) e^{ - 2\pi \sigma^2 k^2} w(x-i\sigma^2 k) dk},
\end{eqnarray}
where
\begin{equation}
\label{eq66b1}
F(x,k)=\int\limits_{y \in \mathbb{R}^n} {
e^{-2\pi i (x-y)k } L(\frac{x+y}2,k+\frac{i(x-y)}{\sigma^2}) dy }.
\end{equation}
\end{theorem}

\noindent {\bf Proof: } It is clear that under our assumptions the statement of the Theorem makes sense. For the proof, it suffices to make the change of variables
\[
k = k' +\frac{i(x-y)}{\sigma^2}
\]
in equation (\ref{eq66}):
\begin{equation}
\label{eq66b2}
\begin{array}{c}
\Phi L \Phi^{-1} w (x) =  \int\limits_{y,k \in \mathbb{R}^n} { e^{2\pi i (x-y)k -2\pi \sigma^2 k^2} L \left( \frac{x+y}{2},k \right) w(y-i\sigma^2 k) dydk}= \\ { } \\

=  \int\limits_{y,k' \in \mathbb{R}^n} { e^{2\pi i (x-y)\left(k' + \frac{i(x-y)}{\sigma^2} \right) -2\pi \sigma^2 \left( k +\frac{i(x-y)}{\sigma^2} \right)^2} L \left( \frac{x+y}2,k+\frac{i(x-y)}{\sigma^2} \right) dy w(x-i\sigma^2 k) dk'}= \\ { } \\

=  \int\limits_{y,k \in \mathbb{R}^n} { e^{-2\pi i (x-y)k} L \left( \frac{x+y}2,k+\frac{i(x-y)}{\sigma^2} \right) dy \, e^{-2\pi \sigma^2 k^2} w(x-i\sigma^2 k) dk}.
\end{array}
\end{equation}

If we have a differential operator as in equation (\ref{ddop}) the last expression is equal to
\begin{equation}
\label{eq66b3}
\begin{array}{c}
\Phi L \Phi^{-1} w (x) =  \int\limits_{k \in \mathbb{R}^n} {\left[ \sum\limits_{m=0}^N \int\limits_{y \in \mathbb{R}^n} {{ e^{-2\pi i (x-y)k} A_m \left( \frac{x+y}2\right) \left( k+\frac{i(x-y)}{\sigma^2} \right)^m dy}} \right] \, e^{-2\pi \sigma^2 k^2} w(x-i\sigma^2 k) dk }.
\end{array}
\end{equation}

The proof is complete.

\section{Smoothed Wigner homogenization}
\label{sec3}

\subsection{Smoothed Wigner calculus}
\label{sec3_1}

In this Subsection we will derive the smoothed Wigner calculus, which, as we briefly described in Section \ref{sec1_2}, allows for the derivation of smoothed Wigner equations and a smoothed trace formula. This work essentially follows the same lines as Theorem \ref{thrm3} and its proof, being somewhat more complicated due to the specifics of the Wigner calculus.

\begin{definition}[The Wigner transform]
\label{DefWT}
We define the Wigner transform (WT) as the sesquilinear transform
\begin{equation}
\label{eq20}
W: f(x),g(x) \mapsto W[f,g](x,k)=\int\limits_{y \in \mathbb{R}^n} { e^{-2\pi i k y} f\left({ x+\frac{y}{2}}\right) \bar{g}\left({ x-\frac{y}{2}}\right) dy}.
\end{equation}
The generalization to vectors is straightforward, i.e. if $f(x),g(x):\mathbb{R}^n \shortrightarrow \mathbb{C}^d$ then
\begin{equation}
\label{eq21}
\left[{ W[f,g](x,k) }\right]_{i,j}=W[f_i,g_j](x,k)  \,\,\,\,\,\, i,j \in \lbrace 1,...,d \rbrace
\end{equation}
The WT is well defined and continuous as a bilinear mapping
\begin{eqnarray}
W: \mathcal{S}(\mathbb{R}^n) \times \mathcal{S}(\mathbb{R}^n) \shortrightarrow \mathcal{S}(\mathbb{R}^{2n}), \\
W: L^2(\mathbb{R}^n) \times L^2(\mathbb{R}^n) \shortrightarrow L^2(\mathbb{R}^{2n}).
\end{eqnarray}
\end{definition}

The WT has a number of properties which allow the interpretation of its quadratic version $W[u](x,k)=W[u,u](x,k)$ (often called {\em the Wigner  distribution of $u$} to avoid confusion) as a ``time-frequency energy quasi-density''. That is
\begin{equation}
\label{eq22}
\int\limits_{(x,k) \in A} {W[u](x,k)dxdk}
\end{equation}
is somehow {\em ``proportional to the energy ($L^2$ norm density) corresponding to the to the wavenumbers $k$ at the locations $x$ for $(x,k) \in A$''}. Making precise this interpretation (and understanding its limitations) is a classic topic in time-frequency analysis \cite{Fl1,Groc}, and there is no need to stay on it too long here. 

One of the first findings however, is that the WT exhibits so-called ``interference terms'', i.e. fast oscillations in phase-space, which severely limit its numerical and intuitive use. The interference terms are due to the non-linearity of the transform; for example in certain many-component signals (such as finite sums of Gaussian wavepackets) the ``bad terms'' can be exactly isolated as the cross-terms,
\begin{equation}
\label{eq23}
\begin{array}{c}
W \left[{ \sum\limits_{m=1,...,M} { g_m }}\right](x,k)=\sum\limits_{m=1,...,M} { W[g_m](x,k) }+\\
+2 Re \left[{ \sum\limits_{p=1,...,M} { \sum\limits_{ q<p }{ W[g_p,g_q](x,k) }} }\right].
\end{array}
\end{equation}

In most cases however, isolating explicitly the ``bad part'' is not possible; the term ``auto-interference'' is used to emphasize that. The oscillations in phase-space are in general at least as fast as the oscillations in $u$ (i.e. comparable wavelengths), {\em but can be arbitrarily faster} \footnote{Take $f(x)=e^{-x^2}, g(x)=e^{-(x-a)^2}$. Then $W[f+g](x,k)$ has oscillations with wavelengths of order $\frac{1}{a}$, while the function $f+g$ itself is not really oscillatory at all.}. This makes absolutely necessary some step of regularization; indeed in most applications of the WT some additional regularization device is proposed, be it convolution with a smooth kernel (similar to what we do) \cite{Jan}, an appropriate scaled limit (in which the oscillations vanish) \cite{LionsPaul,GMMP}, or the introduction of a stochastic averaging \cite{Ryz}. For a more complete discussion of the WT's interference terms and their interpretation, the interested reader can see \cite{Hla,Fl1}.

\begin{definition}[Smoothing in phase-space]
\label{DefSmoo2}
Denote by $\Phi$ the operator
\begin{equation}
\label{PHI2}
\begin{array}{c}
\Phi: w(x,k)\mapsto \mathcal{F}^{-1}_{X,K \shortrightarrow x,k} \left[{e^{-\frac{\pi}{2}\left[{ \sigma_x^2 X^2+\sigma_k^2 K^2 }\right] }  \mathcal{F}_{a,b \shortrightarrow X,K} \left[{ w(a,b) }\right] }\right]=\\
=\frac{2^n}{\sigma_x^n \sigma_k^n} \int\limits_{x \in \mathbb{R}^n} { e^{-2\pi \frac{(x-x')^2}{\sigma_x^2}-2\pi \frac{(k-k')^2}{\sigma_k^2}} w(x',k')dx'dk'}.
\end{array}
\end{equation}
\end{definition}

We use the same symbol as in Definition \ref{DefSmoo1}, although technically it is a different operator. Still, we will go on with this abuse of notation, because they are essentially very similar operators, and it is very easy to understand which one is used from the context: the one of Definition \ref{DefSmoo1} acts on functions of $x \in \mathbb{R}^n$, while the one of Definition \ref{DefSmoo2} acts on functions of $(x,k) \in \mathbb{R}^{2n}$.

\begin{definition}[The smoothed Wigner transform]
The SWT is the sesquilinear transform
\begin{equation}
\label{eq24}
\begin{array}{c}
\tilde{W}: f,g \mapsto
\tilde{W}[f,g](x,k)=\\
= \left({ \frac{\sqrt{2}}{\sigma_x} }\right)^n \int\limits_{u,y \in \mathbb{R}^n} {
e^{ -2 \pi i k y - \frac{\pi \sigma_k^2 y^2}{2} -\frac{2 \pi (u-x)^2}{\sigma_x^2} } f(u+\frac{y}{2}) \bar{g}(u-\frac{y}{2}) du dy }=\\
=\Phi W[f,g](x,k).
\end{array}
\end{equation}
The generalization for vectors is the same as for the WT. 

Moreover, it is well defined as a bilinear mapping
\begin{eqnarray}
\tilde{W}&:&\mathcal{S}(\mathbb{R}^n) \times \mathcal{S}(\mathbb{R}^n) \shortrightarrow 
G^{\frac{1}{2},(\sigma_x^2,\sigma_k^2)} (\mathbb{R}^{2n}) \subseteq 
\mathcal{S}^{\frac{1}{2},\frac{1}{\sqrt{\pi e}}(\frac{1}{\sigma_x},\frac{1}{\sigma_k})}(\mathbb{R}^{2n}) \subseteq 
\mathcal{S}(\mathbb{R}^{2n}), \\
\tilde{W}&:&L^2(\mathbb{R}^n) \times L^2(\mathbb{R}^n) \shortrightarrow L^2(\mathbb{R}^{2n}).
\end{eqnarray}
 \end{definition}

Of course we will be working a lot with the (quadratic) {\em smoothed Wigner distribution} $\tilde{W}[u](x,k)$.
The parameters $\sigma_x, \sigma_k$ control the length scales of the smoothing. The motivation is to smooth out any oscillations at length-scales finer than the oscillations of $u(x)$ itself ($\sigma_x^2$ is scaled with them) or those of $\hat{u}(k)$ ($\sigma_k^2$ is scaled with these). It should be mentioned here that if
\begin{equation}
\label{eq25}
\sigma_x \sigma_k =1
\end{equation}
then
\begin{equation}
\label{eq26}
\tilde{W}[u](x,k)= \left({ \frac{2}{ \sigma^2_x} }\right)^{\frac{n}{2}}  \vert \int\limits_{y \in \mathbb{R}^n} { e^{-2\pi i y k -\frac{\pi}{\sigma_x^2}(x-y)^2 }u(y)dy } \vert^2,
\end{equation}
i.e. $\tilde{W}$ coincides with a spectrogram (also known as Husimi transform) with window $g(y) = \left({ \frac{2}{\sigma_x^2} }\right)^ {\frac{n}{4}} e^{-\frac{\pi}{\sigma_x^2}y^2}$, and is therefore nonnegative. We will say that when $\sigma_x \sigma_k =1$ we have critical smoothing, while if $\sigma_x \sigma_k < 1$ the smoothing is sub-critical. Generally speaking, critical smoothing is pretty strong, and over-critical choices $\sigma_x \sigma_k > 1$ are not interesting. So $\sigma_x^2, \sigma_k^2$ are measured in units of typical wavelengths of $u(x)$, $\hat{u}(k)$ respectively, and the strength of the smoothing is gauged by the number $\sigma_x\sigma_k \in \left({ 0,1 }\right]$.  This automatically puts some structure in the parameter space, which is found to be sufficient in many practical applications -- although clearly there is room for more quantitative results in this respect. For more discussion and examples on the calibration of the smoothing see \cite{Ath,Ath0}.

Before we go on to the smoothed Wigner calculus, let us formulate, in our notation, the Wigner calculus:

\begin{lemma}[Wigner calculus]
\label{lmm6}
Let $f(x),g(x) \in \mathcal{S}(\mathbb{R}^n)$, $L(x,k) \in \mathcal{S}'(\mathbb{R}^{2n})$ and $L$ the pseudodifferential operator with $L(x,k)$ as its Weyl symbol. Then
\begin{equation}
\label{eq93}
W[Lf,g](x,k)=\mathcal{L} W[f,g](x,k)
\end{equation}
where
\begin{equation}
\label{eq94}
\begin{array}{l}
\mathcal{L} w(x,k)=\\
=\int\limits_{a,b} {
e^{2\pi i \left[{ ak+bx }\right]} \hat{L}(a,b) w\left({ x+\frac{a}{2},k-\frac{b}{2} }\right)dadb}=\\
=2^{2n} \mathcal{F}^{-1}_{X,K \shortrightarrow x,k} \left[{ 
\int\limits_{S,T} {
e^{ 2\pi i \left[{ S(x-K)+T(k+X) }\right] }\hat{L}(2S,2T) \hat{w}(X-S,K-T) dSdT} }\right],
\end{array}
\end{equation}
or, equivalently, its Weyl symbol is
\begin{equation}
\label{eq94a}
\mathcal{L}(x,k,X,K)=L(x-\frac{K}{2},k+\frac{X}{2})
\end{equation}
\end{lemma}

We state Lemma \ref{lmm6} for completeness and motivation; it is a standard result, and the proof is also contained as a special case of Theorem \ref{thrm4}.

\begin{theorem}[Smoothed Wigner calculus]
\label{thrm4}
Let $f(x),g(x) \in \mathcal{S}(\mathbb{R}^n)$, $L(x,k) \in \mathcal{S}'(\mathbb{R}^{2n})$ and $L$ the operator with Weyl symbol $L(x,k)$. Then
\begin{equation}
\label{eq70}
\tilde{W}[Lf,g](x,k)=\tilde{\mathcal{L}} \tilde{W}[f,g](x,k)
\end{equation}
where
\begin{equation}
\label{eq71}
\begin{array}{l}
\tilde{\mathcal{L}} w(x,k)=\\
=2^{2n}
\int\limits_{X,K,S,T} {
e^{2 \pi i \left[{ S(x-K+i \sigma_x^2X) + T(k+X+i\sigma_k^2 K) +xX+kK }\right] }
\hat{L}(2S,2T) \hat{w}(X-S,K-T) dSdTdXdK
},
\end{array}
\end{equation}
or
\begin{equation}
\begin{array}{l}
\label{eq71a}
\tilde{\mathcal{L}} w(x,k)= \\
=\int\limits_{S,T \in \mathbb{R}^n} {
\hat{L}(S,T) 
e^{2\pi i (Sx+Tk)-\frac{\pi}{2}\left({ \sigma_x^2 S^2 +\sigma^2_k T^2 }\right)} w(x+\frac{T+i\sigma_x^2 S}{2},k-\frac{S-i\sigma_k^2 T}{2})dSdT}.
\end{array}
\end{equation}
\end{theorem}

\noindent {\bf Remark:} {\em First of all, let us remark that Theorem \ref{thrm4} can be seen as an application of Theorem \ref{thrm3}, using the Wigner calculus $W[Lf,g](x,k)=\mathcal{L}W[f,g](x,k)$. In that connection, $\tilde{\mathcal{L}}=\Phi \mathcal{L} \Phi^{-1}$. However, we will prove Theorem \ref{thrm4} similarly to, but nevertheless independently from Theorem \ref{thrm3}; one reason is that computations which are anyway necessary when working with SWTs will be carried out in the process.} \\ \vspace{1mm}

\noindent {\bf Proof:} First, we will see that the operator $\tilde{\mathcal{L}}$ is well-defined on $G^{\frac{1}{2},(\sigma_x^2,\sigma_k^2)} (\mathbb{R}^{2n})$ functions for each of the formulations of equations (\ref{eq71}), (\ref{eq71a}).
As we saw in Theorem \ref{cor1}, 
\begin{equation}
\label{eq73}
\begin{array}{l}
g_1(S,T)=\\
=e^{2\pi i (Sx+Tk)-\frac{\pi}{2}\left({ \sigma_x^2 S^2 +\sigma^2_k T^2 }\right)} w(x+\frac{T+i\sigma_x^2 S}{2},k-\frac{S-i\sigma_k^2 T}{2}) \in \mathcal{S}(\mathbb{R}^{2n})
\end{array}
\end{equation}
as a function of $(S,T)$, and therefore equation (\ref{eq71a}) makes sense. 

To see that equation (\ref{eq71}) is well defined we have to demonstrate that
\begin{equation}
\label{eq87}
\begin{array}{l}
g(S,T)= \\
=\int\limits_{X,K} {
e^{2 \pi i \left[{ S(x-K+i \sigma_x^2X) + T(k+X+i\sigma_k^2 K) +xX+kK }\right] }
\hat{w}(X-S,K-T) dXdK
} \in \mathcal{S}(\mathbb{R}^{2n}).
\end{array}
\end{equation}
But
\begin{equation}
\label{eq88}
\begin{array}{l}
g(S,T)= \\
=\int\limits_{X,K} {
e^{2 \pi i \left[{ S(x-K+i \sigma_x^2X) + T(k+X+i\sigma_k^2 K) +xX+kK }\right] }
\hat{w}(X-S,K-T) dXdK}=\\ { } \\
=e^{2\pi i \left[{ 2Sx+2Tk }\right] -2\pi \left[{ \sigma_x^2 S^2 +\sigma_k^2 T^2 }\right]} \\
\,\,\,\,\,\,\,\,\,\,\,\,\,\,\,\,\,\,\,\,\,\,\,\,\,\,\,\,\,\,
\int\limits_{X,K} {
e^{2\pi i \left[{ X(x+T+i\sigma_x^2 S) + K(k-S+i\sigma_k^2 T) }\right]} \hat{w}(X,K)dXdK }=\\ { } \\
=e^{2\pi i \left[{ 2Sx+2Tk }\right] -2\pi \left[{ \sigma_x^2 S^2 +\sigma_k^2 T^2 }\right]}
w(x+T+i\sigma_x^2 S,k-S+i\sigma_k^2 T).
\end{array}
\end{equation}
The last equality makes use of Lemma \ref{lmm4}. The end result is a Schwartz test-function according to Theorem \ref{cor1}.

We will show that equations (\ref{eq71}), (\ref{eq71a}) are equivalent. Indeed, the passage from equation (\ref{eq71}) to (\ref{eq71a}) is essentially demonstrated in equation (\ref{eq88}). Observe that we only do a change of variables and a Fourier transform, so the reverse course follows as well. 

So we checked that all the formulations in the statement make sense and are equivalent. Now we will finally show that they give the smoothed Wigner calculus, i.e. equation (\ref{eq70}) holds. Like earlier, the way to check it is by showing
\begin{equation}
\label{eq90}
\Phi \mathcal{L} v(x,k) = \tilde{\mathcal{L}} \Phi v(x,k),
\end{equation}
where $w(x,k)=\Phi v(x,k)$, i.e. $v(x,k)=W[f,g](x,k)$. The {\em lhs} of equation (\ref{eq90}) is equal to
\begin{small}
\begin{equation}
\label{eq91}
\begin{array}{l}
\Phi \mathcal{L} v(a,b)= \Phi L(x-\frac{\partial_k}{4\pi i},k+\frac{\partial_x}{4\pi i}) v(a,b)=\\ { } \\
=\mathcal{F}^{-1}_{A,B \shortrightarrow a,b} \left[{
e^{-\frac{\pi}{2} \left({ \sigma_x^2 A^2 + \sigma^2_k B^2 }\right) } \mathcal{F}_{x,k \shortrightarrow A,B} \left[{
2^{2n} \mathcal{F}^{-1}_{X,K \shortrightarrow x,k} \left[{  { }^{{  }^{ }} }\right. }\right. }\right.  \\
\,\,\,\,\,\,\,\,\,\,\,\,\, \left.{ \left.{ \left.{
\int\limits_{S,T} {
e^{ 2\pi i \left[{ S(x-K)+T(k+X) }\right] }\hat{L}(2S,2T) \hat{v}(X-S,K-T) dSdT
} }\right] }\right] }\right] =\\ { } \\
=2^{2n} \int{
e^{ 2\pi i \left[{ S(x-K)+T(k+X)+xX+kK-xA-kB+aA+bB }\right] -\frac{\pi}{2} \left({ \sigma_x^2 A^2 + \sigma^2_k B^2 }\right)  }
} \\
\,\,\,\,\,\,\,\,\,\,\,\,\,\,\,\,\,\,\,\,\,\,\,\,\,\,\,\,\,\,\,\,\,\,\,\,\,\,\,\,
\hat{L}(2S,2T) \hat{v}(X-S,K-T) dSdTdXdKdxdkdAdB= \\ { } \\
=2^{2n} \int{
e^{ 2\pi i \left[{ -SK+TX+x(X+S-A)+k(K+T-B) +aA+bB }\right] -\frac{\pi}{2} \left({ \sigma_x^2 A^2 + \sigma^2_k B^2 }\right)  }dxdk
} \\
\,\,\,\,\,\,\,\,\,\,\,\,\,\,\,\,\,\,\,\,\,\,\,\,\,\,\,\,\,\,\,\,\,\,\,\,\,\,\,\,
\hat{L}(2S,2T) \hat{v}(X-S,K-T) dSdTdXdKdAdB= \\ { } \\
=2^{2n} \int{
\delta(X+S-A) \delta(K+T-B)
e^{ 2\pi i \left[{ -SK+TX+aA+bB }\right] -\frac{\pi}{2} \left({ \sigma_x^2 A^2 + \sigma^2_k B^2 }\right)  }dAdB
} \\
\,\,\,\,\,\,\,\,\,\,\,\,\,\,\,\,\,\,\,\,\,\,\,\,\,\,\,\,\,\,\,\,\,\,\,\,\,\,\,\,
\hat{L}(2S,2T) \hat{v}(X-S,K-T) dSdTdXdK= \\ { } \\
=2^{2n} \int{
e^{ 2\pi i \left[{ -SK+TX+a(X+S)+b(K+T) }\right] -\frac{\pi}{2} \left({ \sigma_x^2 (X+S)^2 + \sigma^2_k (K+T)^2 }\right)  }
} \\
\,\,\,\,\,\,\,\,\,\,\,\,\,\,\,\,\,\,\,\,\,\,\,\,\,\,\,\,\,\,\,\,\,\,\,\,\,\,\,\,
\hat{L}(2S,2T) \hat{v}(X-S,K-T) dSdTdXdK.
\end{array}
\end{equation}
\end{small}
The {\em rhs} of equation (\ref{eq90}) is equal to
\begin{small}
\begin{equation}
\label{eq92}
\begin{array}{l}
\tilde{\mathcal{L}} \Phi v(a,b)=\\ { } \\
=2^{2n} \int{
e^{2\pi i \left[{ S(a-K+i\sigma_x^2 X)+T(b+X+i\sigma_k^2 K)+aX+bK }\right]} \hat{L}(2S,2T)}\\
\,\,\,\,\,\,\,\,\,\,\,\,\,\,\,\,\,\,\,\,\,\,\,\,\,\,\,\,\,\,\,\,\,\,\,\,\,\,\,\,
e^{-\frac{\pi}{2}\left({
\sigma_x^2(X-S)^2 + \sigma_k^2(K-T)^2 }\right)} \hat{v}(X-S,K-T)dXdKdSdT=\\ { } \\
=2^{2n} \int{
e^{2\pi i \left[{ S(a-K)+T(b+X)+aX+bK }\right]} \hat{L}(2S,2T)}\\
\,\,\,\,\,\,\,\,\,\,\,\,\,\,\,\,\,\,\,\,\,\,\,\,\,\,\,\,\,\,\,\,\,\,\,\,\,\,\,\,
e^{-\frac{\pi}{2}\left({
\sigma_x^2(X+S)^2 + \sigma_k^2(K+T)^2 }\right)} \hat{v}(X-S,K-T)dXdKdSdT=\\ { } \\
=2^{2n} \int{
e^{2\pi i \left[{ -SK+TX+a(X+S)+b(K+T) }\right]} \hat{L}(2S,2T)}\\
\,\,\,\,\,\,\,\,\,\,\,\,\,\,\,\,\,\,\,\,\,\,\,\,\,\,\,\,\,\,\,\,\,\,\,\,\,\,\,\,
e^{-\frac{\pi}{2}\left({
\sigma_x^2(X+S)^2 + \sigma_k^2(K+T)^2 }\right)} \hat{v}(X-S,K-T)dXdKdSdT,
\end{array}
\end{equation}
\end{small}
which is the same as the last member of equation (\ref{eq91}).

The proof is complete. \\

\begin{theorem}[Weyl symbols for the smoothed Wigner calculus]
\label{thrm4b}
Consider $f(x),g(x) \in \mathcal{S}(\mathbb{R}^n)$, and $L(x,k)$ to be the Fourier transform of a compactly supported tempered distribution (in particular it is  the restriction to the real numbers of) an entire analytic function.
Then the operator $\tilde{\mathcal{L}}$, defined in equation (\ref{eq70}), has Weyl symbol
\begin{equation}
\label{eq72}
\begin{array}{c}
\tilde{\mathcal{L}}(x,k,X,K)
=L\left({
x-\frac{ K-i\sigma_x^2X}{2} , k+\frac{X+i\sigma_k^2K}{2}  }\right).
\end{array}
\end{equation}
\end{theorem}

\noindent {\bf Proof: } This proof follows along the exact same lines as the proof of theorem \ref{thrm3b}.

Indeed, observe that
\begin{equation}
\label{eq89}
\begin{array}{l}
L\left({
x-\frac{ \partial_k-i\sigma_x^2 \partial_x}{4\pi i} , k+\frac{\partial_x+i\sigma_k^2 \partial_k}{4\pi i}  }\right) w(x,k)=\\ { } \\
=\int{ e^{2\pi i \left[{ xX+kK-aX-bK }\right]}
L\left({
x-\frac{ K-i\sigma_x^2 X}{2} , k+\frac{X+i\sigma_k^2 K}{2}  }\right)  w(a,b) dadbdXdK }=\\ { } \\
=\int{
e^{2\pi i \left[{ S\frac{x+a-K+i\sigma_x^2 X}{2}+T\frac{k+b+X+i\sigma_k^2 K}{2} -sS-tT }\right]} L(s,t)dsdtdSdT}\\
\,\,\,\,\,\,\,\,\,\,\,\,\,\,\,\,\,\,\,\,\,\,\,\,\,\,\,\,\,\,\,\,\,\,\,\,\,\,\,\,\,\,\,\,\,\,\,\,\,\,\,\,\,\,\,\,\,\,\,\,\,\,\,\,\,\,\,\,\,\,\,\,\,\,\,\,\,\,\,\,
{e^{2\pi i \left[{ xX+kK-aX-bK }\right]} w(a,b) dadbdXdK}=\\ { } \\
=\int{
e^{2\pi i \left[{
X \left({ x+\frac{T}{2}+i\sigma_x^2\frac{S}{2} }\right) +K \left({ k-\frac{S}{2}+i\sigma_k^2 \frac{T}{2} }\right)+
S\frac{x}{2}+T\frac{k}{2}-a\left({X-\frac{S}{2}}\right)-b \left({ K-\frac{T}{2} }\right)
}\right]} } \\
\,\,\,\,\,\,\,\,\,\,\,\,\,\,\,\,\,\,\,\,\,\,\,\,\,\,\,\,\,\,\,\,\,\,\,\,\,\,\,\,\,\,\,\,\,\,\,\,\,\,\,\,\,\,\,\,\,\,\,\,\,\,\,\,\,\,\,\,\,\,\,\,\,\,\,\,\,\,\,\, {\hat{L}(S,T)w(a,b)dSdTdadbdXdK}=\\ { } \\
=\int{
e^{2\pi i \left[{
X \left({ x+\frac{T}{2}+i\sigma_x^2\frac{S}{2} }\right) +K \left({ k-\frac{S}{2}+i\sigma_k^2 \frac{T}{2} }\right)+
S\frac{x}{2}+T\frac{k}{2}
}\right]} } \\
\,\,\,\,\,\,\,\,\,\,\,\,\,\,\,\,\,\,\,\,\,\,\,\,\,\,\,\,\,\,\,\,\,\,\,\,\,\,\,\,\,\,\,\,\,\,\,\,\,\,\,\,\,\,\,\,\,\,\,\,\,\,\,\,\,\,\,\,\,\,\,\,\,\,\,\,\,\,\,\, {\hat{L}(S,T)\hat{w}\left({X-\frac{S}{2}, K-\frac{T}{2} }\right)dSdTXdK}=\\ { } \\
=2^{2n} \int{
e^{2\pi i \left[{ X(x+T+i\sigma_x^2 S)+K(k-S+i\sigma_k^2 T) +xS+kT }\right]} \hat{L}(2S,2T)} \\
\,\,\,\,\,\,\,\,\,\,\,\,\,\,\,\,\,\,\,\,\,\,\,\,\,\,\,\,\,\,\,\,\,\,\,\,\,\,\,\,\,\,\,\,\,\,\,\,\,\,\,\,\,\,\,\,\,\,\,\,\,\,\,\,\,\,\,\,\,\,\,\,\,\,\,\,\,\,\,\, {\hat{w}(X-S,K-T)dXdKdSdT}=\\ { } \\
=2^{2n}
\int{
e^{2 \pi i \left[{ S(x-K+i \sigma_x^2X) + T(k+X+i\sigma_k^2 K) +xX+kK }\right] }
\hat{L}(2S,2T)} \\
\,\,\,\,\,\,\,\,\,\,\,\,\,\,\,\,\,\,\,\,\,\,\,\,\,\,\,\,\,\,\,\,\,\,\,\,\,\,\,\,\,\,\,\,\,\,\,\,\,\,\,\,\,\,\,\,\,\,\,\,\,\,\,\,\,\,\,\,\,\,\,\,\,\,\,\,\,\,\,\, { \hat{w}(X-S,K-T) dSdTdXdK
},
\end{array}
\end{equation}
which is exactly  the rhs of equation (\ref{eq71}). 

We have seen the justification of imaginary translations through Fourier transforms for entire functions before. Here it is applied on a smoothed function function and therefore Lemma \ref{lmm4} applies.

The other step that needs justification is the interchange of the $dadb$ and $dSdT$ integrations, passing from the third to the fourth line. Remember that we have assumed $\exists M>0$ such that $supp \, \hat{L}(S,T) \, \subseteq \, [-M,M]^{2n}$, therefore the real exponential terms $e^{-\pi \sigma_x^2 X S-\pi \sigma_k^2 K T}$ can be substituted by $e^{-\pi \sigma_x^2 X S-\pi \sigma_k^2 K T} \chi_{[-2M,2M]^{2n}}(S,T)$ without changing anything. The result then follows by the standard tempered distribution calculus.

The proof is complete. \\ \vspace{0.5mm}

\subsection{Coarse-scale dynamics in phase-space}
\label{secCSPS}

Theorem \ref{thrm4} allows us to carry out in a precise manner the basic steps of SWT homogenization outlined in Section \ref{sec1_2}: \\

\begin{corollary}[Smoothed Wigner equations]
\label{thrm14}
Let $L(x,k) \in \left({ \mathcal{S}'(\mathbb{R}^{2n})}\right)^{d \times d}$, $L=L(x,\partial_x)$, $u_0(x) \in \left({ \mathcal{S}(\mathbb{R}^n)}\right)^{d}$. Consider the IVP
\begin{equation}
\label{eq111}
\begin{array}{c}
u_t (x,t)+Lu(u,t)=0,\\
u(x,0)=u_0(x).
\end{array}
\end{equation}
Then the SWT of $u$,
\begin{equation}
\label{eq112}
\tilde{W}(x,k,t)=\tilde{W}[u(\cdot,t)](x,k)
\end{equation}
satisfies the IVP
\begin{equation}
\label{eq113}
\begin{array}{c}
\tilde{W}_t (x,t)+2 \mathcal{H} \left({ \tilde{\mathcal{L}} W(u,t) }\right)=0,\\
\tilde{W}(x,k,0)=\tilde{W}[u_0](x,k),
\end{array}
\end{equation}
where $\tilde{\mathcal{L}}$ is defined in terms of $L$ as in Theorem \ref{thrm4}, and $\mathcal{H}(A)=\frac{A+A^*}{2}$ denotes the Hermitian part of a matrix.
\end{corollary}

The proof is obvious, and consists in the application of Theorem \ref{thrm4}, and the observation that
\begin{equation}
\label{eq114}
\begin{array}{c}
\frac{\partial}{\partial t} \tilde{W}(x,k,t)=\tilde{W}[u_t,u](x,k,t)+\tilde{W}[u,u_t](x,k,t)=\\
=-\tilde{W}[Lu,u](x,k,t)-\tilde{W}[u,Lu](x,k,t)=\\
=-\tilde{\mathcal{L}} \tilde{W}[u](x,k,t)-\left({ \tilde{\mathcal{L}} \tilde{W}[u](x,k,t)}\right)^*=
\end{array}
\end{equation}

\begin{corollary}[Smoothed trace formula]
\label{thrm12}
Let $M_{i,j}(x,k) \in \mathcal{S}'(\mathbb{R}^{2n})$, $i,j=1,...,d$.
Any quadratic observable of a wavefield $u_i(x,t) \in \mathcal{S}(\mathbb{R}^n)$, $i=1,...,d$, corresponding to the operator $M=M \left(  {x,\partial _x  } \right)$ (defined as in equation (\ref{eq2}))
 can be directly expressed in terms of the (Hermitian-matrix-valued) smoothed Wigner distribution of $u$ as
\begin{equation}
\label{eq104}
\mathbb{M}(t)
=\int\limits_{x,k \in \mathbb{R}^n} { tr \left({ \tilde{\mathcal{M}} \tilde{W}[u](x,k,t) }\right) dxdk },
\end{equation}
where $\tilde{\mathcal{M}}_{i,j}$ is defined in terms of $M_{i,j}(x,k)$ in the same way as $\tilde{\mathcal{L}}$ is defined in terms of $L(x,k)$ in Theorem \ref{thrm4}.
Moreover, observables can be resolved over phase-space at coarse-scale, 
\begin{equation}
\label{eq105}
\tilde{\mathbb{M}}(x,k,t)=tr \left({ \tilde{\mathcal{M}} \tilde{W}(x,k,t)}\right)
\end{equation}
consistently with their natural resolutions,
\begin{equation}
\label{eq106}
\int\limits_{k\in \mathbb{R}^n} {\tilde{\mathbb{M}}(x,k,t)dk}=\frac{\sqrt{2}^n}{\sigma^n_x} \int\limits_{x' \in \mathbb{R}^n} { e^{-\frac{2\pi (x-x')^2}{\sigma_x^2}}  tr \left({ \bar{u}^T(x',t)M u(x',t) }\right) dx'}.
\end{equation}
and
\begin{equation}
\label{eq107}
\int\limits_{x\in \mathbb{R}^n} {\tilde{\mathbb{M}}(x,k,t)dx}= \frac{\sqrt{2}^n}{\sigma^n_k} \int\limits_{k' \in \mathbb{R}^n} { e^{-\frac{2\pi (k-k')^2}{\sigma_k^2}}  tr \left({ \bar{\hat{u}}^T(k',t) \widehat{M u}(k',t) }\right) dk'}.
\end{equation}
\end{corollary}

The proof is obvious, since
\begin{equation}
\label{eq108}
\tilde{\mathcal{M}} \tilde{W}[u](x,k,t)=\tilde{W}[Mu,u](x,k,t),
\end{equation}
and, $\forall u,v \in \left({ \mathcal{S}(\mathbb{R}^n) }\right)^d$
\begin{equation}
\label{eq109}
\int\limits_{x,k \in \mathbb{R}^n} {tr \left({\tilde{W}[u,v](x,k) }\right)dxdk} = 
\int\limits_{x \in \mathbb{R}^n} {\bar{v}^T(x) u(x) dx},
\end{equation}
\begin{equation}
\label{eq110}
\int\limits_{x\in \mathbb{R}^n} {\tilde{W}[u,v](x,k)dx}= \frac{\sqrt{2}^n}{\sigma^n_k} \int\limits_{k' \in \mathbb{R}^n} { e^{-\frac{2\pi (k-k')^2}{\sigma_k^2}}  tr \left({ \bar{\hat{v}}^T(k',t) \hat{ u}(k',t) }\right) dk'},
\end{equation}
and similarly for the $dk$ marginal. 

In particular, all observables corresponding to polynomial Weyl symbols (which typically include energy and energy flux) can be recovered from the SWT of the wavefunction in terms of finite-order operators. \\ \vspace{2mm}

Corollaries \ref{thrm14} and \ref{thrm12} show how the smoothed Wigner calculus allows us to reformulate problems, originally formulated for ``waves'' (i.e. for an oscillating wavefunction on $\mathbb{R}^n$), to problems for ``phase-space densities'' (i.e. smooth / simple functions on $\mathbb{R}^{2n}$). Indeed, many well known paradigms fit in this general description, with semiclassical limits and Wigner measures being the most relevant from a technical point of view \cite{Ben,Fil,Ge91,LionsPaul,MM1,GMMP,Wig}. The introduction of a ``fundamental length'', controlled by $\sigma_x^2$, $\sigma_k^2$ is a distinctively different feature from the Wigner measure approach; the concept of a fundamental length has also been discussed from a physical point of view as well. For more treatments that can also be described as ``phase-space homogenization'' in a wide sense -- but not as closely related to what we do here from a technical standpoint -- see also \cite{Ryz} and the relevant survey in the introduction therein. A different problem (although similar in the sense that it involves infinite-order equations governing a smooth density) is treated in \cite{Luc}. \\ \vspace{1mm}

Naturally, a concrete example is in place here: 

\begin{example}[Schr\"{o}dinger equation]
\label{ex1}
Consider a wavefunction satisfying the Schr\"{o}dinger equation,
\begin{equation}
\label{eq116}
\begin{array}{c}
\frac{\partial }{\partial t}u \left( {x,t} 
\right)-\frac{i}{2}\Delta u\left( {x,t} 
\right)+iV\left( x \right)u\left( {x,t} \right)=0, \\
u(x,0)=u_0(x).
\end{array}
\end{equation}
Then its SWT $\tilde{W}(x,k,t)=\tilde{W}[u(\cdot,t)](x,k)$ satisfies the equation
\begin{equation}
\label{eq116a}
\begin{array}{c}
\frac{\partial }{\partial t}  \tilde{W}(x,k,t) +
\left(  {2\pi k \cdot \nabla_x+\frac{\sigma _k^2 
}{2} \nabla_x \cdot \nabla_k }\right)\tilde{W}\left( 
{x,k,t} \right)+\\ 
+2 Re \left({i
\int\limits_{s \in \mathbb{R}^n} {
e^{2\pi i s \left({x+\frac{i \sigma_x^2 }{4}s}\right)} \hat{V}(s) \tilde{W}(x+\frac{i\sigma_x^2 s}{2},k-\frac{s}{2})ds
} }\right), \\ { } \\
\tilde{W}(x,k,0)=\tilde{W}[u_0](x,k).
\end{array}
\end{equation}
\end{example}

The algebra of the smoothed calculus yields a very fortunate ``accident'': we can have exact coarse-scale reformulations for certain nonlinear equations with no additional work:

\begin{example}[Cubic non-linear Schr\"odinger equation]
\label{ex2}
Consider a wavefunction satisfying the cubic NLS equation,
\begin{equation}
\label{eq180}
\begin{array}{c}
\frac{\partial }{\partial t}u \left( {x,t} 
\right)-\frac{i}{2}\Delta u\left( {x,t} 
\right)+i \left( V_1\left( x \right) + \beta |u(x,t)|^2 \right) u\left( {x,t} \right)=0, \\ 
u(x,0)=u_0(x).
\end{array}
\end{equation}
Then its SWT $\tilde{W}(x,k,t)=\tilde{W}[u(\cdot,t)](x,k)$ satisfies the equation
\begin{equation}
\label{eq181}
\begin{array}{c}
\frac{\partial }{\partial t}  \tilde{W}(x,k,t) +
\left(  {2\pi k \cdot \nabla_x+\frac{\sigma _k^2 
}{2} \nabla_x \cdot \nabla_k }\right)\tilde{W}\left( 
{x,k,t} \right)+\\ 
+2 Re \left({i
\int\limits_{s \in \mathbb{R}^n} {
e^{2\pi i s x} \hat{F}(s,t) \tilde{W}(x+\frac{i\sigma_x^2 s}{2},k-\frac{s}{2})ds
} }\right), \\ { } \\
\tilde{W}(x,k,0)=\tilde{W}[u_0](x,k),
\end{array}
\end{equation}
where
\begin{equation}
\begin{array}{c}
\label{eq182}
\hat{F}(s,t)=\mathcal{F}_{x \shortrightarrow s} \left[ F(x,t) \right], \\ { } \\
F(x,t)=\frac{\sqrt{2}^n}{\sigma_x^n} \int\limits_{x' \in \mathbb{R}^n} {e^{-\frac{2\pi (x-x')^2}{\sigma_x^2}} V_1(x')dx'}+ \beta \int\limits_{k \in \mathbb{R}^n} {\tilde{W}(x,k,t)dk }.
\end{array}
\end{equation}
\end{example}

\begin{example}[Hartree equation (to smoothed Vlasov)]
\label{ex3}
Consider a wavefunction satisfying the Hartree equation,
\begin{equation}
\label{eq180v}
\begin{array}{c}
\frac{\partial }{\partial t}u \left( {x,t} 
\right)-\frac{i}{2}\Delta u\left( {x,t} 
\right)+i \left( V_1\left( x \right) + \int\limits_{x' \in \mathbb{R}^n}{V_0(x-x') |u(x',t)|^2 dx'} \right) u\left( {x,t} \right)=0, \\ 
u(x,0)=u_0(x).
\end{array}
\end{equation}
Then its SWT $\tilde{W}(x,k,t)=\tilde{W}[u(\cdot,t)](x,k)$ satisfies the equation
\begin{equation}
\label{eq181v}
\begin{array}{c}
\frac{\partial }{\partial t}  \tilde{W}(x,k,t) +
\left(  {2\pi k \cdot \nabla_x+\frac{\sigma _k^2 
}{2} \nabla_x \cdot \nabla_k }\right)\tilde{W}\left( 
{x,k,t} \right)+\\ 
+2 Re \left({i
\int\limits_{s \in \mathbb{R}^n} {
e^{2\pi i s x} \hat{F}(s,t) \tilde{W}(x+\frac{i\sigma_x^2 s}{2},k-\frac{s}{2})ds
} }\right), \\ { } \\
\tilde{W}(x,k,0)=\tilde{W}[u_0](x,k),
\end{array}
\end{equation}
where
\begin{equation}
\begin{array}{c}
\label{eq182v}
\hat{F}(s,t)=\mathcal{F}_{x \shortrightarrow s} \left[ F(x,t) \right], \\ { } \\
F(x,t)=\frac{\sqrt{2}^n}{\sigma_x^n} \int\limits_{x' \in \mathbb{R}^n} {e^{-\frac{2\pi (x-x')^2}{\sigma_x^2}} V_1(x')dx'}+ \beta \int\limits_{k,x' \in \mathbb{R}^n} { V_0(x-x') \tilde{W}(x',k,t)dkdx'}.
\end{array}
\end{equation}
\end{example}

The derivation for either nonlinear equation follows by observing simply that the potential appears in equation (\ref{eq116a}) not just as $\hat{V}(s)$, but as \[
e^{-\frac{\pi}{2} \sigma_x^2 s^2} \hat{V}(s)= 
\mathcal{F}_{x \shortrightarrow s} \left[ \tilde{V}(x) \right]=
\mathcal{F}_{x \shortrightarrow s} \left[ \frac{\sqrt{2}^n}{\sigma_x^n} \int\limits_{x' \in \mathbb{R}^n} {e^{-\frac{2\pi (x-x')^2}{\sigma_x^2}} V(x')dx'} \right],
\]
i.e. instead of the original potential $V(x)$, it suffices to know the smoothed potential, $\tilde{V}(x)$. This, coupled with the marginals property of the SWT -- equation (\ref{eq106}) for $M=I$ -- makes it possible to have closed smoothed Wigner equations in this case. \\ \vspace{2mm}

As we mentioned earlier, the closest relative to this approach (and an important motivation for it) is the WT / Wigner measure based semiclassical limits technique. In the next Section we study an application of the SWT to semiclassical problems.

\section{The semiclassical regime}
\label{secSemi}
Let us start with a  few words of motivation. We will work in an asymptotic regime, scaled with a parameter $0<\varepsilon<<1$. The intuitive meaning of the small parameter $\varepsilon$ is that we work with signals / functions which exhibit very fast oscillations, e.g. WKB functions
\begin{equation}
\label{eq175}
f^\varepsilon (x)=A(x) e^{\frac{2\pi i}{\varepsilon} S(x)}.
\end{equation}

Under certain conditions (e.g. the envelope $A(x)$ is itself ``smooth'', ``slowly varying'') it can be said that the function of equation (\ref{eq175}) has amplitude $n(x) \approx |A(x)|^2$ and ``instantaneous frequency''/``local wavenumber'' $k(x) \approx \nabla S(x)$. Indeed abstractions like these -- and making them precise -- are at the heart of the (motivation for the) WT and time-frequency analysis. It must be clear already why it is natural that these questions are formulated in an asymptotic regime $\varepsilon <<1$. \footnote{This could be seen as a ``signal-processing-inspired'' introduction for the semiclassical regime, see also \cite{Fl1}. The semiclassical regime, as the name shows, can also be seen as a physical regime of ``large'' quantum systems, as was the original motivation of Wigner \cite{Wig,MM1}.  }

\begin{definition}[Semiclassical scaling of the WT]
\label{def7}
The semiclassically scaled WT is defined as
\begin{equation}
\label{eq125A}
\begin{array}{c}
W^\varepsilon: f,g \mapsto
W^\varepsilon[f,g](x,k)=\\
=\int\limits_{y \in \mathbb{R}^n}
{e^{-2 \pi i k y} f\left({ x+\frac{\varepsilon y}{2} }\right) \bar{g}\left({ x-\frac{\varepsilon y}{2} }\right) dy}=
\frac{1}{\varepsilon^n} W[f,g](\frac{k}{\varepsilon}).
\end{array}
\end{equation}
in agreement to \cite{LionsPaul,GMMP,Wig}.
\end{definition}.

\begin{definition}[Semiclassical scaling of the SWT]
The semiclassically scaled SWT is defined as
\begin{equation}
\label{eq124A}
\begin{array}{c}
\tilde{W}^\varepsilon: f,g \mapsto
\tilde{W}^\varepsilon[f,g](x,k)=\\
=\left({ \frac{\sqrt{2}}{\sigma_x \sqrt{\varepsilon}} }\right)^n \int\limits_{u,y \in \mathbb{R}^n} {
e^{ -2 \pi i k y - \frac{\pi \varepsilon \sigma_k^2 y^2}{2} -\frac{2 \pi (u-x)^2}{ \varepsilon \sigma_x^2} } f(u+\frac{\varepsilon y}{2}) \bar{g}(u-\frac{\varepsilon y}{2}) du dy }=\\
=\frac{1}{\varepsilon^n} \Phi_{ \sqrt{\varepsilon} \sigma_x ,\sqrt{\varepsilon} \sigma_k} W[f,g] \left( x,\frac{k}{\varepsilon} \right)=
\Phi_{ \sqrt{\varepsilon} \sigma_x ,\sqrt{\varepsilon} \sigma_k} W^\varepsilon [f,g] \left( x,k \right).
\end{array}
\end{equation}
\end{definition}

Obviously the choice of the scaling of the smoothing,
\begin{equation}
\label{eq126A}
\sigma_x, \, \sigma_k \,\, \mapsto \sqrt{\varepsilon} \sigma_x, \, \sqrt{\varepsilon} \sigma_k,
\end{equation}
is to some extent arbitrary; in what follows we hope to show that it is a natural choice, at least for some problems. \footnote{Numerical examples also offer important insights in this question. Of course, in certain problems it might be that some other scaling is better. We only propose this as a reasonable, general-purpose starting point.}

A central object in semiclassical problems is the well studied Wigner measure (WM). We recall (adapted to our notation) a well-known and central result (see e.g. Proposition 1.1 and Remark 1.3 in \cite{GMMP}) \, :

\begin{theorem}[Definition of the WM]
\label{thrm18}
Consider a ``semiclassical family of functions'' $\lbrace f^\varepsilon(x) \rbrace_{\varepsilon \in (0,1)}$ satisfying the condition
\begin{equation}
\label{eq176}
\exists M_0 >0 \,\,:\,\, || f^\varepsilon||_{L^2 (\mathbb{R}^n)} \leqslant M_0.
\end{equation}
Then the family of the semiclassical WTs $\lbrace W^\varepsilon [f^\varepsilon](x,k) \rbrace_{\varepsilon \in (0,1)}$ has weak-$*$ accumulation points as a set of functionals on an appropriate space of test-functions on phase space. When the accumulation point is unique (equivalently, up to the extraction of a subsequence) it will be called the WM associated with the semiclassical family  $\lbrace f^\varepsilon(x) \rbrace_{\varepsilon \in (0,1)}$.
\end{theorem}

Usually we will consider families with a unique accumulation point
\begin{equation}
\label{eq177}
W^\varepsilon [f^\varepsilon](x,k) \rightarrow W^0(x,k).
\end{equation}
An important example is given for WKB families $f^\varepsilon (x)=A(x) e^{\frac{2\pi i}{\varepsilon} S(x)}$, where (under appropriate auxiliary assumptions, see e.g. \cite{Fil})
\begin{equation}
\label{eq178}
W^\varepsilon [f^\varepsilon](x,k) \rightarrow |A(x)|^2 \delta\left( k-\nabla S(x) \right).
\end{equation}

The WM (i.e. Theorem \ref{thrm18}, but also concrete examples such as that of equation (\ref{eq178}) above) is the true justification for the semiclassical scaling of the WT.

A very successful technique in semiclassical limits has been to use the WM $W^0 (x,k)$ to keep track of the ``data of the problem'' (i.e. an appropriate family of observables) in an asymptotic problem, i.e. for $\varepsilon << 1$. This has been successful in many cases \cite{LionsPaul,GMMP} to name but some landmark works. However, this approach has its own limitations, such as leading to inconsistent / ill-posed problems in some cases -- see \cite{CFMS} for a recent survey. We will also point out a couple of other issues here -- which exist even when the WM based model can be formulated and is well-posed:
\begin{itemize}
\item
Due to the interference terms, the incorporation of $\varepsilon$-dependent corrections to a WM based model is virtually impossible. This introduces a rigid scheme of the information that can be kept track of or not. (Indeed one might say that morally, this is at the root of some of the problems surveyed in \cite{CFMS}). The SWT offers, as we wish to show, a more flexible way to decide ``how much detail to keep''. Some more quantitative results in this direction have also been presented in \cite{Ath,Ath0}.
\item
The fact that we have to work with a singular object (a measure supported on a low-dimensional manifold) introduces many analytical as well as numerical nuances. We believe that regularizing to a ``nice'' smooth density has the potential to make many things easier, or even possible for the first time. A concrete, quantitative result in that direction is Theorem \ref{thrm15}.  \\ \vspace{2mm}
\end{itemize}

An important fact that we need to mention here; it is well known that in the case of critical smoothing, the WM is preserved. We quote the following result (in adapted notation) from \cite{MM1}:

\begin{theorem}[Husimi has the same weak limit as Wigner]
\label{thrm19}
Consider a semiclassical family $\lbrace f^\varepsilon \rbrace$ with WM
\begin{equation}
W^\varepsilon [f^\varepsilon](x,k) \,\,\, \rightarrow \,\,\, W^0(x,k)
\end{equation}
in $L^2(\mathbb{R}^{2n})$-weak.
Denote $g^\varepsilon(x)=2^\frac{n }{2 } e^{ -\frac{2\pi |x|^2}{\varepsilon} }$,
\begin{equation}
H^\varepsilon [f^\varepsilon](x,k)=\int\limits_{x', k' \in \mathbb{R}^n}
{ g^\varepsilon(x-x') g^\varepsilon(k-k') W^\varepsilon [f^\varepsilon] (x',k')dx'dk'}=\Phi_{\sqrt{\varepsilon},\sqrt{\varepsilon}} W^{\varepsilon}[f^\varepsilon](x,k).
\end{equation}
Then
\begin{equation}
H^\varepsilon [f^\varepsilon](x,k) \,\,\, \rightarrow \,\,\, W^0(x,k)
\end{equation}
in $L^2(\mathbb{R}^{2n})$-weak.
\end{theorem}

By a straightforward adaptation of the same proof as in \cite{MM1} the same can be seen to hold in the more general case of ($\varepsilon$-independent) $\sigma_x^2,\sigma_k^2$ as well, i.e. the SWT has the same weak limit as the WT,
\begin{equation}
\tilde{W}^\varepsilon [f^\varepsilon](x,k) \,\,\, \rightarrow \,\,\, W^0(x,k)
\end{equation}
in $L^2(\mathbb{R}^{2n})$-weak.

This is important, because it shows that working with the SWT is in fact not a different approach than the WM. As soon as we take $\varepsilon \shortrightarrow 0$, working with the WT or with the SWT are indistinguishable. The difference we want to build on, is that the SWT behaves drastically better in several respects -- e.g. numerically -- than the WT in the regime $0< \varepsilon << 1$.  \\ \vspace{1mm}

A useful device in working with WTs in the semiclassical regime is the asymptotic computation of WTs of WKB functions. Indeed, computations of that kind are used in \cite{Fil,Fl1,Hla} to provide valuable insights -- a simple one being equation (\ref{eq178}). We carry out the respective computation for the SWT:

Consider a WKB function of the form
\begin{equation}
\label{eqTh1}
f^\varepsilon(x)=A(x)e^{2\pi i\frac{S(x}{\varepsilon} }.
\end{equation}
Let us suppose moreover that both $A$ and $S$ are analytic in a neighborhood of $q\in\mathbb R^n$.
We know that the SWT of $f^\varepsilon$ will be localized near $(q,p:=\nabla S(q))$. The following result makes this precise:
\begin{theorem}[Asymptotic computation of the SWT of a WKB function]
\label{thrm21}
There is a neighborhood $\Omega$ of the point $(q,p:=\nabla S(q))$ such that, $\forall\ (x,k)\in\Omega$,
\begin{equation}
\label{eqTh2}
\tilde{W}^\varepsilon[f^\varepsilon](x,k) =  \frac{|A(x)|^2} {\sqrt{\left(\frac{\varepsilon\sigma_k^2}2\right)^{2n}\mbox{Det}
\left(\frac{1+3(\frac{\sigma_x}{\sigma_k}D^2S(x))^2}{1+(\frac{\sigma_x}{\sigma_k}D^2S(x))^2}\right)}}
e^{-\frac {2\pi} {\varepsilon\sigma_k^2}(k-\nabla S(x))^T\frac{1+3(\frac{\sigma_x}{\sigma_k}D^2S(x))^2}{1+(\frac{\sigma_x}{\sigma_k}D^2S(x))^2}(k-\nabla S(x))}+O(\varepsilon^{\frac{1}2})
\end{equation}
\end{theorem}

Let us sketch the proof of the Theorem. Let us compute $\tilde{W}^\varepsilon[f^\varepsilon]$:
\begin{equation}
\tilde{W}^\varepsilon[f^\varepsilon](x,k)=\frac{2^n}{\sigma_x^n\sigma_k^n}\int A(x'+\delta)\overline{A(x'-\delta)}
e^{2\pi i\frac{S(x'+\delta)-S(x'-\delta)-2k'\delta}\varepsilon}e^{-2\pi\frac{(x'-x)^2}{\varepsilon\sigma_x^2}-
2\pi\frac{(k'-k)^2}{\varepsilon\sigma_k^2}}dx' dk' d\delta
\end{equation}
The stationary points of the phase are real only if $p=\nabla S(q)$. When $p\sim\nabla S(q)$, following the method of \cite{PU} and since $A$ and $S$ are analytic we can change the path 
of integration in order to catch the complex stationary points which are given by the equations:
\begin{eqnarray}
-2k'+\nabla S(x'+\delta)+\nabla S(x'-\delta)=0\\
i\nabla S(x'+\delta)-i\nabla S(x'-\delta)-2\frac{x'-x}{\sigma_x^2}=0\\
-2i\delta -2\frac{k'-k}{\sigma_k^2}=0
\end{eqnarray}
Let us compute everything for $k+\nabla S(x)$, $k-k'$ and $x-x'$ small.
We get
\begin{equation}
\delta=i\frac{k'-k}{\sigma_k^2}
\end{equation}
and
\begin{eqnarray}
-D^2S(x')\left(\frac{k'-k}{\sigma_k^2}\right)-\frac{x'-x}{\sigma_x^2}=0\\
-k+\nabla S(x)+D^2S(x)(x'-x)-(k'-k)=0
\end{eqnarray}
So
\begin{equation}
x'-x=-\frac{\sigma_x^2}{\sigma_k^2}D^2S(x)(k'-k)
\end{equation}
and
\begin{equation}
k'-k=-[1+\frac{\sigma_x^2}{\sigma_k^2}(D^2S(x))^2]^{-1}(k-\nabla S(x))
\end{equation} therefore:
\begin{equation}
x'-x=\frac{\sigma_x^2}{\sigma_k^2}D^2S(x)[1+ (\frac{\sigma_x}{\sigma_k}D^2S(x) )^2]^{-1}(k-\nabla S(x))
\end{equation}
It is easy to check the non singularity of the Hessian of the phase.

Finally we get the result, for $k\sim \nabla S(x)$, that is, $|\delta|\sim 0$, by expanding the phase around the critical point.
\vskip 0.5cm

Now we are ready to go to our results regarding the smoothed Wigner calculus and equations. First of all,
by obvious adaptation of the respective proof, we readily see that Theorem \ref{thrm4} scales as follows:
\begin{theorem}[Semiclassical smoothed Wigner calculus]
\label{thrm16}
Let $f(x),g(x) \in \mathcal{S}(\mathbb{R}^n)$, $L(x,k) \in \mathcal{S}'(\mathbb{R}^{2n})$ and $L^\varepsilon=L^\varepsilon(x,\varepsilon \partial_x)$ \footnote{See Section \ref{AppA} for the scaled Weyl calculus.}. Then
\begin{equation}
\label{eq127A}
\tilde{W}^\varepsilon[L^\varepsilon f,g](x,k)=\tilde{\mathcal{L}}^\varepsilon \tilde{W}^\varepsilon[f,g](x,k)
\end{equation}
where
\begin{equation}
\label{eq128A}
\begin{array}{l}
\tilde{\mathcal{L}}^\varepsilon \tilde{w}^\varepsilon(x,k)=\\
=\int\limits_{S,T \in \mathbb{R}^n} {
\hat{L}^\varepsilon (S,T) 
e^{2\pi i (Sx+Tk)-\frac{\varepsilon \pi}{2}\left({ \sigma_x^2 S^2 +\sigma^2_k T^2 }\right)} 
\tilde{w}^\varepsilon (x+\varepsilon \frac{T+i\sigma_x^2 S}{2},k-\varepsilon \frac{S-i\sigma_k^2 T}{2})dSdT},
\end{array}
\end{equation}
where $\tilde{w}^\varepsilon (x,k)=\tilde{W}^\varepsilon[f,g](x,k)$ for brevity. In the special case $L^\varepsilon(x,k)=V(x)$, the respective expression is
\begin{equation}
\label{eq131A}
\tilde{W}^\varepsilon[Vf,g](x,k)=
\int\limits_{S \in \mathbb{R}^n} {
\hat{V}(S) 
e^{2\pi i Sx-\frac{\varepsilon \pi}{2} \sigma_x^2 S^2 } 
\tilde{w}^\varepsilon (x+ \frac{i \varepsilon \sigma_x^2 }{2}S,k-\frac{\varepsilon }{2}S)dS}
\end{equation}
\end{theorem}

All the functional analytic framework we constructed for smoothed functions should be scaled correctly with $\varepsilon$; for the most part this is a very predictable exercise. The guideline is naturally the substitution of equation (\ref{eq126A}). 

The following result is an elaboration which will be particularly useful in the sequel:

\begin{theorem}[Semiclassical estimates for smoothed Wigner distributions]
\label{thrm17}
Let $f^\varepsilon(x) \in L^2 (\mathbb{R}^n)$, 
\begin{equation}
\label{eq133}
\tilde{w}^\varepsilon (x,k)=\tilde{W}^\varepsilon [f^\varepsilon](x,k).
\end{equation}
Then the following estimate holds:
\begin{equation}
\label{eq134}
|\tilde{w}^\varepsilon(x+iy,k+iz)| \leqslant \frac{2^n ||f^\varepsilon||^2_{L^2 (\mathbb{R}^n)} }{ \varepsilon^n \sigma_x^n \sigma^n_k} e^{ \frac{2 \pi}{\varepsilon} 
\left({ \frac{|y|^2}{\sigma_x^2} + \frac{|z|^2}{\sigma_k^2} }\right) } \,\,.
\end{equation}

Moreover, let $m_1, m_2, m_3, m_4 \in \left( \mathbb{N} \cup \lbrace 0 \rbrace \right)^n $, and denote $|m_1+m_2+m_3+m_4|=m$.
If $|y|,|z| \leqslant \sqrt{\varepsilon}$, we have
\begin{equation}
\label{eq135}
| \partial_x^{m_1}\partial_y^{m_2}\partial_k^{m_3}\partial_k^{m_4} \tilde{w}^\varepsilon(x+iy,k+iz)| \leqslant F(m_1,m_2,m_3,m_4) \frac{ ||f^\varepsilon||^2_{L^2 (\mathbb{R}^n)} \,\,  e^{ \frac{2 \pi}{\varepsilon} 
\left({ \frac{|y|^2}{\sigma_x^2} + \frac{|z|^2}{\sigma_k^2} }\right) } }{\varepsilon^{\frac{m}{2}+n} \sigma_x^{|m_1+m_2|+n} \sigma^{|m_3+m_4|+n}_k} 
\end{equation}
\end{theorem}

\noindent {\bf Proof: } First we will prove equation (\ref{eq134}). The starting point is the observation that
\begin{equation}
\label{eq136}
\hat{W}^\varepsilon [f^\varepsilon](X,K)=\mathcal{F}_{ x,k \shortrightarrow X,K } \left[ W^\varepsilon [f^\varepsilon](x,k) \right] = \int\limits_{y \in \mathbb{R}^n} { e^{-2\pi i x X} f\left({ x-\frac{\varepsilon K}{2} }\right) \bar{g} \left({ x+\frac{\varepsilon K}{2} }\right) dx}.
\end{equation}
This gives us the uniform in $\varepsilon$ bound
\begin{equation}
\label{eq137}
| \hat{W}^\varepsilon [f^\varepsilon](X,K) | \leqslant ||f^\varepsilon||^2_{ L^2 (\mathbb{R}^n) } .
\end{equation}
Now we have
\begin{equation}
\label{eq138}
\begin{array}{c}
| \tilde{w}^\varepsilon(x+iy,k+iz) |=| \int\limits_{X,K \in \mathbb{R}^n} { e^{2\pi i \left[ (x+iy)X+(k+iz)K \right]} \hat{\tilde{w}}^\varepsilon(X,K) dXdK}|=\\ { } \\
=| \int\limits_{X,K \in \mathbb{R}^n} { e^{2\pi i \left[ (x+iy)X+(k+iz)K \right] - \frac{\pi \varepsilon}{2} \left({ \sigma_x^2 X^2 + \sigma_k^2 K^2 }\right)} \hat{W}^\varepsilon [f^\varepsilon] (X,K) dXdK}|\leqslant \\ { } \\
\leqslant ||f^\varepsilon||^2_{ L^2 (\mathbb{R}^n) } \int\limits_{X,K \in \mathbb{R}^n} { e^{ -\frac{\varepsilon \pi}{2} \left({ \sigma_x^2 X^2 + \sigma_k^2 K^2 }\right)-2\pi \left({ yX+zK }\right) }dXdK }=\frac{2^n ||f^\varepsilon||^2_{ L^2 (\mathbb{R}^n) } }{\varepsilon^n \sigma^n_x \sigma^n_k} e^{ \frac{2 \pi}{\varepsilon} 
\left({ \frac{|y|^2}{\sigma_x^2} + \frac{|z|^2}{\sigma_k^2} }\right) } .
\end{array}
\end{equation}
For equation (\ref{eq135}), observe that
\begin{equation}
\label{eq138a}
\begin{array}{l}
| \partial_x^{m_1}\partial_y^{m_2}\partial_k^{m_3}\partial_k^{m_4} \tilde{w}^\varepsilon(x+iy,k+iz) | \leqslant \\ { } \\
\leqslant (2\pi)^{m} | \int\limits_{X,K \in \mathbb{R}^n } { X^{m_1+m_2} K^{m_3+m_4} e^{2\pi i \left[ (x+iy)X+(k+iz)K \right]} \hat{\tilde{w}}^\varepsilon(X,K) dXdK }| \leqslant \\ { } \\
\leqslant ||f^\varepsilon||^2_{ L^2 (\mathbb{R}^n) } (2\pi)^{m} \int\limits_{X \in \mathbb{R}^n } { \left| X^{m_1+m_2} \right| e^{ -\frac{\varepsilon \pi}{2}  \sigma_x^2 X^2 -2\pi yX } }
\int\limits_{K \in \mathbb{R}^n } { \left| K^{m_3+m_4} \right| e^{ -\frac{\varepsilon \pi}{2}  \sigma_k^2 K^2 -2\pi zK } } \leqslant \\ { } \\  

\leqslant
||f^\varepsilon||^2_{ L^2 (\mathbb{R}^n) } (2\pi)^{m} 
\frac{e^{ \frac{2 \pi}{\varepsilon} \left({ \frac{|y|^2}{\sigma_x^2} + \frac{|z|^2}{\sigma_k^2} }\right) } }
{\left({ \frac{\varepsilon \pi}{2} }\right)^{\frac{m}{2}+n} \sigma_x^{|m_1+m_2|+n} \sigma^{|m_3+m_4|+n}_k} \, \cdot \\

\,\,\,\,\,\,\,\,\,\,\,\,\,\,\,\,\,\,\,\,\,\,\,\,\,\,\,\,\,\,\,\,\,\,\,\,\,\,\,\,\,\,\,\,\,\,\,\,\,\,\,\,\,\,\,\,
\cdot \, \prod\limits_{d=1}^n { \left( 2 \sum\limits_{l=0}^{|m_1+m_2|_d} { \binom{|m_1+m_2|_d}{l} \left({ \frac{\sqrt{2 \pi} |y_d|}{\sqrt{\varepsilon} \sigma_x} }\right)^{|m_1+m_2|_d-l} \Gamma \left({ \frac{l+1}{2} }\right) } \right)} \, \cdot \\

\,\,\,\,\,\,\,\,\,\,\,\,\,\,\,\,\,\,\,\,\,\,\,\,\,\,\,\,\,\,\,\,\,\,\,\,\,\,\,\,\,\,\,\,\,\,\,\,\,\,\,\,\,\,\,\,
\cdot \, \prod\limits_{d=1}^n {\left( 2 \sum\limits_{l=0}^{|m_3+m_4|_d} { \binom{|m_3+m_4|_d}{l} \left({ \frac{\sqrt{2 \pi} |z_d|}{\sqrt{\varepsilon} \sigma_k} }\right)^{|m_3+m_4|_d-l} \Gamma \left({ \frac{l+1}{2} }\right) } \leqslant \right) } \\ { } \\
\leqslant
F(m_1,m_2,m_3,m_4) \frac{ ||f^\varepsilon||^2_{ L^2 (\mathbb{R}^n) } \,  e^{ \frac{2 \pi}{\varepsilon} 
\left({ \frac{|y|^2}{\sigma_x^2} + \frac{|z|^2}{\sigma_k^2} }\right) } }{\varepsilon^{\frac{m}{2}+n} \sigma_x^{|m_1+m_2|+n} \sigma^{|m_3+m_4|+n}_k}.
\end{array}
\end{equation}

The elementary computation which allows us to pass from the third line to the fourth, is
\begin{equation}
\label{eq139}
\begin{array}{c}
\int\limits_{r=0}^{+\infty} { x^m e^{-ax^2-bx} dx } =
\frac{ e^{\frac{b^2}{4a}} }{a^{\frac{m+1}{2}}} \sum\limits_{l=0}^m {
\binom{m}{l} \left( \frac{b}{2\sqrt{a}} \right)^{m-l} \int\limits_{u=\frac{b}{2\sqrt{a}}}^{+\infty} { u^l e^{-u^2} du} } \leqslant \\ { } \\

\leqslant 2 \frac{ e^{\frac{b^2}{4a}} }{a^{\frac{m+1}{2}}} \sum\limits_{l=0}^m {
\binom{m}{l} \left( \frac{b}{2\sqrt{a}} \right)^{m-l} \Gamma \left( \frac{l+1}{2} \right) } 
\end{array}
\end{equation}

The assumption $|y|,|z| \leqslant \sqrt{\varepsilon}$ implies that $\frac{|y_d|}{\sqrt{\varepsilon}}, \frac{|z_d|}{\sqrt{\varepsilon}} \, \leqslant \, 1$, and therefore we can pass to the last line (equivalently, $F(m_1,m_2,m_3,m_4)$ is independent of $\varepsilon$).

The proof is complete. \\ \vspace{2mm}

Now we are ready to see how the smoothed Wigner calculus can be approximated by differential operators in the semiclassical regime. This is the kind of computation necessary for the formulation of asymptotic SWT-based models: \\ \vspace{1mm}

\begin{theorem}[Semiclassical finite-order approximations to the smoothed Wigner calculus]
\label{thrm15}
Let $N \in \mathbb{N}$, $V(x):\mathbb{R}^n \shortrightarrow \mathbb{R}$, $\varepsilon \in (0,1)$. Assume that
\begin{itemize}
\item[\textbf{(A1)}]
$\hat{V}(k)$ has no singular support outside $\lbrace 0 \rbrace$ \footnote{In fact we could handle singular support away from $0$ with no big problems; we exclude it here for simplicity.}.
\item[\textbf{(A2)}]
$\exists C_1>0, 0 \leqslant M_1\leqslant N+1$ such that
\begin{equation}
\label{eq140}
|\hat{V}(k)| \leqslant C_1 |k|^{-M_1} \,\,\,\, \forall |k| \leqslant 1, k \neq 0
\end{equation}
\item[\textbf{(A3)}]
For an appropriate\footnote{This is not the full assumption for $M_2$. See remark 3 below, and the remarks at the end of the proof.} (finite) constant $ M_2=M_2 (n,N) \leqslant min \lbrace -n-1, -3 \rbrace$, $\exists  C_2 >0$  such that
\begin{equation}
\label{eq143}
|\hat{V}(k)| \leqslant C_2 |k|^{M_2} \,\,\,\, \forall |k| > 1 .
\end{equation}
\end{itemize}
Moreover, consider a ``semiclassical family of wavefunctions'' $\lbrace f^\varepsilon \rbrace \subset \mathcal{S}(\mathbb{R}^n)$  generating the SWTs
\begin{equation}
\label{eq150}
\tilde{w}^\varepsilon(x,k)=\tilde{W}^\varepsilon[f^\varepsilon](x,k),
\end{equation}
for which we assume that $\exists M_0 >0$ such that
\begin{equation}
\label{eq149}
||f^\varepsilon ||_{L^2 (\mathbb{R}^n) } \leqslant M_0 \,\,\,\, \forall \varepsilon >0 .
\end{equation}

According to Theorem \ref{thrm16},
\begin{equation}
\label{eq141}
\tilde{W}^\varepsilon [V f^\varepsilon, f^\varepsilon](x,k)=
\int\limits_{S \in \mathbb{R}^n} 
{ e^{ 2\pi i S x - \frac{\varepsilon \pi}{2} \sigma_x^2 S^2 } \hat{V}(S) \tilde{w}^\varepsilon(x+\frac{i\varepsilon \sigma_x^2}{2}S,k-\frac{\varepsilon}{2}S)dS } .
\end{equation}
Assuming in addition that $\sigma_x^2 \leqslant 2$, this expression can be approximated by differential operators,
\begin{equation}
\label{eq1510}
\begin{array}{l}
\tilde{W}^\varepsilon [V f^\varepsilon, f^\varepsilon](x,k)= \\ { } \\
%
=\sum\limits_{m=0}^N \frac{\varepsilon^m}{(4\pi i)^m} \sum\limits_{l=0}^m \left( i\sigma_x^2 \right)^l (-1)^{m-l}
\sum\limits_{
\begin{scriptsize}
\begin{array}{c}
A \in (\mathbb{N} \cup \lbrace 0 \rbrace )^n \\
|A|=l
\end{array}
\end{scriptsize} } 
\sum\limits_{
\begin{scriptsize}
\begin{array}{c}
B \in (\mathbb{N} \cup \lbrace 0 \rbrace )^n \\
|B|=m-l
\end{array}
\end{scriptsize} }  
\frac{\partial_x^{A+B}  \tilde{V} (x)}{A!B!} 
\prod\limits_{d=1}^n { \prod\limits_{d'=1}^n { \partial^{A_d}_{x_d} \partial^{B_{d'} }_{k_d'}  } } \tilde{w}^\varepsilon(x,k)+\\ { } \\
\,\,\,\,\,\,\,\,\,\,\,\,\,\,\,\,\,\,\,\,\,\,\,\,\,\,\,\,\,\,\,\,\,\,\,\,\,\,\,\,\,\,\,\,\,\,\,\,\,\,\,\,\,\,\,\,\,\,\,\,\,\,\,\,\,\,\,\,\,\,\,\,\,\,\,\,\,\,\,\,\,\,\,\,\,\,\,\,\,\,\,\,\,\,\,\,\,\,\,\,\,\,\,\,\,\,\,\,\,\,\,\,\,\,\,\,\,\,\,\,\,\,\,\,\,\,\,\,\,\,\,\,\,\,\,\,\,\,\,\,\,\,\,\,\,\,\,\,\,\,\,\,\,\,\,
+r_\varepsilon(x,k),
\end{array}
\end{equation}
where
\begin{equation}
\label{err}
|| r_\varepsilon ||_{L^\infty (\mathbb{R}^{2n})} \, = \, O\left( \varepsilon^{\frac{N+1}{2}-n} \right),
\end{equation}
and $\tilde{V}(x)$ is the potential $V(x)$ smoothed at scale $\sigma_x^2$,
\begin{equation}
\label{eq144m}
\tilde{V}(x)=
\left( \frac{\sqrt{2}}{\sqrt{\varepsilon} \sigma_x} \right)^n \int\limits_{x' \in \mathbb{R}^n} {e^{-\frac{2\pi |x-x'|^2}{\varepsilon \sigma_x^2}} V(x')dx'}
=\Phi_{\sqrt{\varepsilon} \sigma_x} V (x).
\end{equation}

%
\end{theorem}

\noindent {\bf Remarks:} 
\begin{enumerate}
\item
Equation (\ref{eqTh2}) shows clearly that the $L^\infty$ norm of the SWT of a WKB function is of the order $\varepsilon^{-n}$. This makes the estimate of equation (\ref{err})  significant, since it can be written as
\begin{equation}
|| r_\varepsilon ||_{L^\infty (\mathbb{R}^{2n})} \, = \, O\left( \varepsilon^{\frac{N+1}{2}}  \,\, || \tilde{W}^\varepsilon[f^\varepsilon] ||_{L^\infty (\mathbb{R}^{2n})} \,\,\, \right).
\end{equation}
Moreover the same estimates as (\ref{eqTh2}) can be proved to be valid for the suitably scaled derivatives of $\tilde{W}^\varepsilon[f^\varepsilon]$, making the estimate (\ref{err}) sharp.

\item
A qualitative description of the result: for appropriately (but finitely in any case) smooth potentials, the smoothed Wigner calculus (and therefore the ``scattering term'' in the smoothed Wigner equation (\ref{eq116a}) ) can be approximated uniformly by a differential operator in the semiclassical regime. Of course $\tilde{w}^\varepsilon (x,k)$ and $\tilde{W}^\varepsilon [V f^\varepsilon, f^\varepsilon](x,k)$ themselves become unbounded pointwise as $\varepsilon \shortrightarrow 0$, but still we can approximate them strongly.

This should be compared of course to the weak approximation of the Wigner calculus that is the standard device for constructing asymptotic equations for Wigner measures. Indeed, this result is a precise quantification of the argument that ``the SWT is better suited to keep track of the wavefield in the semiclassical regime than the WT''.

\item
A note must be made on the selection of $M_2$: at several instances along the proof, a condition of the type $M_2 \leqslant s_0$ will appear. Some of the conditions originally appear not in that form, but in all cases they can be satisfied by choosing $M_2$ small enough.  The collection of these conditions (which depend on $n,N$ as well) is the actual assumption which has to be satisfied by $M_2$. References to all the conditions are gathered in a remark in the end of the proof.

\end{enumerate}

\noindent {\bf Proof: } The central idea of the proof is actually very simple: we break the $dS$ integral over a neighbourhood of zero and its complement, we  Taylor-expand $\tilde{w}^\varepsilon(x+\frac{i\varepsilon \sigma_x^2}{2}S,k-\frac{\varepsilon}{2}S)$ in equation (\ref{eq141}) around $(x,k)$ up to order $N$ and keep the remainder. Then we compute the approximation error of the Taylor expansion, and the contribution of the $dS$ integral away from zero.

First of all let us fix notations on the Taylor expansion: if $g:\mathbb{R}^n \shortrightarrow \mathbb{R}$ is a sufficiently smooth function, then
\begin{equation}
\label{eq148}
g(x)=\sum\limits_{m=0}^N \frac{ \left[ x_1 \partial_{x_1}+...+x_n \partial_{x_n} \right]^m }{m!} g(0)+ R_N(x),
\end{equation}
where the remainder can be described as follows: $\exists \theta=\theta(x) \in (0,1)$ such that
\begin{equation}
\label{eq153}
R_N(x)=
\frac{ 
\left( \sum\limits_{d=1}^n { x_d \partial_{x_d} }\right)^{N+1} 
}{
(N+1)!
} g(\theta x).
\end{equation}
Now define
\begin{equation}
\label{eq154}
g(S)=\tilde{w}^\varepsilon(x+\frac{i\varepsilon \sigma_x^2}{2}S,k-\frac{\varepsilon}{2}S) .
\end{equation}
It is clear that the Taylor Theorem is applicable \footnote{Indeed, $g(S)$ can be seen is actually an entire function; observe however that we are only interested in $S \in \mathbb{R}^n$, so we can use the Taylor expansion for real functions.}. Observe moreover that
\begin{equation}
\label{eq158}
\partial_{S_d} g(S)=\left[ \frac{i\varepsilon \sigma_x^2}{2} \partial_{x_d} - \frac{\varepsilon}{2} \partial_{k_d} \right]  \tilde{w}^\varepsilon(x+\frac{i\varepsilon \sigma_x^2}{2}S,k-\frac{\varepsilon}{2}S) .
\end{equation}
Now we have the Taylor expansion
\begin{equation}
\label{eq155}
g(S)=\sum\limits_{m=0}^N \frac{ \left[ S_1 \partial_{S_1}+...+S_n \partial_{S_n} \right]^m }{m!} g(0)+ R_N(S),
\end{equation}
where, for each $S$, the remainder is given by
\begin{equation}
\label{eq156}
\begin{array}{l}
R_N(S)=
\frac{\varepsilon^{N+1}}{ (N+1)! 2^{N+1}} 
\left[ \sum\limits_{d=1}^n  i\sigma_x^2 S_d \partial_{x_d} - S_d \partial_{k_d}  \right]^{N+1}
\tilde{w}^{\varepsilon} \left( x+\theta \frac{i \varepsilon\sigma_x^2 S}{2},k-\theta\frac{\varepsilon S}{2}  \right) = \\ { } \\

=\frac{\varepsilon^{N+1}}{ (N+1)! 2^{N+1}} 
\sum\limits_{l=0}^{N+1} \binom{N+1}{l}
\left( \sum\limits_{d=1}^n  i\sigma_x^2 S_d \partial_{x_d} \right)^l
\left( \sum\limits_{d=1}^n  - S_d \partial_{k_d} \right)^{N+1-l}
\tilde{w}^{\varepsilon} \left( x+\theta \frac{i \varepsilon\sigma_x^2 S}{2},k-\theta\frac{\varepsilon S}{2}  \right) = \\ { } \\

=\frac{\varepsilon^{N+1}}{ 2^{N+1}} 
\sum\limits_{l=0}^{N+1} \left( i\sigma_x^2 \right)^l (-1)^{N+1-l} \\

\left( 
\sum\limits_{
\begin{scriptsize}
\begin{array}{c}
A \in (\mathbb{N} \cup \lbrace 0 \rbrace )^n \\
|A|=l
\end{array}
\end{scriptsize} } 
{ \frac{\binom{l}{A}}{l!} \prod\limits_{d=1}^n S^{A_d}_{d} \partial^{A_d}_{x_d} }
\right)

\left( 
\sum\limits_{
\begin{scriptsize}
\begin{array}{c}
B \in (\mathbb{N} \cup \lbrace 0 \rbrace )^n \\
|B|=N+1-l
\end{array}
\end{scriptsize} } 
{ \frac{\binom{N+1-l}{B}}{N+1-l!} \prod\limits_{d=1}^n S^{B_d}_{d} \partial^{B_d}_{k_d} }
\right)
\tilde{w}^{\varepsilon} \left( x+\theta \frac{i \varepsilon\sigma_x^2 S}{2},k-\theta\frac{\varepsilon S}{2}  \right) = \\ { } \\

=\frac{\varepsilon^{N+1}}{ 2^{N+1}} 
\sum\limits_{l=0}^{N+1} \left( i\sigma_x^2 \right)^l (-1)^{N+1-l} \\

\left( 
\sum\limits_{
\begin{scriptsize}
\begin{array}{c}
A \in (\mathbb{N} \cup \lbrace 0 \rbrace )^n \\
|A|=l
\end{array}
\end{scriptsize} } 
{ \frac{1}{A!} \prod\limits_{d=1}^n S^{A_d}_{d} \partial^{A_d}_{x_d} }
\right)

\left( 
\sum\limits_{
\begin{scriptsize}
\begin{array}{c}
B \in (\mathbb{N} \cup \lbrace 0 \rbrace )^n \\
|B|=N+1-l
\end{array}
\end{scriptsize} } 
{ \frac{1}{B!} \prod\limits_{d=1}^n S^{B_d}_{d} \partial^{B_d}_{k_d} }
\right)
\tilde{w}^{\varepsilon} \left( x+\theta \frac{i \varepsilon\sigma_x^2 S}{2},k-\theta\frac{\varepsilon S}{2}  \right) .
\end{array}
\end{equation}
It follows therefore that, if $|S| \leqslant r$,
\begin{equation}
\label{eq157}
\begin{array}{l}
|R_N(S)|\leqslant \frac{(\varepsilon |S|)^{N+1} }{2^{N+1}}
\sum\limits_{l=0}^{N+1} \sigma_x^{2l} \,\,
\left( 
\sum\limits_{
\begin{scriptsize}
\begin{array}{c}
A \in (\mathbb{N} \cup \lbrace 0 \rbrace )^n \\
|A|=l
\end{array}
\end{scriptsize} } 
{ \frac{1}{A!}  }
\right)

\left( 
\sum\limits_{
\begin{scriptsize}
\begin{array}{c}
B \in (\mathbb{N} \cup \lbrace 0 \rbrace )^n \\
|B|=N+1-l
\end{array}
\end{scriptsize} } 
{ \frac{1}{B!}  }
\right) \\
\,\,\,\,\,\,\,\,\,\,\,\,\,\,\,\,\,
\,\,\,\,\,\,\,\,\,\,\,\,\,\,\,\,\,\,\,\,\,\,\,\,\,\,\,\,\,\,\,\,\,\,\,\,\,\,\,\,\,\, \,\, \sup\limits_{
\begin{scriptsize}\begin{array}{c}
|S|\leqslant r \\
A' \in (\mathbb{N} \cup \lbrace 0 \rbrace )^{2n}\\
|A'|=N+1
\end{array}
\end{scriptsize}} \, | \partial^{A'}_{x_1...x_nk_1...k_n} \tilde{w}^{\varepsilon} \left( x+ \frac{i \varepsilon\sigma_x^2 S}{2},k-\frac{\varepsilon S}{2}  \right) | \leqslant \\ { } \\

\leqslant \frac{(\varepsilon |S|)^{N+1} }{2^{N+1}} \, max \lbrace 1, \sigma_x^{2(N+1)} \rbrace \,
\sum\limits_{l=0}^{N+1} \frac{n^{N+1}}{l!(N+1-l)!} \,\,
\sup\limits_{
\begin{scriptsize}\begin{array}{c}
|S|\leqslant r \\
A' \in (\mathbb{N} \cup \lbrace 0 \rbrace )^{2n}\\
|A'|=N+1
\end{array}
\end{scriptsize}} \, | \partial^{A'}_{x_1...x_nk_1...k_n} \tilde{w}^{\varepsilon} \left( x+ \frac{i \varepsilon\sigma_x^2 S}{2},k-\frac{\varepsilon S}{2}  \right) |=\\ { } \\

=\frac{(n \varepsilon |S|)^{N+1} }{(N+1)!} \, max \lbrace 1, \sigma_x^{2(N+1)} \rbrace \,
\,
\sup\limits_{
\begin{scriptsize}\begin{array}{c}
|S|\leqslant r \\
A' \in (\mathbb{N} \cup \lbrace 0 \rbrace )^{2n}\\
|A'|=N+1
\end{array}
\end{scriptsize}} \, | \partial^{A'}_{x_1...x_nk_1...k_n} \tilde{w}^{\varepsilon} \left( x+ \frac{i \varepsilon\sigma_x^2 S}{2},k-\frac{\varepsilon S}{2}  \right) |.
\end{array}
\end{equation}
At this point we need to use Lemma \ref{thrm17}. To do that, we have to check that
\begin{equation}
\label{eq160}
\frac{\varepsilon \sigma_x^2}{2}|S| \leqslant \sqrt{\varepsilon}.
\end{equation}
One thing we will do is use the assumption $\sigma_x^2 \leqslant 2$; moreover, (for reasons that will become more clear below), we will set $r=\varepsilon^{ -\frac{1}{4}+\frac{1}{2(M_2+1)} }$, and therefore $\varepsilon |S| \leqslant \varepsilon^{\frac{3}{4} +\frac{1}{2(M_2+1)}}$. So using finally the constraint 
\begin{equation}
\label{eq172}
M_2 \leqslant -3,
\end{equation}
it follows that $\varepsilon r \leqslant \varepsilon^{\frac{1}{2}}$ and
equation (\ref{eq160}) holds. In particular, the assumptions of Theorem \ref{thrm17} are satisfied. 

Now, using Theorem \ref{thrm17}, and more precisely equation (\ref{eq135}), as well as equation (\ref{eq149}), it follows that there is a constant $\tilde{C}=\tilde{C}(N,n,\sigma_x^2,\sigma_k^2)$ such that
\begin{equation}
\label{eq159}
\sup\limits_{
\begin{scriptsize}\begin{array}{c}
|S|\leqslant r \\
A' \in (\mathbb{N} \cup \lbrace 0 \rbrace )^{2n}\\
|A'|=N+1
\end{array}
\end{scriptsize}} \, | \partial^{A'}_{x_1...x_nk_1...k_n} \tilde{w}^{\varepsilon} \left( x+ \frac{i \varepsilon\sigma_x^2 S}{2},k-\frac{\varepsilon S}{2}  \right) | \leqslant \tilde{C} \frac{ M^2_0 \,\,  e^{ \frac{\varepsilon \pi}{2} 
\sigma_x^2 |S|^2  } }{\varepsilon^{\frac{N+1}{2}+n} }.
\end{equation}

So now it follows that the remainder of the Taylor expansion of
\begin{equation}
g(S)=\tilde{w}^\varepsilon(x+\frac{i\varepsilon \sigma_x^2}{2}S,k-\frac{\varepsilon}{2}S)
\end{equation}
around $S=0$ (and for $S$ in any case not larger than $|S| \leqslant \varepsilon^{ -\frac{1}{4}+\frac{1}{2(M_2+1)} }$), is dominated by
\begin{equation}
\label{eq161}
|R_N(S)| \leqslant H \,\, \varepsilon^{\frac{N+1}{2}-n} |S|^{N+1} e^{\frac{\pi \varepsilon \sigma_x^2}{2}|S|^2}
\end{equation}
for some constant $H=H(N,n,\sigma_x^2,\sigma_k^2)$.

The next part of the proof is simple (if a little tedious): we break the integral of equation (\ref{eq141}) to $I_1 = \int\limits_{|S| \leqslant \varepsilon^{ -\frac{1}{4}+\frac{1}{2(M_2+1)} }}$ and $I_2 = \int\limits_{|S| > \varepsilon^{ -\frac{1}{4}+\frac{1}{2(M_2+1)} }}$. We use the $N$-order Taylor expansion of $g(S)$ in $I_1$, and bound the error using equation (\ref{eq161}).  For the contribution of $I_2$, we will use the estimate of equation (\ref{eq134}) . Of course some more auxiliary assumptions (described in the statement) will come up along the way. (See also the remarks at the end of the proof).

The first contribution to the error comes from
\begin{equation}
\label{eq162}
E_1=\int\limits_{|S| \leqslant \varepsilon^{ -\frac{1}{4}+\frac{1}{2(M_2+1)} }}
{ e^{ 2\pi i S x - \frac{\varepsilon \pi}{2} \sigma_x^2 S^2 } \hat{V}(S) R_N(S)dS }.
\end{equation}
Using the previous results and the assumptions of the Theorem we have
\begin{equation}
\label{eq163}
\begin{array}{l}
|E_1| \leqslant  H \,  \varepsilon^{\frac{N+1}{2}-n} \,
\int\limits_{|S| \leqslant  \varepsilon^{ -\frac{1}{4}+\frac{1}{2(M_2+1)} }}
{  |\hat{V}(S)| |S|^{N+1} dS } \leqslant \\ { } \\
\leqslant H \,  \varepsilon^{\frac{N+1}{2}-n} \,
\left[
C_1 \int\limits_{|S| \leqslant  1}
{ |S|^{-M_1+N+1}dS}
+C_2 \int\limits_{1 < |S| \leqslant  \varepsilon^{ -\frac{1}{4}+\frac{1}{2(M_2+1)} }}
{|S|^{M_2+N+1}  dS}
\right]= \\ { } \\
=H \,  \varepsilon^{\frac{N+1}{2}-n} \, \frac{2 \pi ^{\frac{n}{2}} }{\Gamma \left({\frac{n}{2}}\right) } \,
\left[
C_1 \int\limits_{r=0}^1
{ r^{-M_1+N+n}dS}
+C_2 \int\limits_{r=1}^{\varepsilon^{ -\frac{1}{4}+\frac{1}{2(M_2+1)} }}
{r^{M_2+N+n}  dS}
\right].
\end{array}
\end{equation}
Here we make use of the assumption $M_1 \leqslant N+n$, and moreover we assume
\begin{equation}
\label{eq164}
M_2+N+n+1 \neq 0.
\end{equation}
Now we are able to proceed to
\begin{equation}
\label{eq165}
|E_1| \leqslant H \,  \varepsilon^{\frac{N+1}{2}-n} \, \frac{2 \pi ^{\frac{n}{2}} }{\Gamma \left({\frac{n}{2}}\right) } \,
\left[ 
\frac{C_1}{N+n+1-M_1} + C_2 \left| \frac{\varepsilon^{ 
(M_2+N+n+1)\left( -\frac{1}{4}+\frac{1}{2(M_2+1)} \right)
 }-1 }{M_2+N+n+1} \right|
\right].
\end{equation}
Clearly the best we can ask for here is
\begin{equation}
\label{eq166}
(M_2+N+n+1)\left( -\frac{1}{4}+\frac{1}{2(M_2+1)} \right) \geqslant 0,
\end{equation}
so that the total bound for $E_1$ is controlled by $\varepsilon^{\frac{N+1}{2}-n}$ (it is obvious that for $M_2$ small enough the inequality (\ref{eq166}) holds). 

On the contribution of $I_2$, making use of Theorem \ref{thrm17}, we observe that
\begin{equation}
\label{eq167}
\begin{array}{l}
|I_2| \leqslant \int\limits_{|S| > \varepsilon^{ -\frac{1}{4}+\frac{1}{2(M_2+1)} }}
{ |\hat{V}(S)| e^{-\frac{\pi \varepsilon \sigma_x^2 |S|^2}{2}} |\tilde{w}(x+\frac{i\varepsilon \sigma_x^2}{2}S,k-\frac{\varepsilon}{2}S)| dS} \leqslant \\ { } \\
\leqslant \frac{C_2 M_0^2}{\varepsilon} \frac{2 \pi ^{\frac{n}{2}} }{\Gamma \left({\frac{n}{2}}\right) }
\int\limits_{r=\varepsilon^{ -\frac{1}{4}+\frac{1}{2(M_2+1)} }}^{+\infty}
{r^{M_2+n-1}dr}
\end{array}
\end{equation}
We will have to assume
\begin{equation}
\label{eq168}
M_2 \leqslant -1-n
\end{equation}
for the integral to exist (observe however that this is not automatically enough for $\lim\limits_{\varepsilon \shortrightarrow 0} I_2 =0$); with that we get
\begin{equation}
\label{eq169}
\begin{array}{l}
|I_2| \leqslant 
\frac{C_2 M_0^2}{ (-M_2-n)} \frac{2 \pi ^{\frac{n}{2}} }{\Gamma \left({\frac{n}{2}}\right) }
\varepsilon^{ -\frac{1}{2}-\frac{M_2+n}{4}+\frac{n-1}{2(M_2+1)} }.
\end{array}
\end{equation}
Like earlier, we ask not only that $\lim\limits_{\varepsilon \shortrightarrow 0} I_2 =0$, but that $|I_2|$ is controlled by $\varepsilon^{\frac{N+1}{2}-n}$. This amounts to
\begin{equation}
\label{eq170}
-\frac{1}{2}-\frac{M_2+n}{4}+\frac{n-1}{2(M_2+1)} \geqslant \frac{N+1}{2}-n.
\end{equation}

Now observe that equation (\ref{eq156}) gives the structure of all the terms in the Taylor expansion, not only the remainder. That is, the order-$m$ term of the Taylor expansion (\ref{eq155}) is given by
\begin{equation}
\label{eq1562}
\begin{array}{l}
T_m(S)=
\frac{\varepsilon^{m}}{ 2^{m}} 
\sum\limits_{l=0}^{m} \left( i\sigma_x^2 \right)^l (-1)^{m-l} \\

\left( 
\sum\limits_{
\begin{scriptsize}
\begin{array}{c}
A \in (\mathbb{N} \cup \lbrace 0 \rbrace )^n \\
|A|=l
\end{array}
\end{scriptsize} } 
{ \frac{1}{A!} \prod\limits_{d=1}^n S^{A_d}_{d} \partial^{A_d}_{x_d} }
\right)

\left( 
\sum\limits_{
\begin{scriptsize}
\begin{array}{c}
B \in (\mathbb{N} \cup \lbrace 0 \rbrace )^n \\
|B|=m-l
\end{array}
\end{scriptsize} } 
{ \frac{1}{B!} \prod\limits_{d=1}^n S^{B_d}_{d} \partial^{B_d}_{k_d} }
\right)
\tilde{w}^{\varepsilon} \left( x,k  \right).
\end{array}
\end{equation}

So what we have shown so far is that
\begin{equation}
\label{eq151o}
\begin{array}{l}
\,\,\,\,\,\,\,\,\,
\tilde{W}^\varepsilon [V f^\varepsilon, f^\varepsilon](x,k)=\\ { } \\
=\int\limits_{|S| \leqslant \, \varepsilon^{-\frac{1}{4}+\frac{1}{2(M_2+1)}  } }
{
e^{ 2\pi i Sx - \frac{\varepsilon \pi}{2} \sigma_x^2S^2 } \hat{V}(S)
\sum\limits_{m=0}^N {
\frac{\varepsilon ^m}{m!  2^m} \left[ \sum\limits_{d=1}^n { i \sigma_x^2 S_d \partial_{x_d}-S_d \partial_{k_d} } \right]^m
dS \tilde{w}^\varepsilon (x,k)} } + \tilde{r}_\varepsilon(x,k) = \\ { } \\
=\sum\limits_{m=0}^N \frac{\varepsilon^m}{2^m} \sum\limits_{l=0}^m \left( i\sigma_x^2 \right)^l (-1)^{m-l}
\sum\limits_{
\begin{scriptsize}
\begin{array}{c}
A \in (\mathbb{N} \cup \lbrace 0 \rbrace )^n \\
|A|=l
\end{array}
\end{scriptsize} } 
\sum\limits_{
\begin{scriptsize}
\begin{array}{c}
B \in (\mathbb{N} \cup \lbrace 0 \rbrace )^n \\
|B|=m-l
\end{array}
\end{scriptsize} }  
\frac{D(A+B)}{A!B!} 
\prod\limits_{d=1}^n { \prod\limits_{d'=1}^n { \partial^{A_d}_{x_d} \partial^{B_{d'} }_{k_d'}  } } \tilde{w}^\varepsilon(x,k)
+\tilde{r}_\varepsilon(x,k),
\end{array}
\end{equation}
where the coefficients $D(A)$ are given by
\begin{equation}
\label{eq144o}
D(A)(x)=\int\limits_{  |S| \leqslant \varepsilon ^{-\frac{1}{4}+\frac{1}{2(M_2+1)}} }{
e^{-2\pi i Sx-\frac{\varepsilon \pi}{2} \sigma_x^2S^2} \hat{V}(S) \prod\limits_{d=1}^n S_d^{A} dS},
\end{equation}
and 
\begin{equation}
\label{erro}
|| \tilde{r}_\varepsilon ||_{L^\infty (\mathbb{R}^{2n})} \, = \, O\left( \varepsilon^{\frac{N+1}{2}-n} \right).
\end{equation}

Observe that the coefficients $D(A+B)$ are truncated versions of the derivatives of the smoothed potential. Indeed, denote
\begin{equation}
\label{wec1}
\begin{array}{c}
\tilde{D}(A)(x)=\int\limits_{  S \in \mathbb{R}^n }{
e^{-2\pi i Sx-\frac{\varepsilon \pi}{2} \sigma_x^2S^2} \hat{V}(S) \prod\limits_{d=1}^n S_d^{A_d} dS}= \\ { } \\

=\frac{\partial_x^{A}}{\left( 2\pi i \right)^{|A|} } \left( \frac{\sqrt{2}}{\sqrt{\varepsilon} \sigma_x} \right)^n \int\limits_{x' \in \mathbb{R}^n} {e^{-\frac{2\pi |x-x'|^2}{\varepsilon \sigma_x^2}} V(x')dx'}
=\frac{\partial_x^{A}}{\left( 2\pi i \right)^{|A|} } \Phi_{\sqrt{\varepsilon} \sigma_x} V (x) ;
\end{array}
\end{equation}
then
\begin{equation}
\label{wec2}
\begin{array}{c}
I_3(A)(x)=\tilde{D}(A)(x)-D(A)(x)=\int\limits_{ |S| > \varepsilon ^{-\frac{1}{4}+\frac{1}{2(M_2+1)}} }{
e^{-2\pi i Sx-\frac{\varepsilon \pi}{2} \sigma_x^2S^2} \hat{V}(S) \prod\limits_{d=1}^n S_d^{A_d} dS}.
\end{array}
\end{equation}
This difference is small, and the final part of the proof consists in showing that we can substitute $\tilde{D}(A)$ for $D(A)$ and have a similar error estimate as in equations (\ref{eq151o}), (\ref{erro}). To that end, it suffices to show that $\exists C>0$ such that $\forall x \in \mathbb{R}^n$ and $\forall \, A \in (\mathbb{N} \cup \lbrace 0 \rbrace )^n$, $|A| \leqslant{N}$, 
\begin{equation}
\label{est1}
|I_3(A)(x)| \leqslant C \varepsilon^{\frac{N+1}{2}}.
\end{equation}
Then (using Theorem \ref{thrm17} once more) the error introduced by substituting $\tilde{D}(A)$ for $D(A)$ in equation (\ref{eq151o}) will be dominated by
\begin{equation}
\label{eq151o1}
\begin{array}{c}
\sum\limits_{m=0}^N \frac{\varepsilon^m}{2^m} \sum\limits_{l=0}^m \left( \sigma_x^2 \right)^l 
\sum\limits_{
\begin{scriptsize}
\begin{array}{c}
A \in (\mathbb{N} \cup \lbrace 0 \rbrace )^n \\
|A|=l
\end{array}
\end{scriptsize} } 
\sum\limits_{
\begin{scriptsize}
\begin{array}{c}
B \in (\mathbb{N} \cup \lbrace 0 \rbrace )^n \\
|B|=m-l
\end{array}
\end{scriptsize} }  
\frac{|I_3(A+B)|}{A!B!} 
\prod\limits_{d=1}^n { \prod\limits_{d'=1}^n { |\partial^{A_d}_{x_d} \partial^{B_{d'} }_{k_d'}  } } \tilde{w}^\varepsilon(x,k)|
+|\tilde{r}_\varepsilon(x,k)| \leqslant \\ { } \\
\leqslant C' \sum\limits_{m=0}^N \varepsilon^m \varepsilon^{\frac{N+1}{2}} \varepsilon^{-\frac{m}{2}-n} + 
|\tilde{r}_\varepsilon(x,k)|= O\left( \varepsilon^{\frac{N+1}{2}-n} \right)
\end{array}
\end{equation}

So let us prove equation (\ref{est1}):
\begin{equation}
\begin{array}{c}
\label{est2}
|I_3(A)(x)| \leqslant \int\limits_{ |S| > \varepsilon ^{-\frac{1}{4}+\frac{1}{2(M_2+1)}} }{
e^{-\frac{\varepsilon \pi}{2} \sigma_x^2S^2} |\hat{V}(S)| \prod\limits_{d=1}^n |S_d|^{A_d} dS} \leqslant \\ { } \\

\leqslant C_2 \int\limits_{ |S| > \varepsilon ^{-\frac{1}{4}+\frac{1}{2(M_2+1)}} }{
e^{-\frac{\varepsilon \pi}{2} \sigma_x^2S^2} |S|^{M_2+n|A|}  dS} = \\ { } \\

= C_2 \frac{2 \pi^{\frac{n}2}}{\Gamma\left( \frac{n}2 \right)} 
\int\limits_{ r = \varepsilon ^{-\frac{1}{4}+\frac{1}{2(M_2+1)}} }^{+\infty}
{e^{-\frac{\varepsilon \pi}{2} \sigma_x^2 r^2} r^{M_2+n|A|+n-1}  dr} = \\ { } \\

= \frac{2 \pi^{\frac{n}2}}{\Gamma\left( \frac{n}2 \right)} \frac{C_2}{\left( -\varepsilon \pi \sigma_x^2 \right)^{M_2+n|A|+n-1}}
\int\limits_{ r = \varepsilon ^{-\frac{1}{4}+\frac{1}{2(M_2+1)}} }^{+\infty}
{ \partial_r^{M_2+n|A|+n-1} e^{-\frac{\varepsilon \pi}{2} \sigma_x^2 r^2}  dr} = \\ { } \\

= \frac{2 \pi^{\frac{n}2}}{\Gamma\left( \frac{n}2 \right)} \frac{C_2}{\left( \varepsilon \pi \sigma_x^2 \right)^{M_2+n|A|+n-1}}
\,\,
\partial_r^{M_2+n|A|+n-2} e^{-\frac{\varepsilon \pi}{2} \sigma_x^2 r^2} |_{r = \varepsilon ^{-\frac{1}{4}+\frac{1}{2(M_2+1)}}} \leqslant \\  { } \\

\leqslant C'_2 \, \varepsilon^{-M_2-n|A|-n+1} 
\left| H_{M_2+n|A|+n-2} \left( \frac{ \pi \sigma_x^2 \varepsilon^{ \frac{3}{4}+\frac{1}{2(M_2+1)} } }{2} \right) \right| 
e^{-\frac{\varepsilon \pi}{2} \sigma_x^2 \varepsilon^{\frac{1}{2} +\frac{1}{M_2+1}}},
\end{array}
\end{equation}
where of course $H_s(x)$ is the Hermite polynomial of order $s$. By strengthening the assumption $M_2 \leqslant -3$ (which already has appeared in equation (\ref{eq172}) ) to 
\begin{equation}
\label{eq172o}
M_2<-3
\end{equation}
we get that $\varepsilon^{ \frac{3}{4}+\frac{1}{2(M_2+1)} } =o(1)$, and therefore only the zero order term of $H_{M_2+n|A|+n-2}$ has to be considered in the last line of equation (\ref{est2}). Observe moreover that
\[
e^{-\frac{\varepsilon \pi}{2} \sigma_x^2 \varepsilon^{\frac{1}{2} +\frac{1}{M_2+1}}} \leqslant 1.
\]
Using these observations, equation (\ref{est2}) implies
\begin{equation}
\begin{array}{c}
\label{est3}
|I_3(A)(x)| \leqslant
C''_2 \, \varepsilon^{-M_2-n|A|-n+1} \leqslant
C''_2 \, \varepsilon^{-M_2-nN-n+1}
\end{array}
\end{equation}
(remember that $|A| \leqslant N $). Now asking that equation (\ref{est1}) holds is equivalent to asking
\begin{equation}
\begin{array}{c}
\label{est4}
\varepsilon^{-M_2-nN-n+1}  \leqslant \varepsilon^{\frac{N+1}{2}} \,\, \Leftrightarrow \,\,
-M_2-nN-n+1 \geqslant \frac{N+1}{2}  \,\, \Leftrightarrow \\ { } \\
 \Leftrightarrow \,\, M_2 \leqslant -N \left( n+\frac{1}2 \right)-n+\frac{1}2.
\end{array}
\end{equation}

The proof is complete. \\ \vspace{0.5mm}

\noindent {\bf Remarks on the choice of the parameters: }

The constraints for $M_2$ come up during the proof in equations (\ref{eq172}), (\ref{eq164}), (\ref{eq166}), (\ref{eq170}), (\ref{eq172o}), (\ref{est4}). Clearly, each of these constraints can be satisfied for  $M_2$ small enough, depending on $N,n$.

The constraint for $\sigma_x^2$ came up right after equation (\ref{eq160}).

The constraint for $M_1$ appeared right after equation (\ref{eq163}). \\ \vspace{0.5mm} 

Theorem \ref{thrm15} result must be compared with its counterpart for the WT, quoted here (in adapted notation) from \cite{GMMP}:

\begin{theorem}[Semiclassical finite-order approximations to the Wigner calculus]
\label{thrm20}
Let $p(x,k) \in C^\infty(\mathbb{R}^{2n})$ satisfy, for some $M'_1>0$
\begin{equation}
\forall \,\,\, a \in \mathbb{N}^{2n} \, |\partial_{x_1...x_nk_1...k_n}^a p(x,k)| \leqslant \, C_a (1+|k|)^{M'_1}.
\end{equation}
Assume also $||f^\varepsilon ||_{L^2(\mathbb{R}^n)},  ||g^\varepsilon ||_{L^2(\mathbb{R}^n)} \leqslant M_0 \,\,\, \forall \varepsilon >0$. Then
\begin{equation}
\begin{array}{l}
W^\varepsilon[p(x,\varepsilon \partial_x) f^\varepsilon,g^\varepsilon](x,k) = 
p(x,k) W^\varepsilon[ f^\varepsilon,g^\varepsilon](x,k) +\\ { } \\
+\frac{\varepsilon}{4 \pi i} \sum\limits_{d=1}^n \left[{ \partial_{k_d} p(x,k) \partial_{x_d} W^\varepsilon[ f^\varepsilon,g^\varepsilon](x,k) - \partial_{x_d} p(x,k) \partial_{k_d} W^\varepsilon[ f^\varepsilon,g^\varepsilon](x,k) }\right] + \varepsilon^2 r_\varepsilon ,
\end{array}
\end{equation}
where $r_\varepsilon$ is bounded in $\mathcal{S}' (\mathbb{R}^{2n})$ as $\varepsilon \shortrightarrow 0$.
\end{theorem}

\textbf{Acknowledgements.} 
AGA was supported by the INRIA via the ERCIM
``Alain Bensoussan fellowship"; he would also like to thank the Wolfgang Pauli Institut for its hospitality during research visits in the course of this work.
TP would like to thank the Wolfgang Pauli Institut (UMI of CNRS) for its hospitality at the beginning of this project.
Support by the Austrian Ministry of Science via its grant for the
WPI, by the Viennese Fund for Technology and Science (WWTF), and by the EU funded Marie Curie project DEASE (contract MEST-CT-2005-021122) is acknowledged.


\end{document}